\documentclass[12pt, twoside]{article}
\usepackage{times}
\usepackage{graphics}
\usepackage{theorem}
\usepackage{graphicx}
\usepackage{amsmath}
\usepackage{latexsym}
\usepackage{amssymb}
\usepackage{titlesec}
\usepackage[flushmargin]{footmisc}

\theorembodyfont{\itshape}
\newtheorem{thm}{Theorem}[section]
\newtheorem{lem}{Lemma}[section]
\newtheorem{prop}{Proposition}[section]
\newtheorem{coro}{Corollary}[section]

\theorembodyfont{\rmfamily}
\newtheorem{defn}{Definition}[section]{\bf}{\rm}
\newtheorem{assumpt}{Assumption}[section]{\bf}{\rm}

\newtheorem{rem}{Remark}[section]{\itshape}{\rmfamily}
\newenvironment{proof}{\noindent{\it Proof.~~}}{\qed}

\makeatletter
\def\eqnarray{\stepcounter{equation}\let\@currentlabel=\theequation
\global\@eqnswtrue
\global\@eqcnt\z@\tabskip\@centering\let\\=\@eqncr
$$\halign to \displaywidth\bgroup\@eqnsel\hskip\@centering
  $\displaystyle\tabskip\z@{##}$&\global\@eqcnt\@ne 
  \hfil$\;{##}\;$\hfil
  &\global\@eqcnt\tw@ $\displaystyle\tabskip\z@{##}$\hfil 
   \tabskip\@centering&\llap{##}\tabskip\z@\cr}
\makeatother
\makeatletter
    \renewcommand{\theequation}{%
    \thesection.\arabic{equation}}
    \@addtoreset{equation}{section}
  \makeatother

\def\vc#1{\mbox{\boldmath $#1$}}

\newcommand{\down}[2]{\smash{\lower#2\hbox{#1}}}
\newcommand{\up}[2]{\smash{\lower-#2\hbox{#1}}}
\newcommand{\dm}{\displaystyle}

\newcommand{\qed}{\hspace*{\fill}$\Box$\rule[-10pt]{0pt}{10pt}}

\newcommand{\EE}{\mathsf{E}}
\newcommand{\PP}{\mathsf{P}}
\newcommand{\calC}{\mathcal{C}}
\newcommand{\calD}{\mathcal{D}}
\newcommand{\calL}{\mathcal{L}}
\newcommand{\calR}{\mathcal{R}}
\newcommand{\calS}{\mathcal{S}}
\newcommand{\SC}{\mathcal{SC}}
\newcommand{\bbD}{\mathbb{D}}

\newcommand{\bbR}{\mathbb{R}}
\newcommand{\bbZ}{\mathbb{Z}}

\newcommand{\rd}{{\rm d}}
\newcommand{\re}{{\rm e}}
\newcommand{\rmi}{{\rm i}}
\newcommand{\on}{{\rm on}}
\newcommand{\off}{{\rm off}}
\def\simhm#1{\stackrel{#1}{\sim}}
\def\lesimhm#1{\lesssim_{#1}}
\def\gesimhm#1{\gtrsim_{#1}}

\renewcommand{\labelenumi}{(\roman{enumi})}
\def\simhm#1{\stackrel{#1}{\sim}}
\def\lesimhm#1{\lesssim_{#1}}
\def\gesimhm#1{\gtrsim_{#1}}
\newcommand{\dd}[1]{\if#11 1\!\!1 
\else {\if#1C I\!\!\!C
\else {\if#1G I\!\!\!G 
\else {\if#1J J\!\!\!J 
\else {\if#1S S\!\!\!S
\else {\if#1Z Z\!\!\!Z
\else {\if#1Q O\!\!\!\!Q
\else I\!\!#1
\fi} 
\fi}
\fi}
\fi} 
\fi} 
\fi} 
\fi} 

\makeatother

\begin{document}\thispagestyle{plain} 

\hfill

{\Large{\bf
\begin{center}
Tail asymptotics for cumulative processes sampled at heavy-tailed
random times with applications to queueing models
in Markovian environments%
%
\end{center}
}
}

\begin{center}
{
Hiroyuki Masuyama%
\footnote[2]{E-mail: masuyama@sys.i.kyoto-u.ac.jp}
}

\medskip

{\small
Department of Systems
Science, Graduate School of Informatics, Kyoto University\\
Kyoto 606-8501, Japan
}

\bigskip
\medskip

{\small
\textbf{Abstract}

\medskip

\begin{tabular}{p{0.85\textwidth}}
This paper considers the tail asymptotics for a cumulative process
$\{B(t); t \ge 0\}$ sampled at a heavy-tailed random time $T$. The
main contribution of this paper is to establish several sufficient
conditions for the asymptotic equality $\PP(B(T) > bx) \sim \PP(M(T) >
bx) \sim \PP(T>x)$ as $x \to \infty$, where $M(t) = \sup_{0 \le u \le
  t}B(u)$ and $b$ is a certain positive constant. The main results of
this paper can be used to obtain the subexponential asymptotics for
various queueing models in Markovian environments. As an example,
using the main results, we derive subexponential asymptotic
formulas for the loss probability of a single-server finite-buffer
queue with an on/off arrival process in a Markovian environment.
\end{tabular}
}
\end{center}

\begin{center}
\begin{tabular}{p{0.90\textwidth}}
{\small
{\bf Keywords:} %
Queue; 
Markovian environment;
heavy-tailed; 
tail asymptotics;
sampling; 
cumulative process
%
%

\medskip

{\bf Mathematics Subject Classification:} %
60G50 (Sums of independent random variables; random walks) $\cdot$
60F10 (Large deviations) $\cdot$ 
60K25 (Queueing theory)
%
}
\end{tabular}

\end{center}

\section{Introduction}\label{introduction}

The main purpose of this paper is to provide mathematical tools for
obtaining the heavy-tailed asymptotic behavior of queueing models in
Markovian environments. Many researchers have studied the heavy-tailed
asymptotics of the random sum of random variables (r.v.s), and several
interesting results have been reported in the literature. However,
those results cannot be applied directly to queueing models in
Markovian environments, such as queues with batch Markovian arrival
processes (BMAPs) \cite{Luca91} and general semi-Markovian arrival
processes. Therefore in this paper, we construct a framework to
study the heavy-tailed asymptotics for such queueing models.

Let $\{B(t);t\ge0\}$ denote a (possibly delayed) cumulative process on
$\bbR := (-\infty,\infty)$, where $|B(0)| < \infty$ with probability
one (w.p.1) (see, e.g., \cite[Section 2.11]{Wolf89}). By definition,
there exist regenerative points $0 \le \tau_0 < \tau_1 < \tau_2 <
\cdots$ such that $\{B(t+\tau_n)-B(\tau_n); t \ge 0\}$ ($n=0,1,\dots$)
is stochastically equivalent to $\{B(t+\tau_0)-B(\tau_0); t \ge 0\}$
and is independent of $\{B(u);0 \le u < \tau_n\}$. Let
\begin{eqnarray}
\Delta B_n 
&=& \left\{
\begin{array}{ll}
B(\tau_0), & n = 0,
\\
B(\tau_n)-B(\tau_{n-1}), & n=1,2,\dots,
\end{array}
\right.
~~
\Delta \tau_n 
= \left\{
\begin{array}{ll}
\tau_0, & n = 0,
\\
\tau_n-\tau_{n-1}, & n=1,2,\dots,
\end{array}
\right.
\label{defn-Delta-tau_n}
\\
\Delta B_n^{\ast} 
&=& 
\left\{
\begin{array}{ll}
\dm \sup_{0 \le t \le \tau_0} \max(B(t), 0), & n = 0,
\\
\dm\sup_{\tau_{n-1} \le t \le \tau_n} B(t) - B(\tau_{n-1}), & n=1,2,\dots.
\end{array}
\right.
\nonumber
\end{eqnarray}
Clearly, $\Delta B_n^{\ast} \ge \Delta B_n$ for $n=0,1,\dots$. Further
$\{\Delta \tau_n;n=1,2,\dots\}$ (resp.\ $\{\Delta B_n;n=1,2,\dots\}$
and $\{\Delta B_n^{\ast};n=1,2,\dots\}$) is a sequence of independent
and identically distributed (i.i.d.)~r.v.s, which is independent of
$\Delta \tau_0$ (resp.\ $\Delta B_0$ and $\Delta B_0^{\ast}$).

Throughout this paper, we assume that
\begin{eqnarray}
&&
\PP(0\le \Delta \tau_n < \infty)   
= \PP(0 \le \Delta B_n^{\ast} < \infty) = 1~~(n=0,1),
\nonumber
\\
&&
\EE[|\Delta B_1|] <\infty,~~~
0 < \EE[\Delta \tau_1] < \infty,~~~
b := {\EE[\Delta B_1] \big/ \EE[\Delta \tau_1] } > 0.
\label{cond-SLLN-01}
\end{eqnarray}
Under these basic conditions, we study the heavy-tailed asymptotics of
$B(T)$, where $T$ is a nonnegative r.v.\ representing the sampling
time of $\{B(t)\}$.  More specifically, we establish sufficient
conditions for a simple asymptotic formula:
\begin{equation}
\PP(B(T) > bx) 
\simhm{x}
\PP(M(T) > bx) \simhm{x} \PP(T > x),
\label{Jele-asymp}
\end{equation}
where $M(t) = \sup_{0 \le u \le t}B(u)$ for $t\ge0$, and where for any
functions $f$ and $g$, $f(x) \simhm{x} g(x)$ represents
$\lim_{x\to\infty}f(x)/g(x) = 1$ (if the limit holds).

We now give a brief discussion of the conditions for
(\ref{Jele-asymp}) to hold. Note that if $\{B(t)\}$ has no deviation,
i.e., $B(t) = bt$ for all $t \ge 0$, then $\PP(B(T) > bx) = \PP(T >
x)$. In general, however, $B(t) - bt$ has a deviation from zero, which
is caused by the distributions of $\Delta B_n$ and $\Delta
\tau_n$~($n=0,1$).  Thus $\PP(B(T) > bx)$ may be decomposed in an
intuitive way:
\begin{eqnarray}
\PP(B(T) > bx)
&\approx& \PP(T > x) 
+ \big(\mbox{remainder term associated with $\Delta B_n$ and $\Delta \tau_n$)}.
\label{eqn-P(B(T)>bx)}
\end{eqnarray}
If the remainder term of (\ref{eqn-P(B(T)>bx)}) is negligible compared
with $\PP(T > x)$ as $x \to \infty$, then (\ref{Jele-asymp})
holds. Asmussen et~al.~\cite{AsmuKlupSigm99} show that if $T$ is
independent of $\{B(t)\}$, then an important necessary condition for
(\ref{Jele-asymp}) is that $\sqrt{T}$ is heavy-tailed, i.e., $\PP(T >
x) = e^{-o(\sqrt{x})}$, where for any functions $f$ and $g$, $f(x) =
o(g(x))$ represents $\lim_{x\to\infty}f(x)/g(x) = 0$ (if the limit
holds). On the other hand, if the remainder term of
(\ref{eqn-P(B(T)>bx)}) is not negligible, then it is likely that the
asymptotic behavior of $\PP(B(T) > bx)$ is complicated. Indeed,
Asmussen et al.~\cite{AsmuKlupSigm99} and Foss and
Korshunov~\cite{Foss00} consider such cases, and they present some
asymptotic formulas with implicit functions for two special cumulative
processes: the Poisson counting process \cite{AsmuKlupSigm99} and the
sum of nonnegative r.v.s~\cite{Foss00}. Although it is challenging to
generalize those results, we leave it for future work. In this paper,
we focus on the case where (\ref{Jele-asymp}) holds.

As mentioned at the beginning, this study is motivated by the
heavy-tailed asymptotics for queueing models in Markovian
environments. A typical example of the application of this study is as
follows. Consider a stationary BMAP/GI/1 queue. Suppose that $B(t)$ is
the total number of stationary BMAP arrivals in the interval $(0,t]$,
  which is a cumulative process. Further suppose that $T$ is the
  service time of one customer and is independent of $\{B(t)\}$. In
  this setting, $b$ is the arrival rate and $b\EE[T]$ is the traffic
  intensity. Note here (see, e.g., Proposition~3.1 in Masuyama et
  al.~\cite{Masu09}) that the subexponential asymptotics of the
  stationary queue length $L$ is connected to that of $B(T)$ as
  follows:
\[
\PP(L > x) \simhm{x} {1 \over 1 - b\EE[T]} \int_x^{\infty} \PP(B(T) > y) \rd y.
\]
Therefore, if the subexponential asymptotics of $\PP(B(T) > x)$ is
given, we can obtain an asymptotic formula for the stationary queue
length $L$. Especially, when (\ref{Jele-asymp}) holds, we have the
following simple and explicit formula:
\[
\PP(L > x) \simhm{x} {b \EE[T] \over 1 - b\EE[T]} \cdot \PP(T_{\rm e} > x/b),
\]
where  $T_{\re}$ denotes the equilibrium
r.v.\ of $T$, i.e., $\PP(T_{\re} \le x) = (1/\EE[T])\int_0^x\PP(T >
y)\rd y$ for $x \ge 0$.

Next we review related work. For this purpose, we introduce two
classes of distributions (for details, see
Appendix~\ref{higher-order}).
\begin{defn}\label{defn-subclass-L}
A nonnegative r.v.\ $X$ and its distribution function (d.f.) $F_X$
belong to the $p$th-order long-tailed class $\calL^p$ ($p \ge 1$) if
$X^{1/p} \in \calL$, i.e., $\PP(X^{1/p} > x) > 0$ for all $x \ge 0$
and $\PP(X^{1/p} > x+y) \simhm{x} \PP(X^{1/p} > x)$ for some (thus
all) $y > 0$. Further if $X \in \calL^{1/\theta}$ (resp.\ $F_X \in
\calL^{1/\theta}$) for any $0 < \theta \le 1$, we write $X \in
\calL^{\infty}$ (resp.\ $F_X \in \calL^{\infty}$) and call $X$
(resp.\ $F_X$) infinite-order long-tailed.
\end{defn}

\begin{defn}\label{defn-class-C}
A nonnegative r.v.\ $X$ and its d.f.\ $F_X$ belong to the consistent
variation class $\calC$ if $\overline{F}_X(x) > 0$ for all $x \ge 0$
and
\[
\lim_{v\downarrow1}\liminf_{x\to\infty}
{\overline{F}_X(vx) \over \overline{F}_X(x)}
= 1
~~
\mbox{or equivalently,}~~
\lim_{v\uparrow1}\limsup_{x\to\infty}
{\overline{F}_X(vx) \over \overline{F}_X(x)}
= 1,
\]
where $\overline{F}_X(x) = 1 - F_X(x)$ for all $x \in \bbR$. Note that
every distribution with a consistently varying tail is infinite-order
long-tailed (i.e., $\calC \subset \calL^{\infty}$; see
Lemma~\ref{lem-C-in-L^{infty}}).
\end{defn}

The related work is classified into two cases: (i) $T$ is independent
of $\{B(t)\}$; and (ii) $T$ may depend on $\{B(t)\}$. The former is
called {\it independent-sampling case}, and the latter is called {\it
  dependent-sampling case}. The dependent-sampling case includes a
case where $T$ is a stopping time with respect to $\{B(t)\}$.

To the best of our knowledge, there are a few results for the
dependent-sampling case. Robert and Segers~\cite{Robe08}
consider a special case where
\begin{equation}
B(t)=\sum_{n=1}^{\lfloor t \rfloor}X_n~~
\mbox{with the $X_n$'s being i.i.d.~nonnegative r.v.s}.
\label{cond-iid-sums}
\end{equation}
Note here that the summation over the empty set is defined as zero,
e.g., $\sum_{n=k}^l \cdot = 0$ for $k > l$. Thus if
(\ref{cond-iid-sums}) holds, then
\begin{align}
\Delta \tau_0 = 0,~~
\Delta B_0 = \Delta B_0^{\ast} = 0,~~
\Delta \tau_n=1,~~
\Delta B_n = \Delta B_n^{\ast} = X_n~(n=1,2,\dots).
\label{eqns-special-case}
\end{align}
For this special case, Robert and Segers~\cite{Robe08} present the
following:
\begin{prop}[Theorem 4.1 in Robert and Segers~\cite{Robe08}]
\label{prop-Robe08-thm4.1}
Suppose that $X,X_1,X_2,\dots$ are i.i.d.\ nonnegative r.v.s. Further
suppose that (i) $T$ satisfies
\begin{equation}
\lim_{x\to\infty}{\PP(T > x + ya(x)) \over \PP(T > x)} = e^{-y},
\qquad y \in \bbR,
\label{prop-robert-cond-02-b}
\end{equation}
for some function $a(x)$ ($x \ge 0$) such that
$x^{2/3} = o(a(x))$; and (ii) $\EE[e^{\gamma X}] <
\infty$ for some $\gamma > 0$. Under these conditions, we have
\begin{equation}
\PP(X_1+\cdots+X_{\lfloor T \rfloor} > \EE[X] x) \simhm{x} \PP(T > x).
\label{eqn-Robe08-depend}
\end{equation}
\end{prop}
\begin{prop}[Theorem 3.1 in Robert and Segers~\cite{Robe08}]
\label{prop-Robe08-thm3.1}
Suppose that $X,X_1,X_2,\dots$ are i.i.d.\ nonnegative r.v.s.  Further
suppose that (i) $T \in \calC$; (ii) $\EE[X^{\gamma}] < \infty$ for
some $\gamma > 1$; and (iii) $x\PP(X>x)=o(\PP(T>x))$. Under these
conditions, (\ref{eqn-Robe08-depend}) holds.
\end{prop}
Compared with Proposition~\ref{prop-Robe08-thm4.1},
Proposition~\ref{prop-Robe08-thm3.1} requires a heavier tail of $T$
but relaxes the condition on $X$, which is implied by
(\ref{eqn-P(B(T)>bx)}).

For the independent-sampling case, several results have been
reported. However, as far as we know, only Jelenkovi\'{c}
et~al.~\cite{Jele04} consider the general cumulative process
$\{B(t)\}$:
\begin{prop}[Proposition~3 in Jelenkovi\'{c} et~al.~\cite{Jele04}]\label{prop-Jele04}
Suppose that $T$ is independent of $\{B(t);t\ge0\}$. Further suppose
that (i) $T \in \calL^2$ (i.e., $\sqrt{T} \in \calL$); (ii)
$\EE[(\Delta \tau_1)^2] < \infty$ and $\Delta B_n \ge 0$ ($n=0,1$)
w.p.1; and (iii) $\EE[\exp\{\eta\sqrt{\Delta B_n^{\ast}}\}] < \infty$
($n=0,1$) for some $\eta > 0$. Under these conditions,
(\ref{Jele-asymp}) holds.
\end{prop}

According to Jelenkovi\'{c} et~al.~\cite{Jele04}'s result, the
condition on $\Delta B_n^{\ast}$ and thus $\Delta B_n$ is insensitive
to the tail of $T$, given that $T \in \calL^2$. On the other hand,
(\ref{eqn-P(B(T)>bx)}) implies that the conditions on $\Delta B_n$ and
$\Delta \tau_n$ for (\ref{Jele-asymp}) to hold are weaker as the tail
of $T$ is heavier.  In fact, as with the dependent-sampling case, such
a result has been reported by Ale\v{s}kevi\v{c}ien\.{e}
et~al.~\cite{Ales08}.
\begin{prop}[Theorem 1.2 in Ale\v{s}kevi\v{c}ien\.{e} et~al.~\cite{Ales08}]
\label{prop-Ales08-thm2.1}
Suppose that $X,X_1,X_2,\dots$ are i.i.d.\ nonnegative r.v.s and $T$
is independent of $\{X_n;n=1,2,\dots\}$. Further suppose that (i) $T
\in \calC$; (ii) $\EE[X] < \infty$; and (iii) $\EE[T] < \infty$ and
$\PP(X > x) = o(\PP(T > x))$. Under these conditions,
(\ref{eqn-Robe08-depend}) holds.
\end{prop}
Note here that Proposition~\ref{prop-Jele04} does not allow that the
tail distribution of $X$ is heavier than $e^{-\eta \sqrt{x}}$; whereas
Proposition~\ref{prop-Ales08-thm2.1} does.

Lin and Shen~\cite{Lin11} extend Proposition~\ref{prop-Ales08-thm2.1}
to the case where the $X_n$'s are asymptotically quadrant
sub-independent and identically distributed (see Theorem~2.1~(I)
therein).  Robert and Segers~\cite{Robe08} present a theorem result
similar to Proposition~\ref{prop-Ales08-thm2.1} (see Theorem 3.2
therein). The theorem states that (\ref{eqn-Robe08-depend}) requires
$\EE[X^r] < \infty$ for some $r > 1$, which is more restrictive than
condition (ii) of Proposition~\ref{prop-Ales08-thm2.1}. However, the
theorem also presents a sufficient condition for
(\ref{eqn-Robe08-depend}) with $\EE[T] = \infty$, which is described
in the following:
\begin{prop}[Theorem 3.2 in Robert and Segers~\cite{Robe08}]
\label{prop-Robe08-thm3.2}
Suppose that $X,X_1,X_2,\dots$ are i.i.d.\ nonnegative r.v.s and $T$
is independent of $\{X_n;n=1,2,\dots\}$. Further suppose that (i) $T
\in \calC$ and $\EE[T] = \infty$; (ii) $\EE[X^r] < \infty$ for some $r
> 1$; and (iii) for some $1 \le q < r$,
\[
\limsup_{x\to\infty}{\EE[T\cdot \dd{1}(T \le x)] \over x^q \PP(T > x)} < \infty,
\]
where $\dd{1}(\chi)$ denotes the indicator
function of event (or condition) $\chi$. Under these conditions,
(\ref{eqn-Robe08-depend}) holds.
\end{prop}

In what follows, we summarize the contributions of this paper.  For
the dependent-sampling case, we assume that $\{B(t)\}$ is
nondecreasing with $t$ (e.g., $\{B(t)\}$ is the counting process of
BMAP arrivals).  Under this assumption, we present two theorems:
Theorems~\ref{thm-dependent-01} and \ref{thm-dependent-02}, which are
extensions of Propositions~\ref{prop-Robe08-thm4.1} and
\ref{prop-Robe08-thm3.1}, respectively, to the general cumulative
process. In addition, the two theorems are still more general than the
corresponding propositions even if (\ref{cond-iid-sums}) holds, i.e.,
$B(T)$ is reduced to the random sum of i.i.d.\ nonnegative r.v.s.

As for the independent-sampling case, we do not necessarily assume
that $\{B(t)\}$ is nondecreasing with $t$, which means that $\Delta
B_n$ can take negative values. We first present two theorems:
Theorems~\ref{thm-independ-01} and
\ref{thm-independ-02}. Theorem~\ref{thm-independ-01} provides a weaker
sufficient condition for (\ref{Jele-asymp}) than that in
Proposition~\ref{prop-Jele04}. Theorem~\ref{thm-independ-02} is an
extension of Propositions~\ref{prop-Ales08-thm2.1} and
\ref{prop-Robe08-thm3.2} to the general cumulative process. However,
unfortunately, when $\{B(t)\}$ satisfies (\ref{cond-iid-sums}), one of
the conditions of Theorem~\ref{thm-independ-02} is more restrictive
than the corresponding ones of Propositions~\ref{prop-Ales08-thm2.1}
and \ref{prop-Robe08-thm3.2}. Thus, instead of the general cumulative
process, we next consider a special case where $B(t) = B(\lfloor t
\rfloor)$ for all $t \ge 0$ and $\{B(n);n=0,1,\dots\}$ is the additive
component of a discrete-time Markov additive process (see, e.g.,
\cite[Chapter~XI, Section~2]{Asmu03}), which implies that $B(T)$ is
the random sum of r.v.s with Markovian correlation. Under this
assumption, we prove Theorems~\ref{thm-MAdP-01} and \ref{thm-MAdP-02},
which completely include Propositions~\ref{prop-Ales08-thm2.1} and
\ref{prop-Robe08-thm3.2} as special cases. Further the two theorems
are readily extended to the case where $\{B(t)\}$ is the additive
component of a continuous-time Markov additive process.

As mentioned above, our results for the independent-sampling case are
more general than those in the literature and thus can be applied to
derive new asymptotic formulas for queueing models in Markovian
environments. Indeed, Masuyama~\cite{Masu13a} derives some new
subexponential asymptotic formulas for the BMAP/GI/1 queue by using
the results of this paper.  Masuyama \cite{Masu13b} also presents
subexponential asymptotic formulas for the BMAP/GI/1 queue with
retrials by combining the results of \cite{Masu13a} with the
subexponential tail equivalence of the queue length distributions of
BMAP/GI/1 queues with and without retrials.  In addition, unlike the
previous studies, our results for the independent-sampling case can be
applied to queues with negative customers (see, e.g., \cite{Baye96})
because the results do not necessarily require the monotonicity of
$\{B(t)\}$.

To demonstrate the utility of our results for the dependent-sampling
case as well as the independent-sampling case, we discuss their
application to the subexponential asymptotics of the loss probability
of a discrete-time single-server queue with a finite buffer fed by an
on/off arrival process in a Markovian environment. In the on/off
arrival process, the lengths of on-periods (resp.\ off-periods) are
i.i.d.\ with a general distribution, and arrivals in each on-period
follow a discrete-time BMAP started with some initial distribution at
the beginning of the on-period. We call the arrival process {\it
  on/off batch Markovian arrival process (ON/OFF-BMAP)}, which is a
generalization of the batch-on/off process \cite{Galm01} and is
closely related to a platoon arrival process (PAP)
\cite{Alfa95,Breu05} (see also Remarks~\ref{rem-ON/OFF-BMAP-01} and
\ref{rem-PAP}). For analytical convenience, we assume that service
times are all equal to the unit of time. The queueing model is denoted
by (ON/OFF-BMAP)/D/1/$K$ in Kendall's notation. For this queue, we
derive subexponential asymptotic formulas for the loss probability by
combining our results with the existing one on a finite GI/GI/1 queue
\cite{Jele99}.

The rest of this paper is organized as
follows. Section~\ref{sec-preliminary} introduce some
definitions. Section~\ref{sec-main-result} presents the main results
of this paper, and Section~\ref{sec-application-queue} discusses their
application to the (ON/OFF-BMAP)/D/1/$K$
queue. Appendix~\ref{sec-preliminary-lem} is devoted to technical
lemmas. The proofs of all the lemmas and the main results are given in
Appendices~\ref{sec-proof-lemmas} and \ref{sec-proof-main-results}.

\section{Basic Definitions}\label{sec-preliminary}

In this section, we provide the definitions of the subexponential
distribution and some related classes of distributions.  For later
use, we first introduce the following notations. Let $C$ (resp.~$c$)
denote a special symbol representing a sufficiently large
(resp.~small) positive constant, which takes an appropriate value
according to the context. Thus $C$ (resp.~$c$) can take different
values in different places. For example, $C$ in a place may be equal
to $C+1$, $2C$ and $C^2$, etc. in other places.  For any $x \in \bbR$,
let $x^+ = \max(x,0)$. For any r.v.\ $U$ in $\bbR$, let $F_U$ denote
the d.f.\ of $U$, i.e., $F_U(x) = \PP(U \le x)$ for $x \in
\mathbb{R}$, which is assumed to be right-continuous. Further let
$\overline{F}_U = 1 - F_U$ and $Q_U = - \log \overline{F}_U$. The
latter is called the cumulative hazard function of $U$. Finally, for
any nonnegative functions $f$ and $g$, $f(x) = O(g(x))$, $f(x)
\lesimhm{x} g(x)$ and $f(x) \gesimhm{x} g(x)$ represent
\[
\limsup_{x\to\infty}f(x)/g(x)
< \infty,
\qquad
\limsup_{x\to\infty}f(x)/g(x)
\le 1,\qquad \liminf_{x\to\infty}f(x)/g(x) \ge 1,
\]
respectively.

\subsection{Subexponential distributions}

We begin with the definition of the subexponential class.
\begin{defn}\label{defn-class-S}
A nonnegative r.v.\ $X$ and its d.f.\ $F_X$ belong to the
subexponential class $\calS$ if $\PP(X>x) > 0$ for all $x \ge 0$ and
$\PP(X_1+X_2 > x) \simhm{x} 2\PP(X>x)$, where $X_1$ and $X_2$ are
independent copies of $X$.
\end{defn}

\begin{rem}
The class $\calS$ was first introduced by Chistyakov~\cite{Chis64},
and it was shown that $\calS$ is a strictly subclass of class $\calL$,
i.e., $\calS \subset \calL$ (see \cite{Pitm80}). \
\end{rem}

Next we introduce two subclasses of $\calS$. The first one is class
$\calS^{\ast}$, which is a well-known subclass of $\calS$.
\begin{defn}
A nonnegative r.v.\ $X$ and its d.f.\ $F_X$ belong to class
$\calS^{\ast}$ if $\EE[X] < \infty$ and
\[
\lim_{x\to\infty}\int_0^x {\overline{F}_X(x-y) \over \overline{F}_X(x)}
\overline{F}_X(y)\rd y = 2\EE[X].
\]
\end{defn}

\begin{rem}\label{rem-S^*}
An important property of $\calS^{\ast}$ is that $F \in \calS^{\ast}$
implies $F,F_{\re} \in \calS$, where $F_{\re}$ denotes the equilibrium
distribution (or integrated tail distribution) of $F$, i.e.,
$F_{\re}(x) =\int_0^x \overline{F}(y) \rd y/\int_0^{\infty}
\overline{F}(y) \rd y$ for $x \ge 0$ (see \cite[Theorem 3.2]{Klup88}).
\end{rem}

The second one is the subexponential concave class $\SC$, which is a
subclass of $\calS^{\ast}$, i.e., $\SC \subset \calS^{\ast}$ (see
\cite[Lemma 1]{Shne04}). The class $\SC$ plays a key role in
establishing large deviation bounds for a cumulative process. The
definition of $\SC$ is as follows:
\begin{defn}\label{defn-SC'}
A nonnegative r.v.\ $X$ and its d.f.\ $F_X$ and cumulative hazard
function $Q_X$ belong to the subexponential concave class $\SC$ if the
following are satisfied: (i) $Q_X$ is eventually concave; (ii)
$\lim_{x \to \infty}Q_X(x)/\log x = \infty$; and (iii) there exist
some $0 < \alpha < 1$ and $x_0 > 0$ such that $Q_X(x)/x^{\alpha}$ is
nonincreasing for all $x \ge x_0$, i.e.,
\begin{equation}
{Q_X(x) \over Q_X(u)}  \le \left( {x \over u} \right)^{\alpha},
\qquad x \ge u \ge x_0.
\label{ineqn-Q_X-03}
\end{equation}
We may use the notation $\SC_{\alpha}$ to emphasize the parameter
$\alpha$.
\end{defn}

\begin{rem}
Typical examples of the cumulative hazard function in $\SC$ are (i)
$(\log x)^{\gamma}x^{\alpha}$ and (ii) $(\log x)^{\beta}$ for
sufficiently large $x$, where $0 < \alpha < 1$, $\beta > 1$ and
$\gamma \in \bbR$.
\end{rem}

\begin{rem}\label{remark-SC'-02}
If a nonnegative r.v.\ $X$ satisfies $\EE[e^{Q(X)}] < \infty$ for some
cumulative hazard function $Q \in \SC$, then $\EE[X^p] < \infty$ for
any $p \ge 0$ because $e^{Q(x)} \ge x^p$ for sufficiently large $x >
0$ (see condition (ii) of Definition~\ref{defn-SC'}).
\end{rem}

Appendix~\ref{appendix-SC} provides some lemmas and further remarks on
$\SC$.

\subsection{Dominatedly varying distributions}\label{subsec-class-D}

The definition of the dominated variation class is as follows:
\begin{defn}\label{defn-class-D}
A nonnegative r.v.\ $X$ and its d.f.\ $F_X$
belong to the dominated variation class $\calD$ if $\overline{F}_X(x)
> 0$ for all $x \ge 0$ and
\[
\limsup_{x\to\infty}{\overline{F}_X(vx) \over \overline{F}_X(x)} < \infty,
\]
for some (thus for all) $v \in (0,1)$.
\end{defn}

\begin{rem}\label{rem-class-D}
$\calC \subset \calL \cap \calD \subset \calS^{\ast} \subset \calS$
  (see \cite[Theorem 3.2]{Klup88} and \cite{Clin94,Embr84}).
\end{rem}

For any d.f.\ $F$, let
\[
\overline{F}{}_{\ast}(v)
= \liminf_{x\to\infty}{\overline{F}(vx) \over \overline{F}(x)},
\quad
\overline{F}{}^{\ast}(v)
= \limsup_{x\to\infty}{\overline{F}(vx) \over \overline{F}(x)},
\quad v > 0,
\]
and let
\[
r_+(F)
= -\lim_{v\to\infty}{\log\overline{F}{}_{\ast}(v) \over \log v},
\quad
r_-(F)
= -\lim_{v\to\infty}{\log\overline{F}{}^{\ast}(v) \over \log v}.
\]
Strictly, $r_+(F)$ and $r_-(F)$ are called the upper and lower
Matuszewska indices of the function $1/\overline{F}(x)$ on
$[0,\infty)$ (see, e.g., Section~2.1 in \cite{Bing89}). For
  simplicity, however, they are sometimes called the upper and lower
  Matuszewska indices of d.f.\ $F$.

\begin{prop}[Proposition~2.2.1 in \cite{Bing89}]\label{prop-Bing89}
If $F \in \calD$, then for any $\alpha_1 < r_-(F) $ and $\alpha_2 >
r_+(F)$ there exist positive numbers $x_i > 0$ and $C_i > 0$ ($i=1,2$)
such that
\begin{align*}
&&
{\overline{F}(x) \over \overline{F}(y)} 
&\le C_1 \left( {x \over y} \right)^{-\alpha_1},
& \forall x &\ge \forall y \ge x_1,
&&
\\
&&
{\overline{F}(x) \over \overline{F}(y)} 
&\ge C_2 \left( {x \over y} \right)^{-\alpha_2},
& \forall x &\ge \forall y \ge x_2.
&&
\end{align*}
The second inequality implies that $x^{-\alpha} = o(\overline{F}(x))$
for all $\alpha > r_+(F)$.
\end{prop}

\section{Main Results}\label{sec-main-result}

This section consists of three subsections. In
subsection~\ref{subsec-cumulative}, we present four sets of conditions
under which (\ref{Jele-asymp}) holds for the general cumulative
process. Unfortunately, the last set of conditions is not completely
weaker than the corresponding ones in the literature if $\{B(t)\}$
satisfies (\ref{cond-iid-sums}), i.e., $B(T)$ is reduced to the random
sum of nonnegative r.v.s. Thus in subsection~\ref{subsec-MAdP}, we
discuss a special case where $B(t) = B(\lfloor t \rfloor)$ for all $t
\ge 0$ and $\{B(n);n=0,1,\dots\}$ is the additive component of a
discrete-time Markov additive process. For the special case, we have
two sets of conditions, which are weaker than the known ones even if
$\{B(t)\}$ satisfies (\ref{cond-iid-sums}). Finally in
subsection~\ref{subsec-MAdP-02}, we extend the results presented in
subsection~\ref{subsec-MAdP} to a continuous-time Markov additive
process.

\subsection{General case}\label{subsec-cumulative}

In this subsection, we assume $b=1$, i.e., $\EE[\Delta B_1] =
\EE[\Delta \tau_1]$ without loss of generality. Indeed, $\{
B(t)/b;t\ge0\}$ is a cumulative process with the same regenerative
points as those of $\{B(t)\}$, and the asymptotic equality
(\ref{Jele-asymp}) is rewritten as
\[
\PP(B(T)/b > x)
\simhm{x} \PP(M(T)/b > x) \simhm{x} \PP(T>x).
\]

In what follows, we first consider the dependent-sampling case
and then the independent-sampling case.

\subsubsection{Dependent-sampling case}\label{subsec-dependent}

In the dependent-sampling case, we assume that $\{B(t);t\ge0\}$ is
nondecreasing with $t$. In this case, $M(t) = B(t)$ for all $t \ge 0$
and thus (\ref{Jele-asymp}) is reduced to
\[
\PP(B(T) > bx) \simhm{x} \PP(T>x).
\]

\begin{thm}\label{thm-dependent-01}
Suppose that $\{B(t);t\ge0\}$ is nondecreasing with $t$.  Further
suppose that (i) $T \in \calL^{1/\theta}$ for some $0 < \theta \le
1/3$; and (ii) $\EE[\exp\{Q( (-B(0))^+ + \Delta \tau_0)\}] < \infty$,
$\EE[\exp\{Q(\Delta \tau_1)\}] < \infty$, $\EE[\exp\{Q( (\Delta B_0)^+
  )\}] < \infty$ and $\EE[\exp\{Q(\Delta B_1)\}] < \infty$ ($n=0,1$)
for some $Q \in \SC$ such that
\begin{equation}
x^{3\theta/2} = O( Q(x) ).
\label{add-eqn-28}
\end{equation}
Under these conditions, $\PP(B(T) > x) \simhm{x} \PP(T > x)$.
\end{thm}

\proof See Appendix~\ref{proof-thm-dependent-01}. \qed

\begin{rem}\label{rem-thm-dependent-01}
We prove Theorem~\ref{thm-dependent-01} by using
Lemma~\ref{prop-LD-B(t)}~(i) and (ii), which require condition (ii)
(see Remark~\ref{add-remark2}). In addition to the nondecreasingness
of $\{B(t)\}$, we assume $B(0) \ge 0$. It then follows that $(-B(0))^+
= 0$ and $(\Delta B_0)^+ = \Delta B_0$. Therefore condition (ii) is
reduced to $\EE[\exp\{Q(\Delta \tau_n)\}] < \infty$ and
$\EE[\exp\{Q(\Delta B_n)\}] < \infty$ ($n=0,1$) for some $Q \in \SC$
such that $x^{3\theta/2} = O( Q(x) )$.
\end{rem}

Theorem~\ref{thm-dependent-01} is a generalization of
Proposition~\ref{prop-Robe08-thm4.1}. To compare the two results, we
suppose that $\{B(t)\}$ satisfies (\ref{cond-iid-sums}). We then have
(\ref{eqns-special-case}) and $B(0) = 0$. Therefore conditions (i) and
(ii) of Theorem~\ref{thm-dependent-01} are reduced to the following
(see Remark~\ref{rem-thm-dependent-01}):

\begin{narrow}
\begin{enumerate}
\renewcommand{\labelenumi}{(\Roman{enumi})}
\item $T \in \calL^{1/\theta}$ for some $0 < \theta \le 1/3$; and
\item $\EE[\exp\{Q(X)\}] < \infty$ for some $Q \in \SC$
  satisfying (\ref{add-eqn-28}).
\end{enumerate}
\end{narrow}

\smallskip

Condition (I) is equivalent to $T \in \calL^3$ (see
Lemma~\ref{lem-01}~(ii)). On the other hand, condition (i) of
Proposition~\ref{prop-Robe08-thm4.1} implies that $T$ belongs to the
maximum domain of attraction of the Gumbel distribution (see, e.g.,
Theorem~3.3.27 in \cite{Embr97}). It further follows from
(\ref{prop-robert-cond-02-b}) and $x^{2/3}=o(a(x))$ that
\[
1 \ge \lim_{x\to\infty}{\PP(T > x + x^{2/3}) \over \PP(T > x)}
\ge \lim_{x\to\infty}{\PP(T > x + \varepsilon a(x)) \over \PP(T > x)}
= e^{-\varepsilon} \to 1 
\qquad \mbox{as } \varepsilon \to 0,
\]
which shows that $T \in \calL^3$. Thus condition (I) is weaker than
condition (i) of Proposition~\ref{prop-Robe08-thm4.1}. In addition,
condition (II) is satisfied by condition (ii) of
Proposition~\ref{prop-Robe08-thm4.1} due to $Q(x) = o(x)$ (see
Definition~\ref{defn-SC'}). As a result, the conditions of
Theorem~\ref{thm-dependent-01} are weaker than those of
Proposition~\ref{prop-Robe08-thm4.1}.

\begin{thm}\label{thm-dependent-02}
Suppose that $\{B(t);t\ge0\}$ is nondecreasing with $t$. Further
suppose that (i) $T \in \calC$; (ii) $\EE[(\Delta \tau_1)^2] <
\infty$; (iii) $\PP(-B(0) > x) = o(\PP(T > x))$, $\PP(\Delta \tau_n >
x) = o(\PP(T > x))$ and $\PP(\Delta B_n > x) = o(\PP(T > x))$
($n=0,1$); (iv) $x\PP(|\Delta B_1 - \Delta \tau_1| > x) = o(\PP(T >
x))$; and (v) either of the following is satisfied:

\begin{narrow}
\begin{enumerate}
\renewcommand{\labelenumi}{(\alph{enumi})}
\item $\EE[|\Delta B_1 - \Delta \tau_1|^r] < \infty$ for some $r > 1$;
  or
\item $\int_y^{\infty}x^{-1}\PP(T>x)\rd x < \infty$ for some $y \in
  (0,\infty)$.
\end{enumerate}
\end{narrow}
Under these conditions, $\PP(B(T) > x) \simhm{x} \PP(T > x)$.

\end{thm}

\proof See Appendix~\ref{proof-thm-dependent-02}. \qed

\begin{rem}\label{rem-thm-dependent-02}
The asymptotic upper bound $\PP(B(T)>x) \lesimhm{x} \PP(T>x)$ is
proved under the conditions that (iii$'$) $\PP(\Delta B_n > x) =
o(\PP(T > x))$ ($n=0,1$),~ (iv$'$) $x\PP(\Delta B_1 - \Delta \tau_1 >
x) = o(\PP(T > x))$ and (v$'$) either of the following holds:

\begin{narrow}
\begin{enumerate}
\renewcommand{\labelenumi}{(\alph{enumi})}
\item $\EE[\{(\Delta B_1 - \Delta \tau_1)^+\}^r] < \infty$ for some $r
  > 1$ or
\item $\int_y^{\infty}x^{-1}\PP(T>x)\rd x < \infty$ for some $y \in
  (0,\infty)$;
\end{enumerate}
\end{narrow}
whereas the asymptotic lower bound $\PP(B(T)>x) \gesimhm{x} \PP(T>x)$
is proved under the conditions that (iii$''$) $\PP(-B(0) > x) =
o(\PP(T > x))$ and $\PP(\Delta \tau_n > x) = o(\PP(T > x))$
($n=0,1$),~ (iv$''$) $x\PP(\Delta \tau_1 - \Delta B_1 > x) = o(\PP(T >
x))$ and (v$''$) either of the following holds:

\begin{narrow}
\begin{enumerate}
\renewcommand{\labelenumi}{(\alph{enumi})}
\item $\EE[\{(\Delta \tau_1 - \Delta B_1)^+\}^r] < \infty$ for some $r
  > 1$ or
\item $\int_y^{\infty}x^{-1}\PP(T>x)\rd x < \infty$ for some $y \in
  (0,\infty)$.
\end{enumerate}
\end{narrow}
These conditions are integrated into conditions (iii), (iv) and
(v). Further in the proof of the two asymptotic bounds, we use
Lemmas~\ref{add-prop-N(x)}, which requires condition (ii).
\end{rem}

Theorem~\ref{thm-dependent-02} is a generalization of
Proposition~\ref{prop-Robe08-thm3.1}. We compare them, assuming that
$\{B(t)\}$ satisfies (\ref{cond-iid-sums}). Under this assumption,
conditions (i)--(v) of Theorem~\ref{thm-dependent-02} are reduced to
the following:

\begin{narrow}
\begin{enumerate}
\renewcommand{\labelenumi}{(\Roman{enumi})}
\item $T \in \calC$; 
\item $x\PP(X > x) = o(\PP(T > x))$; and 
\item either of the following is satisfied:

\begin{enumerate}
\renewcommand{\labelenumii}{(\Alph{enumii})}
\item $\EE[X^r] < \infty$ for some $r > 1$; or 
\item $\int_y^{\infty}x^{-1}\PP(T > x) \rd x < \infty$ for some $y \in
(0,\infty)$.
\end{enumerate}
\end{enumerate}
\end{narrow}

Clearly, the set of conditions (I), (II) and (III.A) is the same as
that of conditions (i), (ii) and (iii) of
Proposition~\ref{prop-Robe08-thm3.1}.  Further the set of conditions
(I), (II) and (III.B) does not imply that of conditions (I), (II) and
(III.A).  In fact, suppose that $\PP(T>x) \simhm{x} (\log x)^{-2}$ and
$\PP(X>x) \simhm{x} x^{-1}(\log x)^{-3}$. We then have $T \in \calC$
and $x\PP(X > x) = o(\PP(T > x))$. It also holds that $\EE[X] <
\infty$ and $\int_y^{\infty}x^{-1}\PP(T > x) \rd x < \infty$ for some
$y \in (0,\infty)$, which follow from
\[
\int_y^{\infty} {\rd x \over x(\log x)^m}= {1 \over (m-1)(\log y)^{m-1}},
\quad y > 1,~m \neq 1.
\]
Thus conditions (I), (II) and (III.B) are satisfied. However,
condition (III.A) does not hold, i.e., $\EE[X^r] = \infty$ for any
$r>1$ because
\[
\PP(X^r > x) 
\simhm{x} {r^3 \over x^{1/r}(\log x)^3} 
\gesimhm{x} {r^3 \over x^{(1/r) + (r-1)/(2r)}}
= {r^3 \over x^{(r+1)/(2r)}},
\]
where $0 < (r+1)/(2r) < 1$ for $r > 1$. 

Consequently, the conditions of Theorem~\ref{thm-dependent-02} are
still weaker than those of Proposition~\ref{prop-Robe08-thm3.1} in the
context of the random sum of nonnegative i.i.d.\ r.v.s.

\subsubsection{Independent-sampling case}

\begin{thm}\label{thm-independ-01}
Suppose that $T$ is independent of $\{B(t);t\ge0\}$. Further suppose
that (i) $T \in \calL^{1/\theta}$ for some $0 < \theta \le 1/2$; (ii)
$\EE[(\Delta \tau_1)^2] < \infty$ and $\EE[(\Delta B_1)^2] < \infty$;
and (iii) $\EE[\exp\{Q(\Delta B_n^{\ast})\}] < \infty$ ($n=0,1$) for
some $Q \in \SC$ such that $x^{\theta} = O(Q(x))$. Under these
conditions, $\PP(B(T) > x) \simhm{x} \PP(M(T) > x) \simhm{x} \PP(T >
x)$.
\end{thm}

\proof See Appendix~\ref{proof-thm-main}. \qed

\begin{rem}
We use Lemma~\ref{prop-LD-B(t)}~(i) to prove $\PP(M(T) > x)
\lesimhm{x} \PP(T > x)$. For this purpose, conditions (ii) and (iii)
are assumed. Further the proof of $\PP(B(T) > x) \gesimhm{x} \PP(T >
x)$ requires the central limit theorem (CLT) for $\{B(t)\}$, which
holds under condition (ii) (see, e.g., \cite[Chapter~VI,
  Theorem~3.2]{Asmu03}).
\end{rem}

Theorem~\ref{thm-independ-01} is a generalization of
Proposition~\ref{prop-Jele04}. Condition (i) of
Theorem~\ref{thm-independ-01} is equivalent to condition (i) of
Proposition~\ref{prop-Jele04}, i.e., $T \in \calL^2$ (see
Lemma~\ref{lem-01}~(ii)). Condition (ii) of
Theorem~\ref{thm-independ-01} is weaker than the corresponding
condition of Proposition~\ref{prop-Jele04} because the positivity of
$\Delta B_n$ and condition (iii) of Proposition~\ref{prop-Jele04}
imply $\EE[(\Delta B_1)^2] < \infty$ (see
Remark~\ref{remark-SC'-02}). In addition, if $Q(x) = \eta \sqrt{x}$
for some $\eta > 0$, then condition (iii) of
Theorem~\ref{thm-independ-01} is reduced to condition (iii) of
Proposition~\ref{prop-Jele04}. As a results, the conditions of
Theorem~\ref{thm-independ-01} are weaker than those of
Proposition~\ref{prop-Jele04}.

\begin{thm}\label{thm-independ-02}
Suppose that $T$ is independent of $\{B(t);t\ge0\}$. Further suppose
that (i) $T \in \calC$; (ii) $\EE[\sup_{\tau_0 \le t \le \tau_1}|B(t)
  - B(\tau_0)|] < \infty$ and $\EE[(\Delta \tau_1)^2] < \infty$; (iii)
$\PP(\Delta B_n^{\ast} > x) = o(\PP(T > x))$ ($n=0,1$); (iv)
$x\PP(\Delta B_1 - \Delta \tau_1 > x) = o(\PP(T > x))$; and (v) either
of the following is satisfied:

\begin{narrow}
\begin{enumerate}
\renewcommand{\labelenumi}{(\alph{enumi})}
\item $\EE[\{(\Delta B_1 - \Delta \tau_1)^+\}^r] < \infty$ for some $r
  > 1$; or
\item $\int_y^{\infty}x^{-1}\PP(T>x)\rd x < \infty$ for some $y \in
  (0,\infty)$.
\end{enumerate}
\end{narrow}
Under these conditions, $\PP(B(T) > x)
\simhm{x} \PP(M(T) > x) \simhm{x} \PP(T > x)$.
\end{thm}

\proof See Appendix~\ref{proof-thm-independ-02}. \qed

\begin{rem}\label{rem-SLLN}
We prove the asymptotic upper bound $\PP(M(T) > x) \lesimhm{x} \PP(T >
x)$ of Theorem~\ref{thm-independ-02} in a similar way to that of
Theorem~\ref{thm-dependent-02}. To do this, we require condition (i),
$\EE[(\Delta \tau_1)^2] < \infty$ and conditions (iii)--(v). On the
other hand, we prove the asymptotic lower bound $\PP(B(T) > x)
\gesimhm{x} \PP(T > x)$ of Theorem~\ref{thm-independ-02} by using the
strong law of large numbers (SLLN) for $\{B(t)\}$, i.e.,
$\lim_{t\to\infty}B(t)/t = b$ w.p.1, which requires $\EE[\sup_{\tau_0
    \le t \le \tau_1}|B(t) - B(\tau_0)|] < \infty$ in condition (ii)
(see \cite[Chapter~VI, Theorem~3.1]{Asmu03}).
\end{rem}

We make a comparison of Theorem~\ref{thm-independ-02} with
Propositions~\ref{prop-Ales08-thm2.1} and \ref{prop-Robe08-thm3.2}.
Suppose that $\{B(t)\}$ satisfies (\ref{cond-iid-sums}). It then
follows that conditions (i)--(v) of Theorem~\ref{thm-independ-02} are
reduced to the following:

\begin{narrow}
\begin{enumerate}
\renewcommand{\labelenumi}{(\Roman{enumi})}
\item $T \in \calC$; 
\item $x\PP(X > x) = o(\PP(T > x))$; and
\item $\EE[X^r] < \infty$ for some $r > 1$ or
  $\int_y^{\infty}x^{-1}\PP(T > x)\rd x < \infty$ for some $y \in
  (0,\infty)$.
\end{enumerate}
\end{narrow}

\smallskip

Theorem~\ref{thm-independ-02} does not necessarily require either the
condition $\EE[T] < \infty$ of Proposition~\ref{prop-Ales08-thm2.1} or
condition (ii) of Proposition~\ref{prop-Robe08-thm3.2}. On the other
hand, Proposition~\ref{prop-Ales08-thm2.1} does not necessarily
require condition (II) (which is obvious). Further we can confirm that
condition (II) is not necessary for
Proposition~\ref{prop-Robe08-thm3.2}, as follows.

Suppose that $\PP(T > x) \simhm{x} x^{-\alpha}$ for some $0 < \alpha <
1$. In this case, $\EE[T] = \infty$ and $T \in \calC$ (see
Appendix~\ref{subsec-class-R}), which shows that condition (i) of
Proposition~\ref{prop-Robe08-thm3.2} is satisfied. In addition,
$\EE[T\cdot\dd{1}(T \le x)] = O(x\PP(T > x))$ (see Remark below
Theorem~3.2 in \cite{Robe08}). Therefore condition (iii) of
Proposition~\ref{prop-Robe08-thm3.2} holds for $q=1$.  We now assume
that $\PP(X > x) = (x+1)^{-\beta}$ for some $1 < \beta < \alpha +
1$. We then have $\EE[X^r] < \infty$ for all $r < \beta$, and thus
condition (ii) of Proposition~\ref{prop-Robe08-thm3.2} is
satisfied. As a result, all the conditions of
Proposition~\ref{prop-Robe08-thm3.2} hold, whereas condition (II) does
not hold.

The above discussion shows that Theorem~\ref{thm-independ-02} is not a
complete generalization of Propositions~\ref{prop-Ales08-thm2.1} and
\ref{prop-Robe08-thm3.2}.

\subsection{Special case: discrete-time Markov additive process}\label{subsec-MAdP}

In this subsection, we extend Propositions~\ref{prop-Ales08-thm2.1}
and \ref{prop-Robe08-thm3.2} to the random sum of (possibly negative)
r.v.s with Markovian correlation. For this purpose, we introduce a
discrete-time Markov additive process.

Let $\{J_n;n=0,1,\dots\}$ is a discrete-time Markov chain with a
finite state space $\bbD:=\{0,1,\dots,d-1\}$. Let $X_n$'s
($n=0,1,\dots$) denote r.v.s such that for all $i,j \in \bbD$ and $x
\in \bbR$,
\begin{eqnarray*}
\PP(X_0 \le x, J_0 = i) &=& \beta_{i}(x),
\\
\PP(X_{n+1} \le x, J_{n+1} = j \mid J_n = i) &=& H_{i,j}(x),\quad n=0,1,\dots,
\end{eqnarray*}
where $\sum_{i\in\bbD}\beta_{i}(\infty) = 1$ and
$\sum_{j\in\bbD}H_{i,j}(\infty) = 1$ for all $i \in \bbD$.  Let
$S_n=\sum_{\nu=0}^n X_{\nu}$ for $n=0,1,\dots$.  It then follows that
$\{(S_n,J_n);n=0,1,\dots\}$ is a Markov additive process with initial
distribution $\vc{\beta}(x)=(\beta_i(x))_{i\in\bbD}$ and Markov
additive kernel (called ``kernel" for short)
$\vc{H}(x)=(H_{i,j}(x))_{i,j\in\bbD}$ ($x \in \bbR$).  Further let
$\widehat{\vc{\beta}}(\xi)$ and $\widehat{\vc{H}}(\xi)$ denote the
characteristic functions of $\vc{\beta}(x)$ and $\vc{H}(x)$, i.e.,
\[
\widehat{\vc{\beta}}(\xi)=\int_{x \in \bbR} e^{\rmi \xi x} \rd \vc{\beta}(x),
\qquad
\widehat{\vc{H}}(\xi)=\int_{x \in \bbR} e^{\rmi \xi x} \rd \vc{H}(x),
\]
respectively, where $\rmi = \sqrt{-1}$. 

In what follows, we make the following assumption:
\begin{assumpt}\label{assumpt-(S_n,J_n)}

\begin{enumerate}
\item Let $B(t) = S_{\lfloor t \rfloor} =\sum_{n=0}^{\lfloor t
  \rfloor}X_n$ for $t \ge 0$;
\item the background process $\{J_n\}$ is irreducible, i.e.,
  $\vc{H}(\infty)$ is an irreducible stochastic
  matrix; and
\item the mean drift of the additive component $\{S_n\}$ is finite and
  positive, i.e.,
\begin{equation}
h := \vc{\varpi}\int_{x \in \bbR} x \rd \vc{H}(x) \vc{e} \in (0,\infty),
\label{defn-sigma}
\end{equation}
where $\vc{\varpi}=(\varpi_i)_{i\in\bbD}$ is the stationary
probability vector of $\vc{H}(\infty)$, and where $\vc{e}$ is a column
vector of ones with an appropriate dimension.
\end{enumerate}
\end{assumpt}

It is easy to see that $\{B(t);t\ge0\}$ is a cumulative process
because $\{(B(n),J_n);n=0,1,\dots\}$ is a discrete-time Markov
additive process. Let $0 \le \tau_0 < \tau_1 < \cdots$ denote hitting
times of $\{J_n\}$ to state zero, which are regenerative points of the
cumulative process $\{B(t)\}$. Clearly, $\Delta \tau_1 \ge 1$
w.p.1. Further from (\ref{defn-Delta-tau_n}), we have $\tau_0 = \Delta
\tau_0$ and thus $\PP(\Delta \tau_0 = 0) = \PP(J_0 = 0) =
\beta_0(\infty)$.

Let $\widehat{\psi}_0(z,\xi) =\EE[z^{\Delta \tau_0}
  e^{\rmi \xi \Delta B_0}]$ and $\widehat{\psi}_1(z,\xi)
=\EE[z^{\Delta \tau_1} e^{\rmi \xi \Delta B_1}]$. We then have
\begin{eqnarray}
\widehat{\psi}_0(z,\xi)
&=& \widehat{\beta}_{0}(\xi) + \widehat{\vc{\beta}}_{+}(\xi)
\left( \vc{I} - z\widehat{\vc{H}}_{+}(\xi) \right)^{-1}
z\widehat{\vc{h}}_{+}(\xi),
\label{joint-transform-0}
\\
\widehat{\psi}_1(z,\xi)
&=& z\widehat{H}_{0,0}(\xi) + z\widehat{\vc{\eta}}_{+}(\xi)
\left(\vc{I} - z\widehat{\vc{H}}_{+}(\xi)\right)^{-1}
z\widehat{\vc{h}}_{+}(\xi),
\label{joint-transform-1}
\end{eqnarray}
where $\vc{I}$ denotes the identity matrix with an appropriate
dimension and
\[
\widehat{\vc{\beta}}(\xi)
= 
\bordermatrix{
               & \{0\} &   \bbD\setminus\{0\}       
\cr
 & \widehat{\beta}_{0}(\xi) & \widehat{\vc{\beta}}_{+}(\xi)
},
\quad
\widehat{\vc{H}}(\xi)
= 
\bordermatrix{
               & \{0\} &   \bbD\setminus\{0\}       
\cr
\{0\} & \widehat{H}_{0,0}(\xi) & \widehat{\vc{\eta}}_{+}(\xi)
\cr
\bbD\setminus\{0\} & \widehat{\vc{h}}_{+}(\xi) & \widehat{\vc{H}}_{+}(\xi)
}.
\]
The first term of (\ref{joint-transform-0}) corresponds to the event
where $J_0 = 0$ and thus $\Delta \tau_0 = 0$.  The first term of
(\ref{joint-transform-1}) corresponds to the event where a
regenerative cycle lasts only for a unit of time, i.e., the background
process $\{J_n\}$ moves from state zero to state zero in one
transition. As for the second terms of (\ref{joint-transform-0}) and
(\ref{joint-transform-1}), they correspond to the events where
$\{J_n\}$ moves from state zero to a state in $\bbD\setminus\{0\}$ and
then eventually returns to state zero.

Fixing $\xi = 0$ in (\ref{joint-transform-0}) and
(\ref{joint-transform-1}) and taking the inverse of them with respect
to $z$, we have
\begin{align}
\PP(\Delta \tau_0 = k)
&= \dd{1}(k=0) 
\widehat{\beta}_{0}(0) 
+ \dd{1}(k\ge1)\widehat{\vc{\beta}}_{+}(0) 
\left(\widehat{\vc{H}}_{+}(0) \right)^{k-1}
\widehat{\vc{h}}_{+}(0),
\nonumber
\\
\PP(\Delta \tau_1 = k)
&= 
\dd{1}(k=1) 
\widehat{H}_{0,0}(0)
+ \dd{1}(k\ge2)\widehat{\vc{\eta}}_{+}(0)
\left(\widehat{\vc{H}}_{+}(0)\right)^{k-2}
\widehat{\vc{h}}_{+}(0),
\label{dist-Delta-tau_1}
\end{align}
for $k=0,1,2,\dots$.  Therefore $\Delta \tau_0$ and $\Delta \tau_1$
follow discrete phase-type distributions \cite{Lato99}. Further we
have the following result by combining the renewal reward theory (see,
e.g., \cite[Chapter~2, Theorem~2]{Wolf89}) and the discrete-time
version of the ergodic theorem (see, e.g., \cite[Chapter~3,
  Theorem~4.1]{Brem99}):
\begin{prop}\label{prop-b=h}
Under Assumption~\ref{assumpt-(S_n,J_n)},
\[
b 
:= {\EE[\Delta B_1] \over \EE[\Delta \tau_1]} 
= \vc{\varpi}\int_{x \in \bbR}x \rd\vc{H}(x)\vc{e} = h \in (0,\infty).
\]
\end{prop}

In what follows, we present two theorems that supersede
Propositions~\ref{prop-Ales08-thm2.1} and
\ref{prop-Robe08-thm3.2}. Before doing this, we introduce three lemmas
for the proofs of the theorems.
\begin{lem}\label{add-lem-MAdP-01}
Suppose that Assumptions~\ref{assumpt-(S_n,J_n)} holds. Further let
$\overline{\vc{\beta}}(x) = \int_x^{\infty} \rd \vc{\beta}(y)$ and
$\overline{\vc{H}}(x) = \int_x^{\infty} \rd \vc{H}(y)$ for $x \in
\bbR$ and suppose that there exist some $\tilde{c} \in [0,\infty)$ and
  some nonnegative r.v.\ $Y \in \calS$ such that
\[
\limsup_{x \to \infty}{\overline{\vc{\beta}}(x) \over \PP(Y > x)} 
\le \tilde{c}\widetilde{\vc{\beta}},
\qquad
\limsup_{x \to \infty}{\overline{\vc{H}}(x) \over \PP(Y > x)} 
\le \tilde{c}\widetilde{\vc{H}},
\]
where $\widetilde{\vc{\beta}} = (\widetilde{\beta}_i)_{i\in\bbD}$ is a
finite nonnegative vector and $\widetilde{\vc{H}} =
(\widetilde{H}_{i,j})_{i,j\in\bbD}$ is a finite nonnegative matrix.
We then have
\[
\limsup_{x\to\infty}
{
\PP(\Delta B_n > x) \over \PP(Y > x)}
\le\tilde{c} C,\quad n=0,1.
\]
\end{lem}

\proof See Appendix~\ref{proof-add-lem-MAdP-01}. \qed

\begin{lem}\label{add-lem-03}
If the assumptions of Lemma~\ref{add-lem-MAdP-01} are satisfied, then
\[
\limsup_{x\to\infty}
{\PP(\Delta B_1 > x \mid \Delta \tau_1 = k) \over \PP(Y > x)}
\le \tilde{c} C k, \qquad \forall k=1,2,\dots,
\]
where $C$ is independent of $k$.
\end{lem}

\proof See Appendix~\ref{proof-add-lem-03}. \qed

\begin{lem}\label{lem-MAdP-03}
If the assumptions of Lemma~\ref{add-lem-MAdP-01} are satisfied, then
for all $t \ge 0$ and $m=0,1,\dots$,
\begin{equation}
\limsup_{x\to\infty}
{ \PP\left(\left. \sum_{i=1}^m \Delta B_i  > x \,\right| N(t) = m \right)
\over \PP(Y > x)}
\le \tilde{c} C t,
\label{add-eqn-29}
\end{equation}
where  $N(t) =
\max\{k\ge0;\sum_{i=1}^k \Delta \tau_i \le t\}$ for $t \in \bbR$.
\end{lem}

\proof See Appendix~\ref{proof-lem-MAdP-03}. \qed

\smallskip

The following theorems present two sets of conditions for
(\ref{Jele-asymp}). Note here that under
Assumption~\ref{assumpt-(S_n,J_n)}, the asymptotic equality
(\ref{Jele-asymp}) is reduced to $\PP(S_{\lfloor T \rfloor} > h x)
\simhm{x} \PP(M_{\lfloor T \rfloor} > h x) \simhm{x} \PP(T > x)$,
where $M_n = \max_{0 \le k \le n}S_k$.
\begin{thm}\label{thm-MAdP-01}
Suppose that Assumption~\ref{assumpt-(S_n,J_n)} holds and $T$ is
independent of the Markov additive process $\{(S_n,J_n)\}$. Further
suppose that $T \in \calC$, $\EE[T] < \infty$ and
\begin{equation}
\int_{|y| > x} \rd \vc{\beta}(y) = o( \PP(T > x) ),
\qquad
\int_{|y| > x} \rd \vc{H}(y)
 = o( \PP(T > x) ).
\label{MAdP-assumpt-01}
\end{equation}
Under these conditions, $\PP(S_{\lfloor T \rfloor} > h x) \simhm{x}
\PP(M_{\lfloor T \rfloor} > h x) \simhm{x} \PP(T > x)$.
\end{thm}

\proof See Appendix~\ref{proof-thm-MAdP-01}. \qed

\begin{thm}\label{thm-MAdP-02}
Suppose that Assumption~\ref{assumpt-(S_n,J_n)} holds and $T$ is
independent of the Markov additive process $\{(S_n,J_n)\}$. Further
suppose that $T \in \calC$ and there exists some nonnegative r.v.\ $Y
\in \calS$ such that
\begin{eqnarray}
&&
\int_{|y| > x} \rd \vc{\beta}(y) = O( \PP(Y > x) ),
\qquad
\int_{|y| > x} \rd \vc{H}(y) = O( \PP(Y > x) ),
\label{MAdP-assumpt-02}
\\
&&
\lim_{x\to\infty}\EE[T \cdot \dd{1}(T \le x, N(T) \le x/\EE[\Delta \tau_1])]
{\PP(Y>x) \over \PP(T>x)} = 0.
\label{MAdP-assumpt-03}
\end{eqnarray}
Under these conditions, $\PP(S_{\lfloor T \rfloor} > h x) \simhm{x}
\PP(M_{\lfloor T \rfloor} > h x) \simhm{x} \PP(T > x)$.
\end{thm}

\proof See Appendix~\ref{proof-thm-MAdP-02}. \qed

\begin{rem}\label{rem-thm-MAdP-02}
Equation (\ref{MAdP-assumpt-03}) holds if
\[
\lim_{x\to\infty}\EE[T \cdot \dd{1}(T \le x)]
{\PP(Y>x) \over \PP(T>x)} = 0.
\]
\end{rem}

It is easy to see that Theorems~\ref{thm-MAdP-01} and
\ref{thm-MAdP-02} include Propositions~\ref{prop-Ales08-thm2.1} and
\ref{prop-Robe08-thm3.2}, respectively, as special cases.

\subsection{Special case: continuous-time Markov additive process}\label{subsec-MAdP-02}

In this subsection, we consider a continuous-time Markov additive
process $\{(B(t),J(t));t \ge 0\}$ with state space $\bbR \times \bbD$,
where $\{B(t)\}$ is the additive component and $\{J(t)\}$ is the
background process.  Let $\vc{D}(x) = (D_{i,j}(x))_{i,j\in\bbD}$ ($x
\in \bbR$) denote the kernel of $\{(B(t),J(t))\}$ such that $\vc{D}(x)
\ge \vc{O}$ for all $x < 0$ and $\vc{D}(x) - \vc{D}(0) \ge \vc{O}$ for
all $x \ge 0$, where $\vc{O}$ denotes the zero matrix. Further for
later use, let $\widehat{\vc{D}}(\xi) = \int_{x\in\bbR} e^{\rmi\xi x}
\rd \vc{D}(x)$ and $[\,\cdot\,]_{i,j}$ denote the $(i,j)$th element of
the matrix between square brackets.

In what follows, we make the following assumption:
\begin{assumpt}\label{assumpt-continu-MAdP}
(i) For all $t \ge 0$,
\begin{equation}
\EE[\exp\{\rmi\xi B(t)\} \cdot \dd{1}(J(t) = j) \mid J(0) = i]
=\left[ \exp\{\widehat{\vc{D}}(\xi) t\}  \right]_{i,j},
\quad i,j \in \bbD;
\label{eqn-MAdP-01}
\end{equation}
(ii)
  $\widehat{\vc{D}}(0) = \vc{D}(\infty)$ is an irreducible
  infinitesimal generator; and
(iii) $\vc{\pi}\int_{x\in\bbR}x \rd \vc{D}(x)\vc{e} \in (0,\infty)$,
  where $\vc{\pi}=(\pi_i)_{i\in\bbD}$ denotes the stationary
  probability vector of $\widehat{\vc{D}}(0)$.
\end{assumpt}

Under Assumption~\ref{assumpt-continu-MAdP}, $\{B(t)\}$ is a
cumulative process. It thus follows from the renewal reward theory
(see, e.g., \cite[Chapter~2, Theorem~2]{Wolf89}) and the
continuous-time version of the ergodic theorem (see, e.g.,
\cite[Chapter~8, Theorem~6.2]{Brem99}) that
\[
b:= {\EE[\Delta B_1] \over \EE[\Delta \tau_1]}
= \vc{\pi}\int_{x\in\bbR}x \rd \vc{D}(x)\vc{e} \in (0,\infty).
\]
Further it follows from (\ref {eqn-MAdP-01}) that
\begin{eqnarray}
\lefteqn{
\EE[\exp\{\rmi\xi B(T)\} \cdot \dd{1}(J(T) = j) \mid J(0) = i]
}
\quad &&
\nonumber
\\
&=& 
\left[ 
\int_0^{\infty} \exp\{\widehat{\vc{D}}(\xi) t\} \rd \PP(T \le t) 
\right]_{i,j}
\nonumber
\\
&=& \sum_{n=0}^{\infty}
\int_0^{\infty}  e^{-\gamma t}{(\gamma t)^n \over n!}\rd \PP(T \le t) 
\cdot
\left[ 
\left\{ \vc{I} + \gamma^{-1} \widehat{\vc{D}}(\xi) \right\}^n 
\right]_{i,j}
\nonumber
\\
&=& \sum_{n=0}^{\infty} p_n 
\cdot 
\left[ 
\{ \widehat{\vc{K}}(\xi) \}^n 
\right]_{i,j},
\label{eqn-MAdP-02}
\end{eqnarray}
where 
\[
p_n
= \int_0^{\infty}  e^{-\gamma t}{(\gamma t)^n \over n!}\rd \PP(T \le t)~~(n=0,1,\dots),
~~
\widehat{\vc{K}}(\xi)
= \vc{I} + \gamma^{-1} \widehat{\vc{D}}(\xi),
~~
\gamma = \max_{i\in\bbD}|D_{i,i}(\infty)|.
\]
Note here that $\widehat{\vc{K}}(\xi)$ is the characteristic function
of $\vc{K}(x) := \dd{1}(x\ge0)\vc{I} + \gamma^{-1} \vc{D}(x)$ ($x \in
\bbR$), which can be considered as the kernel of a discrete-time
Markov additive process $\{(S_n,J_n)\}$ discussed in the previous
subsection.  Note also that $\{p_n;n=0,1,\dots\}$ is the distribution
of the counting process of Poisson arrivals with rate $\gamma$ during
time interval $(0,T]$. It is easy to see that if $T \in \calC$, then
the counting process satisfies all the conditions of
Theorem~\ref{thm-independ-02} and thus
\[
\sum_{n=k+1}^{\infty}p_n \simhm{k} \PP(T > k/\gamma).
\]

We now define $T'$ as a nonnegative integer-valued r.v.\ such that
$\PP(T' = n) = p_n$ ($n=0,1,\dots$) and $T'$ is independent of a
discrete-time Markov additive process $\{(S_n,J_n)\}$ with initial
condition $S_0=X_0=0$ (i.e., $\int_{\{0\}} \rd \vc{\beta}(x)\vc{e} =
1$) and kernel $\vc{H}(x) = \vc{K}(x)$ ($x \in \bbR$). It then follows
from (\ref{eqn-MAdP-02}) that $(B(T),J(T))$ is stochastically
equivalent to $\{(S_{T'},J_{T'})\}$.  As a result, using
Theorems~\ref{thm-MAdP-01} and \ref{thm-MAdP-02}, we can readily prove
the following corollaries, whose proofs are omitted.
\begin{coro}\label{coro-MAdP-01}
Suppose that Assumption~\ref{assumpt-continu-MAdP} holds and $T$ is
independent of the Markov additive process $\{(B(t),J(t))\}$. Further
suppose that $T \in \calC$, $\EE[T] < \infty$ and
\[
\int_{|y| > x} \rd \vc{D}(y)
 = o( \PP(T > x) ).
\]
Under these conditions, $\PP(B(T) > b x) \simhm{x} \PP(M(T) > b x)
\simhm{x} \PP(T > x)$.
\end{coro}

\begin{coro}\label{coro-MAdP-02}
Suppose that Assumption~\ref{assumpt-continu-MAdP} holds and $T$ is
independent of the Markov additive process $\{(B(t),J(t))\}$. Further
suppose that $T \in \calC$ and there exists some nonnegative r.v.\ $Y
\in \calS$ such that
\begin{eqnarray*}
\int_{|y| > x} \rd \vc{D}(y) = O( \PP(Y > x) ),
\quad
\lim_{x\to\infty}\EE[T \cdot \dd{1}(T \le x)]
{\PP(Y>x) \over \PP(T>x)} = 0.
\end{eqnarray*}
Under these conditions, $\PP(B(T) > b x) \simhm{x} \PP(M(T) > b x)
\simhm{x} \PP(T > x)$.
\end{coro}

\section{Application}\label{sec-application-queue}

In this section, we first introduce a new (discrete-time) on/off
arrival process, ON/OFF-BMAP, mentioned in
Section~\ref{introduction}. We then consider a single-server
finite-buffer queue with an ON/OFF-BMAP and deterministic service
times. For this queueing model, we derive some subexponential
asymptotic formulas for the loss probability by using the main results
presented in Section~\ref{sec-main-result}.

\subsection{ON/OFF batch Markovian arrival process}

We describe the definition of ON/OFF-BMAPs in discrete time.  The time
interval $[n,n+1]$ ($n\in\bbZ$) is called slot $n$, where
$\bbZ=\{0,\pm1,\pm2\dots\}$. The ON/OFF-BMAP is an on/off arrival
process, where on-periods and off-periods are repeated
alternately. For simplicity, slots in on-periods (resp.\ off-periods)
are called {\it on-slots} (resp.\ {\it off-slots}).

The lengths of off-periods are i.i.d., and no arrivals occur in any
off-slot.  On the other hand, at least one arrival occurs in each
on-slot w.p.1 and the number of arrivals during each on-period follows
a BMAP started with some initial distribution at the beginning of the
on-period. Further the lengths of on-periods are i.i.d., but the
length of each on-period may depend on the BMAP in the on-period. In
what follows, the BMAP in the $m$th ($m\in\bbZ$) on-period is called
the $m$th BMAP.

To describe the ON/OFF-BMAP more precisely, we define some
notations. Let $N_{m,n}$ ($m\in\bbZ$, $n=0,1,\dots$) denote the number
of arrivals in the $n$th slot of the $m$th on-period. Let
$J_{m,0},J_{m,1},J_{m,2},\dots$ ($m\in\bbZ$) denote the background
states of the $m$th BMAP, which belong to $\bbD=\{0,1,\dots,d-1\}$. We
then assume that
\begin{equation}
\PP(N_{m,0} = k, J_{m,0} = i) 
= \alpha_i(k), \qquad  i \in \bbD,~k =1,2,\dots,
\label{cond-N_{m,0}}
\end{equation}
where $\vc{\alpha}(k)=(\alpha_i(k))_{i\in\bbD}$ is a $1 \times d$
nonnegative vector such that
$\underline{\vc{\alpha}}:=\sum_{k=1}^{\infty}\vc{\alpha}(k)$ is a
probability vector. We also assume that for $n=1,2,\dots$,
\begin{equation}
\PP(N_{m,n} = k, J_{m,n} = j \mid J_{m,n-1} = i) 
= \varLambda_{i,j}(k), \qquad i,j \in \bbD,~k =1,2,\dots,
\label{cond-N_{m,n}}
\end{equation}
where $\vc{\varLambda}(k)=(\varLambda_{i,j}(k))_{i,j\in\bbD}$ is a $d
\times d$ substochastic matrix such that
$\underline{\vc{\varLambda}}:=\sum_{k=1}^{\infty}\vc{\varLambda}(k)$
is an irreducible stochastic matrix.

Let $I_m^\on$ ($m\in\bbZ$) denote the length of the $m$th on-period.
Let $\Phi_m$ ($m \in \bbZ$) denote
\begin{equation}
\Phi_m =\{I_m^\on,(N_{m,0},J_{m,0}),(N_{m,1},J_{m,1}),
\dots,(N_{m,I_m^\on-1},J_{m,I_m^\on-1})\}.
\label{defn-Phi_m}
\end{equation}
We then assume that the $\Phi_m$'s ($m \in \bbZ$) are i.i.d. Thus the
$I_m^\on$'s ($m\in\bbZ$) are i.i.d.\ r.v.s, though each $I_m^\on$ may
depend on the $m$th BMAP, i.e., $\{(N_{m,n},J_{m,n});n=0,1,\dots\}$.

For later use, let $\lambda$ denote the arrival rate during on
periods, i.e., the time-average number of arrivals in an on-slot. It
follows from the ergodic theorem (see, e.g., \cite[Chapter~3,
  Theorem~4.1]{Brem99}) that
\begin{equation}
\lambda = \vc{\phi}\sum_{k=1}^{\infty}k\vc{\varLambda}(k)\vc{e} \ge 1,
\label{defn-lambda}
\end{equation}
where $\vc{\phi}=(\phi_i)_{i\in\bbD}$ denotes the
stationary probability vector of $\underline{\vc{\varLambda}}$.

\begin{rem}\label{rem-ON/OFF-BMAP-01}
The ON/OFF-BMAP is a generalization of the batch-on/off process
introduced by Galm\'es and Puigjaner~\cite{Galm01}. In the
batch-on/off process, the numbers of arrivals in individual on-slots
are i.i.d.\ and independent of the lengths of on-periods. Based on the
Wiener-Hopf factorization (see, e.g., \cite[Chapter~VIII,
  Section~3]{Asmu03}), Galm\'es and Puigjaner~\cite{Galm03,Galm05}
study the response time distribution of a single-server queue with a
batch-on/off process and deterministic service times.
\end{rem}

\begin{rem}\label{rem-PAP}
The ON/OFF-BMAP is similar to the PAP proposed by Alfa and
Neuts~\cite{Alfa95} and Breuer and Alfa~\cite{Breu05}. The PAP can be
considered as a special case of the ON/OFF-BMAP in the sense that the
lengths of on-periods (resp.\ off-periods) follow a phase-type
distribution. However, the PAP allows that no arrivals occur in an
{\it on-slot}.
\end{rem}

\subsection{Loss probability of (ON/OFF-BMAP)/D/1/$\vc{K}$ queue}

We begin with the description of our queueing model.  Customers arrive
at the system according to an ON/OFF-BMAP.  The system has a single
server and a buffer of finite capacity $K-1$ (thus the system capacity
is equal to $K$).  The service times of customers are all equal to the
length of one slot.  According to Kendall's notation, our queueing
model is denoted by (ON/OFF-BMAP)/D/1/$K$.

For analytical convenience, we assume that arrivals in each on-slot
occur at the same time, immediately after the beginning of the
on-slot. We also assume that departure points are located immediately
before the ends of slots. Under these assumptions, we observe the
queue length process immediately after the ends of off-periods.

Let $L_m^{(K)}$ ($m\in\bbZ$) denote the queue length immediately after
the end of the $m$th off-period.  Let $I_m^\off$ ($m\in\bbZ$) denote
the length of the $m$th off-period, where the $I_m^\off$'s are
i.i.d.\ r.v.s. Further let $A_m$ ($m\in\bbZ$) denote the increment in
the queue length during the $m$th on-period, i.e.,
\begin{equation}
A_m = \sum_{n=0}^{I_m^\on-1} (N_{m,n} - 1),
\label{defn-A_m}
\end{equation}
where the $A_m$'s are i.i.d.\ r.v.s because the $\Phi_m$'s in
(\ref{defn-Phi_m}) are i.i.d.  We then have
\[
L_{m+1}^{(K)} = (\min(L_m^{(K)} + A_{m+1}, K) - I_{m+1}^\off)^+.
\]

We now define $P_{\rm loss}^{(K)}$ as the loss probability, which is
the time-average of losses. Note that in the $m$th renewal cycle
consisting of the $m$th on- and off-periods, the numbers of arrivals
and losses are equal to $A_m + I_m^\on$ and $(L_{m-1}^{(K)} + A_m - K)^+$,
respectively.  It then follows from the renewal reward theory (see,
e.g., \cite[Chapter~2, Theorem~2]{Wolf89}) that
\[
P_{\rm loss}^{(K)}
= {\EE[(L_{m-1}^{(K)} + A_m - K)^+] \over \EE[A_m + I_m^\on]}.
\]

\subsection{Subexponential asymptotics of the loss probability}

In this subsection, we derive some subexponential asymptotic formulas
for the loss probability $P_{\rm loss}^{(K)}$.  To achieve this, we
combine our main results with the following proposition:
\begin{prop}[Theorem 5 in \cite{Jele99}]\label{prop-Jele99}
Let $A$, $I^\on$ and $I^\off$ denote generic r.v.s for
i.i.d.\ sequences $\{A_m\}$, $\{I_m^\on\}$ and $\{I_m^\off\}$,
respectively.  Suppose $0 < \EE[A] < \infty$ and let $A_{\re}$ denote
the equilibrium r.v.\ of $A$, i.e., $\PP(A_{\re} \le x) =
(1/\EE[A])\int_0^x\PP(A > y)\rd y$ for $x \ge 0$. If $\EE[A] <
\EE[I^\off]$ and $A_{\re} \in \calS$, then
\[
P_{\rm loss}^{(K)} 
\simhm{K} {\EE[(A-K)^+] \over \EE[A] + \EE[I^\on]}
= {\EE[A] \over \EE[A] + \EE[I^\on]} \PP(A_{\re} > K).
\]
\end{prop}

In the rest of this subsection, we set
\begin{eqnarray}
T &=& I_m^\on - 1,
\quad 
B(t) = \sum_{n=0}^{\lfloor t \rfloor} (N_{m,n} - 1), \quad  t \ge 0,
\label{defn-B(t)}
\end{eqnarray}
where $N_{m,n} - 1$ for $m \in \bbZ$ and $n=0,1,\dots$ due to
(\ref{cond-N_{m,0}}) and (\ref{cond-N_{m,n}}).  It then follows from
(\ref{defn-A_m}) and (\ref{defn-B(t)}) that
\begin{equation}
A \stackrel{\rd}{=} B(T),
\label{eqn-A}
\end{equation}
where the symbol $\stackrel{\rd}{=}$ denotes equality in distribution.

For simplicity, let $X_n = N_{m,n} -
1 \ge 0$ and $J_n = J_{m,n}$ for $n=0,1,\dots$. Further let $S_n =
\sum_{\nu=0}^n X_{\nu}$ for $n=0,1,\dots$. It then follows that $B(t) =
S_{\lfloor t \rfloor}$ for $t \ge 0$ and $\{(S_n,J_n);n=0,1,\dots\}$
is a Markov additive process with state space $\{0,1,\dots\} \times
\bbD$, initial distribution $\vc{\alpha}(k)$ and Markov additive
kernel $\vc{\varLambda}(k+1)$ ($k=0,1,\dots$). Thus the stochastic
process $\{B(t)\}$ defined in (\ref{defn-B(t)}) is a cumulative
process of the same type as that in subsection~\ref{subsec-MAdP}.  As
with subsection~\ref{subsec-MAdP}, let $0 \le \tau_0 < \tau_1 <
\cdots$ denote hitting times of $\{J_n\}$ to state zero, which are
regenerative points of $\{B(t)\}$. From (\ref{defn-lambda}) and
Proposition~\ref{prop-b=h}, we have
\begin{equation}
b := {\EE[\Delta B_1] \big/ \EE[\Delta \tau_1]} = \lambda -1.
\label{eqn-b}
\end{equation}

We now assume the following:
\begin{assumpt}\label{assumpt-alpha(k)-D(k)}
$\lambda > 1$ and there exists some nonnegative r.v.\ $Y \in \calS$
  such that
\[
\sum_{l=k+1}^{\infty}\vc{\alpha}(l) = O( \PP(Y > k) ),
\qquad
\sum_{l=k+1}^{\infty}\vc{\varLambda}(l) = O( \PP(Y > k) ).
\]
\end{assumpt}

In what follows, we present four subexponential asymptotic formulas
for the loss probability $P_{\rm loss}^{(K)}$. The first two formulas
are obtained from the results for the dependent-sampling case, and the
others are from those for the independent-sampling case.

\begin{thm}\label{thm-app1-depend-01}
Suppose that Assumption~\ref{assumpt-alpha(k)-D(k)} holds and $\EE[A]
< \EE[I^\off]$. Further suppose that (i) $I^\on \in \calL^{1/\theta}$
for some $0 < \theta \le 1/3$; (ii) $\EE[I^\on] < \infty$ and the
equilibrium r.v.\ $I^\on_{\re}$ of $I^\on$ is subexponential (i.e.,
$I^\on_{\re} \in \calS$); and (iii) $\EE[\exp\{Q(Y)\}] < \infty$ for
some $Q \in \SC$ such that $x^{3\theta/2} =O( Q(x) )$. We then have
\begin{equation}
P_{\rm loss}^{(K)} 
\simhm{K} {(\lambda - 1)\EE[I^\on] \over \EE[A] + \EE[I^\on]} 
\PP(I^\on_{\re} > K/(\lambda - 1)).
\label{asympt-P_{loss}}
\end{equation}
In addition, if (iv) each for $m \in \bbZ$, $\{I_m^\on \ge n+1\}$ is
independent of $N_{m,n}$ for all $n=0,1,\dots$; and (v)
$\vc{\alpha}(k) = \vc{\phi}\vc{\varLambda}(k)$ for $k=1,2,\dots$, then
\begin{equation}
P_{\rm loss}^{(K)} 
\simhm{K} {\lambda - 1 \over \lambda} 
\PP(I^\on_{\re} > K/(\lambda - 1)).
\label{asympt-P_{loss}-02}
\end{equation}
\end{thm}

\begin{rem}\label{rem-stationary-BMAP}
Condition (v) implies that BMAPs in on-periods are stationary and thus
$\PP(J_{m,n}=j) = \phi_j$ ($j \in \bbD$) for all $m \in \bbZ$ and
$n=0,1,\dots,I_m^\on-1$.
\end{rem}

\noindent
{\it Proof of Theorem~\ref{thm-app1-depend-01}}.~ We first show that
the sampling time $T$ and the cumulative process $\{B(t)\}$ in
(\ref{defn-B(t)}) satisfy the conditions of
Theorem~\ref{thm-dependent-01}.  Condition (i) of
Theorem~\ref{thm-app1-depend-01} yields
\[
\PP(T > x) = \PP(I^\on > x+1) \simhm{x} \PP(I^\on > x),
\]
and thus $T \in \calL^{1/\theta}$ for some $0 < \theta \le 1/3$, i.e.,
condition (i) of Theorem~\ref{thm-dependent-01} is satisfied.

Note here that $\{B(t)\}$ in (\ref{defn-B(t)}) satisfies $B(0) \ge 0$
and thus condition (ii) is reduced to $\EE[\exp\{Q(\Delta \tau_n)\}] <
\infty$ and $\EE[\exp\{Q(\Delta B_n)\}] < \infty$ ($n=0,1$) for some
$Q \in \SC$ such that $x^{3\theta/2} = O( Q(x) )$ (see
Remark~\ref{rem-thm-dependent-01}).  Further it follows from
Assumption~\ref{assumpt-alpha(k)-D(k)} and Lemma~\ref{add-lem-MAdP-01}
that for $n=0,1$,
\begin{equation}
\PP(\Delta B_n > x) \le C\PP(Y > x), \qquad \forall x \ge 0.
\label{add-ineqn-06}
\end{equation}
Therefore condition (iii) of Theorem~\ref{thm-app1-depend-01} implies
$\EE[\exp\{Q(\Delta B_n)\}] < \infty$ ($n=0,1$) for some $Q \in \SC$
such that $x^{3\theta/2} =O( Q(x) )$. Recall that the distribution of
$\Delta \tau_n$ is phase-type and $Q(x) = o(x)$ for any $Q \in \SC$
(see Definition~\ref{defn-SC'}). We then have $\EE[\exp\{Q(\Delta
  \tau_n)\}] < \infty$ ($n=0,1$) for any $Q \in \SC$. As a result,
condition (ii) of Theorem~\ref{thm-dependent-01} holds.

Applying Theorem~\ref{thm-dependent-01} to (\ref{eqn-A}) and using
(\ref{eqn-b}), we obtain
\begin{equation}
\PP(A > x) \simhm{x} \PP(I^\on - 1 > x/(\lambda- 1))
\simhm{x} \PP(I^\on > x/(\lambda- 1)),
\label{asymp-PP(A>k)}
\end{equation}
from which and $\EE[I^\on] < \infty$ we have $\EE[A] < \infty$ and
\begin{equation}
\PP(A_{\re} > x) 
\simhm{x} {(\lambda-1)\EE[I^\on] \over \EE[A]}\PP(I^\on_{\re} > x/(\lambda- 1)).
\label{asymp-PP(A_e>k)}
\end{equation}
Thus $A_{\re} \in \calS$ due to $I_{\re}^\on \in \calS$.  As a result,
(\ref{asymp-PP(A_e>k)}) and Proposition~\ref{prop-Jele99} yield
(\ref{asympt-P_{loss}}).

Finally, we prove (\ref{asympt-P_{loss}-02}). From (\ref{defn-A_m}),
we have
\begin{equation}
\EE[A_m+I_m^\on] 
= \EE\left[\sum_{n=0}^{I_m^\on-1}N_{m,n} \right].
\label{add-eqn-52}
\end{equation}
Note here that condition (v) of Theorem~\ref{thm-app1-depend-01} and
(\ref{defn-lambda}) yield (see Remark~\ref{rem-stationary-BMAP})
\[
\EE[N_{m,n}] = \vc{\phi}\sum_{k=1}^{\infty}k\vc{\varLambda}(k)\vc{e} = \lambda,
\qquad \forall m\in\bbZ,~ \forall n=0,1,\dots.
\]
This equation and condition (iv) of Theorem~\ref{thm-app1-depend-01}
imply that Wald's lemma (see, e.g., \cite[Chapter 1,
  Theorem~3.2]{Brem99}) is applicable to (\ref{add-eqn-52}). We thus
have
\begin{equation}
\EE[A_m+I_m^\on] 
= \EE[N_{m,0}] \EE[I^\on]
= \lambda \EE[I^\on],\qquad m\in\bbZ.
\label{add-eqn-53}
\end{equation}
Substituting (\ref{add-eqn-53}) into (\ref{asympt-P_{loss}}) yields
(\ref{asympt-P_{loss}-02}). \qed

\begin{rem}
If $I_m^\on-1$ is a stopping time
with respect to $\{N_{m,n};n=0,1,\dots\}$, then 
condition (iv) of Theorem~\ref{thm-app1-depend-01} is satisfied. 
\end{rem}

\begin{thm}\label{thm-app1-depend-02}
Suppose that Assumption~\ref{assumpt-alpha(k)-D(k)} holds and $\EE[A]
< \EE[I^\off]$. Further suppose that (i) $I^\on \in \calC$; (ii)
$\EE[I^\on] < \infty$; and (iii) $x\PP(Y > x) = o(\PP(I^\on>x))$. We
then have (\ref{asympt-P_{loss}}).  In addition, if (iv) each for $m
\in \bbZ$, $\{I_m^\on \ge n+1\}$ is independent of $N_{m,n}$ for all
$n=0,1,\dots$; and (v) $\vc{\alpha}(k) = \vc{\phi}\vc{\varLambda}(k)$
for $k=1,2,\dots$, then (\ref{asympt-P_{loss}-02}) holds.
\end{thm}

\proof Suppose that the sampling time $T$ and the cumulative process
$\{B(t)\}$ in (\ref{defn-B(t)}) satisfy the conditions of
Theorem~\ref{thm-dependent-02}. Applying
Theorem~\ref{thm-dependent-02} to (\ref{eqn-A}), we have
(\ref{asymp-PP(A>k)}) and thus (\ref{asymp-PP(A_e>k)}).  Note here
that (\ref{asymp-PP(A>k)}) and $I^\on \in \calC$ imply $A \in \calC
\subset \calS^{\ast}$, which leads to $A_{\re} \in \calS$ (see
Remarks~\ref{rem-S^*} and \ref{rem-class-D}). Therefore we can prove
(\ref{asympt-P_{loss}}) and (\ref{asympt-P_{loss}-02}) by following
the proof of Theorem~\ref{thm-app1-depend-01} (and using
Theorem~\ref{thm-dependent-02} instead of
Theorem~\ref{thm-dependent-01}).

In what follows, we confirm that $T$ and $\{B(t)\}$ in
(\ref{defn-B(t)}) satisfy conditions (ii)--(v) of
Theorem~\ref{thm-dependent-02} (condition (i) is obvious due to
$T\stackrel{\rd}{=}I^\on-1$ and $I^\on \in \calC$). Since the
distribution of $\Delta \tau_n$ ($n=0,1$) is phase-type,
\begin{equation}
\EE[(\Delta \tau_n)^p] < \infty, \qquad \forall p > 0,
\label{finite-moment-tau}
\end{equation}
which implies condition (ii) of Theorem~\ref{thm-dependent-02}.
Further since $T \in \calC \subset \calD$,
Proposition~\ref{prop-Bing89} shows that $\PP(T > x) = O(x^{-\gamma})$
for some $\gamma > 0$.  From this and (\ref{finite-moment-tau}), we
have
\begin{equation}
\limsup_{x\to\infty}
{x\PP(\Delta \tau_n > x) \over \PP(T > x)}
\le C\limsup_{x\to\infty}x^{\gamma +1}\PP(\Delta \tau_n > x)
= 0.
\label{add-3rd-01}
\end{equation}
Note here that (\ref{add-ineqn-06}) holds due to
Assumption~\ref{assumpt-alpha(k)-D(k)} and
Lemma~\ref{add-lem-MAdP-01}. It then follows from condition (iii) of
Theorem~\ref{thm-app1-depend-02} that for $n=0,1$,
\begin{equation}
x\PP(\Delta B_n > x) = O(x\PP(Y > x)) = o(\PP(T > x)).
\label{add-3rd-02}
\end{equation}
Note also that $\PP(-B(0) > x > 0) = 0$ due to $B(0) \ge 0$.
Therefore (\ref{add-3rd-01}) and (\ref{add-3rd-02}) imply that
condition (iii) of Theorem~\ref{thm-dependent-02} is satisfied. In
addition,
\begin{eqnarray*}
x\PP(|\Delta B_1 - \Delta \tau_1| > x) 
&\le& x[\PP(\Delta B_1 > x) + \PP(\Delta \tau_1 > x)]
= o(\PP(T > x)),
\end{eqnarray*}
which shows that condition (iv) of Theorem~\ref{thm-dependent-02} is
satisfied. Finally, condition (v.b) of Theorem~\ref{thm-dependent-02}
holds due to $\EE[T] < \infty$ .  \qed

\begin{thm}\label{thm-app1-independ-01}
Suppose that Assumption~\ref{assumpt-alpha(k)-D(k)} holds, $\EE[A] <
\EE[I^\off]$ and $I_m^\on$ is independent of the $m$th BMAP for all
$m\in\bbZ$.  Further suppose that (i) $I^\on \in \calL^{1/\theta}$ for
some $0 < \theta \le 1/2$; (ii) $\EE[I^\on] < \infty$ and $I^\on_{\re}
\in \calS$; and (iii) $\EE[\exp\{Q(Y)\}] < \infty$ for some $Q \in
\SC$ such that $x^{\theta} = O(Q(x))$. Under these conditions,
(\ref{asympt-P_{loss}}) holds. In addition, if $\vc{\alpha}(k) =
\vc{\phi}\vc{\varLambda}(k)$ for $k=1,2,\dots$, then
(\ref{asympt-P_{loss}-02}) holds.
\end{thm}

\proof According to the proofs of Theorems~\ref{thm-app1-depend-01}
and \ref{thm-app1-depend-02}, it suffices to show that $T$ and
$\{B(t)\}$ in (\ref{defn-B(t)}) satisfy conditions (i)--(iii) of
Theorem~\ref{thm-independ-01}.

Since $T \stackrel{\rd}{=} I^\on - 1$, condition (i) of
Theorem~\ref{thm-app1-independ-01} implies condition (i) of
Theorem~\ref{thm-independ-01}. Further since $\{B(t)\}$ in
(\ref{defn-B(t)}) is nondecreasing with $t$, we have $\Delta
B_n^{\ast} = \Delta B_n \ge 0$ ($n=0,1$). It thus follows from
(\ref{add-ineqn-06}) and condition (iii) of
Theorem~\ref{thm-app1-independ-01} that $\EE[\exp\{Q(\Delta
  B_n^{\ast})\}] < \infty$ ($n=0,1$) for some $Q \in \SC$ such that
$x^{\theta} = O(Q(x))$, which implies condition (iii) of
Theorem~\ref{thm-independ-01}. Note here that $\EE[\exp\{Q(\Delta
  B_n^{\ast})\}] < \infty$ ($n=0,1$) leads to $\EE[(\Delta B_1)^2] <
\infty$ (see Remark~\ref{remark-SC'-02}). Note also that $\EE[(\Delta
  \tau_1)^2] < \infty$ due to (\ref{finite-moment-tau}).  Therefore
condition (ii) of Theorem~\ref{thm-independ-01} are satisfied. \qed

\begin{thm}\label{thm-app1-independ-02}
Suppose that Assumption~\ref{assumpt-alpha(k)-D(k)} holds, $\EE[A] <
\EE[I^\off]$ and $I_m^\on$ is independent of the $m$th BMAP for all
$m\in\bbZ$.  Further suppose that (i) $I^\on \in \calC$; (ii)
$\EE[I^\on] < \infty$; and (iii) $\PP(Y > x) = o(\PP(I^\on >
x))$. Under these conditions, (\ref{asympt-P_{loss}}) holds.  In
addition, if $\vc{\alpha}(k) = \vc{\phi}\vc{\varLambda}(k)$ for
$k=1,2,\dots$, then (\ref{asympt-P_{loss}-02}) holds.
\end{thm}

\proof It is easy to see that the conditions of
Theorem~\ref{thm-MAdP-01} are satisfied. Thus similarly to the other
theorems in this subsection, we can prove (\ref{asympt-P_{loss}}) and
(\ref{asympt-P_{loss}-02}). \qed

\smallskip

Finally, we mention previous studies related to the results presented
in this subsection. Zwart~\cite{Zwar00} and Jelenkovi\'{c} and
Mom\v{c}ilovi\'{c}~\cite{Jele03b} study the subexponential asymptotics
of the loss probabilities of finite-buffer fluid queues fed by the
superposition of independent on/off sources that generate fluid at
constant rates. These studies assume that the lengths of the
on-periods of each on/off source follow a regularly or consistently
varying distribution, and then present asymptotic formulas for the
loss probability such that the decay of the loss probability is
connected to the tail of the equilibrium distribution of on-period
lengths.

\appendix

\section{Technical Lemmas}\label{sec-preliminary-lem}

This appendix presents technical lemmas, whose proofs are all given
in Appendix~\ref{sec-proof-lemmas}.

\subsection{Higher-order long-tailed distributions}\label{higher-order}

In this section, we consider the class $\calL^p$ ($p\ge1$) of
higher-order long-tailed distributions. By definition, $\calL^1 =
\calL$ (see Definition~\ref{defn-subclass-L}). Further $\calL^2$ is
equivalent to the class of square-root insensitive distributions (see
Lemma~1 in \cite{Jele04}). We can readily confirm that the following
are examples of the distributions in $\calL^p$:

\begin{enumerate}
\item $\PP(X > x) \simhm{x} \exp\{-x^{\alpha}\}$, where $0 <
  \alpha < 1/p$.
\item $\PP(X > x) \simhm{x} \exp\{-x^{1/p}/(\log x)^{\gamma}\}$, where
  $\gamma > 0$.
\end{enumerate}

In what follows, we provide five lemmas, which summarize the basic
properties of $\calL^p$.
\begin{lem}\label{lem-01}
If $X \in \calL^{1/\theta}$ (i.e., $X^{\theta} \in \calL$) for
some $0 < \theta \le 1$, the following are satisfied:

\begin{enumerate}
\item\label{lem-01-a} $\lim_{x \to \infty}e^{\varepsilon
  x^{\theta}}\PP(X > x) = \infty$ for any $\varepsilon > 0$, i.e.,
  $\PP(X > x)=e^{-o(x^{\theta})}$.
\item \label{lem-01-c} $X \in \calL^{1/\eta}$ for all $1 \le 1/\eta <
  1/\theta$.
\end{enumerate}
\end{lem}

\proof See Appendix~\ref{proof-lem-01}. \qed

\begin{rem}\label{rem-property-higher-order-L}
Lemma~\ref{lem-01}~(ii) implies that $\calL^{p_2} \subset
\calL^{p_1}$ for $1 \le p_1 < p_2$.
\end{rem}

\begin{lem}\label{lem-02}
For any $0 < \theta \le 1$, $X \in \calL^{1/\theta}$ if and only if
$\PP(X > x - \xi x^{1-\theta}) \simhm{x} \PP(X > x)$ for all $\xi \in
\bbR$.
\end{lem}

\proof See Appendix~\ref{proof-lem-02}. \qed

\smallskip

Lemma~\ref{lem-02} is an extension of Lemma~1 in \cite{Jele04}.  The
following lemma shows that the ``if" part of Lemma~\ref{lem-02} holds
under a weaker condition.
\begin{lem}\label{add-lem-02}
For any $0 < \theta \le 1$, $X \in \calL^{1/\theta}$ if $\PP(X > x -
\xi x^{1-\theta}) \simhm{x} \PP(X > x)$ for some $\xi \in \bbR
\backslash \{0\}$.
\end{lem}

\proof See Appendix~\ref{proof-add-lem-02}. \qed

\begin{rem}
The statements of Lemmas~\ref{lem-01}--\ref{add-lem-02} are presented
in a slightly different way in a technical report \cite{Miyo11} (see
Lemmas~1--3 therein), where the statements are described in terms of
$h$-insensitivity (see Chapter~2 in \cite{Foss11} for the definition
of $h$-insensitivity).
\end{rem}

Lemma~\ref{lem-C-in-L^{infty}} below shows the inclusion relation
between class $\calL^p$ and the consistent variation class $\calC$.
\begin{lem}\label{lem-C-in-L^{infty}}
$\calC \subset \calL^{\infty}$, i.e., $\calC \subset \calL^{1/\theta}$
  for any $0 < \theta \le 1$.
\end{lem}

\proof See Appendix~\ref{proof-lem-C-in-L^{infty}}. \qed

\smallskip

The following lemma is used to prove Theorem~\ref{thm-independ-01}.
\begin{lem}\label{lem-03}
If $X \in \calL^{1/\theta}$ for some $0 < \theta \le 1$, then for any
$\varepsilon > 0$ there exists $\breve{x}_{\varepsilon} > 0$ such that
for all $x > \breve{x}_{\varepsilon}$ and $0 \le u \le g(x)$,
\[
\PP(X > x - u) \le \PP(X>x)e^{\varepsilon (u^{\theta}+1)},
\]
where $g$ is a nonnegative function on $[0,\infty)$ such that
  $\limsup_{x\to\infty}g(x)/x < 1$.
\end{lem}

\proof See Appendix~\ref{proof-lem-03}. \qed

\subsection{Subexponential concave distributions}\label{appendix-SC}

The subexponential concave class was first introduced by
Nagaev~\cite{Naga77}. According to Nagaev's definition of $\SC$,
condition (iii) of Definition~\ref{defn-SC'} is replaced by the
following condition: ${\rm (iii^{\prime})}$ there exist $x_0 > 0$, $0
< \alpha < 1$ and $0 < \beta < 1$ such that for all $x \ge x_0$ and
$\beta x \le u \le x$,
\begin{equation}
{Q_X(x)-Q_X(u) \over Q_X(x)} \le \alpha{x-u \over x}.
\label{ineqn-Q_X-02}
\end{equation}
Actually, Nagaev's definition is equivalent to
Definition~\ref{defn-SC'}. Lemma~3.1 (i) in \cite{Jele03a} shows that
Nagaev's definition implies Definition~\ref{defn-SC'}. The converse
follows from Theorem~2 in \cite{Shne06}, though the phrase
``$Q(x)/f(x)$ is nondecreasing" should be replaced by ``$Q(x)/f(x)$ is
nonincreasing."

\begin{rem}\label{rem-SC-differentiable}
Suppose that $Q \in \SC$ is differentiable. It then follows from
(\ref{ineqn-Q_X-02}) and (\ref{ineqn-Q_X-03}) that
\begin{equation}
Q'(x) := {\rd \over \rd x}Q(x) \le {\alpha Q(x) \over x}
\le Cx^{\alpha-1},\qquad x > x_0.
\label{ineqn-Q'(x)}
\end{equation}
\end{rem}

Lemma~\ref{add-lem-SC} below establishes the relationship between
class $\SC$ and the higher-order long-tailed class.
\begin{lem}\label{add-lem-SC}
\begin{enumerate}
\item $\SC_{\alpha} \subset \calL^{1/\beta}$ for all $0 < \alpha <
  \beta \le 1$.
\item $X^{\alpha} \in \calL$ if $X \in \SC_{\alpha}$ for some $0 <
  \alpha < 1$ and
\begin{equation}
\lim_{x\to\infty}Q_X(x)/x^{\alpha} = 0.
\label{add-eqn-18}
\end{equation}
\end{enumerate}
\end{lem}

\proof See Appendix~\ref{proof-add-lem-SC}. \qed 

\smallskip

The following lemma plays an important role in the proof of
Theorems~\ref{thm-dependent-01} and \ref{thm-independ-01}.
\begin{lem}\label{prop-LD-B(t)}
Assume $\EE[(\Delta B_1)^2] < \infty$.

\begin{enumerate}
\item If $\EE[(\Delta \tau_1)^2] < \infty$ and $\EE[\exp\{Q(\Delta
  B_n^{\ast})\}] < \infty$ ($n=0,1$) for some $Q \in \SC$, then
\[
\PP\left(\sup_{0\le t \le x}\{B(t) - bt\} > u \right)
\le C\left(e^{-cu^2/x} + e^{-cx} + xe^{-cQ(u)}\right),
\quad \forall x \ge 0,~\forall u \ge 0.
\]
\item Let $\Delta B_0^{\ast\ast} = \sup_{0 \le t \le
  \tau_0}\max(-B(t),0)$ and $\Delta B_n^{\ast\ast}= \sup_{\tau_{n-1}
  \le t \le \tau_n}(B(\tau_{n-1})-B(t))$ for $n=1,2,\dots$. If
  $\EE[\exp\{Q(\Delta B_n^{\ast\ast}+b\Delta \tau_n)\}] < \infty$
  ($n=0,1$) for some $Q \in \SC$, then
\[
\PP\left(\inf_{0\le t \le x}\{B(t) - bt\} < -u \right)
\le C\left(e^{-cu^2/x} + e^{-cx} + xe^{-cQ(u)}\right),
\quad \forall x \ge 0,~\forall u \ge 0.
\]
\end{enumerate}
In the two above inequalities, $C$ and $c$ are independent of $x$ and
$u$.
\end{lem}

\proof See Appendix~\ref{proof-prop-LD-B(t)}. \qed

\begin{rem}\label{add-remark}
Lemma~\ref{prop-LD-B(t)}~(i) is a slight extension of Proposition~1 in
\cite{Jele04}, where the latter assumes that $\Delta B_1 \ge 0$ w.p.1.
\end{rem}

\begin{rem}\label{add-remark2}
Suppose that $\{B(t)\}$ is nondecreasing with $t$.  It then follows
that $\Delta B_0^{\ast} = (\Delta B_0)^+$, $\Delta B_1^{\ast} = \Delta
B_1$, $\Delta B_0^{\ast\ast} = \max(-B(0),0)=(-B(0))^+$ and $\Delta
B_1^{\ast\ast} = 0$. Therefore the condition $\EE[\exp\{Q(\Delta
  B_n^{\ast})\}] < \infty$ ($n=0,1$) is reduced to
\[
\EE[\exp\{Q( (\Delta B_0)^+ )\}] < \infty \quad
\mbox{and} \quad 
\EE[\exp\{Q(\Delta B_1)\}] < \infty;
\]
and the condition $\EE[\exp\{Q(\Delta
  B_n^{\ast\ast}+b\Delta \tau_n)\}] < \infty$ ($n=0,1$) is reduced to
\[
\EE[\exp\{Q( (-B(0))^+ + b\Delta \tau_0 )\}] < \infty \quad
\mbox{and} \quad 
\EE[\exp\{Q(b\Delta \tau_1)\}] < \infty.
\]
It should be noted that $\EE[\exp\{Q(\Delta B_1)\}] < \infty$ and
$\EE[\exp\{Q(b\Delta \tau_1)\}] < \infty$ imply $\EE[(\Delta B_1)^2] <
\infty$ and $\EE[(\Delta \tau_1)^2] < \infty$, respectively (see
Remark~\ref{remark-SC'-02}).
\end{rem}

\subsection{Regular varying distributions}\label{subsec-class-R}

The regular variation class $\calR$ is defined as follows:
\begin{defn}\label{defn-class-R}
A nonnegative r.v.\ $X$ and its d.f.\ $F_X$ belong to class
$\calR(-\alpha)$ ($\alpha \ge 0$) if $\overline{F}_X$ is regularly
varying with index $-\alpha$, i.e.,
\[
\lim_{x\to\infty}{\overline{F}_X(vx) \over \overline{F}_X(x)}
= v^{-\alpha}, \qquad \forall v > 0.
\]
Further let $\calR = \cup_{\alpha \ge 0}\calR(-\alpha)$.
\end{defn}

\begin{rem}
If $F \in \calR$, then $\overline{F}(x) = x^{-\alpha}\tilde{l}(x)$ for some
$\alpha \ge 0$, where $\tilde{l}$ is a slowly varying function, i.e.,
\[
\lim_{x\to\infty}{\tilde{l}(vx) \over \tilde{l}(x)}
= 1, \qquad \forall v > 0.
\]
See \cite{Bing89} for the details of regularly varying functions.
\end{rem}

\begin{rem}\label{rem-calL^{infty}}
It is known that $\calR \subset \calC$ (see, e.g.,
\cite{Clin94,Embr84}). Thus $\calR \subset \calC \subset
\calL^{\infty} \subset \calL^p \subset \calL$ for any $p > 1$ (see
Remark~\ref{rem-property-higher-order-L} and
Lemma~\ref{lem-C-in-L^{infty}}).
\end{rem}

The following lemma is used to prove Theorem~\ref{thm-dependent-02}.
\begin{lem}\label{lem-upper-calC}
Suppose that $U$ is a r.v.\ with $\EE[|U|] < \infty$. If $\PP(U > x) =
o(\PP(Y > x))$ for some nonnegative r.v.\ $Y$ with $\EE[Y] < \infty$,
then for any $\mu > \EE[U]$ there exists some r.v.\ $Z$ in $\bbR$ such
that $\EE[U] < \EE[Z] < \mu$, $\overline{F}_Z(x) \ge
\overline{F}_U(x)$ for all $x \in \bbR$ and
\[
\overline{F}_Z(x) = \tilde{l}(x)\overline{F}_Y(x)
\qquad \mbox{for all sufficiently large $x > 0$},
\]
where $\tilde{l}$ is some slowly varying function such that
$\lim_{x\to\infty}\tilde{l}(x) = 0$.
\end{lem}

\proof See Appendix~\ref{proof-lem-upper-calC}. \qed

\subsection{Bounds on deviation probabilities}

This subsection presents three lemmas on the deviation probabilities
associated with i.i.d.\ r.v.s.  The first one
(Lemma~\ref{add-prop-N(x)}) is used to prove
Theorems~\ref{thm-dependent-02} and \ref{thm-independ-02}, and the
other two are required by the proof of Theorem~\ref{thm-dependent-02}.
\begin{lem}\label{add-prop-N(x)}
If $X,X_1,X_2,\dots$ are i.i.d.\ nonnegative
r.v.s with $\EE[X] > 0$ and $\EE[X^2] < \infty$, then for any
$\delta > 0$ there exist finite constants
$\widetilde{C}:=\widetilde{C}(\delta) > 0$ and
$\tilde{c}:=\tilde{c}(\delta) > 0$ such that
\begin{equation}
\PP\left(N_X(x) -  {x \over \EE[X]} > u \right) 
\le \widetilde{C} \exp\{-\tilde{c}u^2/x\},
\quad
\forall x \ge 0,~~0 \le \forall u \le \delta x,
\label{eqn-prob-N(x)}
\end{equation}
where $N_X(x) =
\max\{k\ge0;\sum_{n=1}^k X_n \le x\}$ for $x \in \bbR$.
\end{lem}

\proof See Appendix~\ref{proof-add-prop-N(x)}. \qed

\begin{rem}
Lemma~\ref{add-prop-N(x)} is very similar to, but not exactly the same
as Lemma~6 in \cite{Jele04}. The latter states that there exists {\it
  some} $\delta > 0$ such that (\ref{eqn-prob-N(x)}) holds.
\end{rem}

Lemmas~\ref{lem-LD-maxima} and \ref{prop-Lin11} below are extensions
of Lemma~2.3 in \cite{Tang06} and Lemma~2.2 in \cite{Lin11},
respectively, to the maxima of partial sums of i.i.d.\ r.v.s.
\begin{lem}\label{lem-LD-maxima}
Suppose that $U,U_1,U_2,\dots$ are i.i.d.\ r.v.s in $\bbR$. If $\EE[U]
= 0$ and $\EE[(U^+)^r] < \infty$ for some $r > 1$, then for any fixed
$\gamma > 0$ and $p > 0$ there exist some $v:=v(r,p) > 0$ and
$\widetilde{C}:=\widetilde{C}(v,\gamma)$ such that for all
$n=1,2,\dots$,
\begin{equation}
\PP\left(\max_{1\le k\le n}\sum_{i=1}^k U_i \ge x \right)
  \le n\PP(U > v x) + \widetilde{C}x^{-p},
\quad \forall x \ge \gamma n.
\label{add-eqn-35}
\end{equation}
\end{lem}

\proof See Appendix~\ref{proof-lem-LD-maxima}. \qed 

\begin{lem}\label{prop-Lin11}
Suppose that $U,U_1,U_2,\dots$ are i.i.d.\ r.v.s in $\bbR$. If $0 \le
\EE[U] < \infty$ and $U^+ \in \calC$, then for any $\gamma > \EE[U]$
there exists some constant $\widetilde{C}:=\widetilde{C}(\gamma) > 0$
such that for all $n=1,2,\dots$,
\begin{equation}
\PP\left(\max_{1\le k \le n}\sum_{i=1}^k U_i > x\right) 
\le \widetilde{C} n\PP(U > x),
\quad \forall x \ge \gamma n.
\label{add-eqn-33}
\end{equation}
\end{lem}

\proof See Appendix~\ref{proof-prop-Lin11}. \qed

\subsection{Convolution tail of matrix-valued functions associated 
with subexponential distributions}\label{subsec-convo-tail}

Let $\vc{F}=(F_{i,j})\genfrac{}{}{0pt}{1}{1 \le i \le m_0}{1 \le j \le
  m_1}$ and $\vc{G}=(G_{i,j})\genfrac{}{}{0pt}{1}{1 \le i \le m_1}{1
  \le j \le m_2}$ denote matrix-valued functions on $\bbR$. Assume
that $F_{i,j}(x)$ and $G_{i,j}(x)$ are nonnegative and nondecreasing
for all $x \in \bbR$ and that
$F_{i,j}(\infty):=\lim_{x\to\infty}F_{i,j}(x) < \infty$ and
$G_{i,j}(\infty):=\lim_{x\to\infty}G_{i,j}(x) < \infty$. We then
define $\overline{\vc{F}}(x) = \vc{F}(\infty) - \vc{F}(x)$ and
$\overline{\vc{G}}(x) = \vc{G}(\infty) - \vc{G}(x)$ for $x \in \bbR$.

Let $\vc{F} \ast \vc{G}$ denote the convolution of
$\vc{F}$ and $\vc{G}$, i.e.,
\[
\vc{F} \ast \vc{G}(x) 
= \int_{y \in \bbR} \vc{F}(x-y)\rd \vc{G}(y)
= \int_{y \in \bbR} \rd\vc{F}(y) \vc{G}(x-y),
\qquad x \in \bbR.
\]
Let $\overline{\vc{F} \ast \vc{G}}(x)$ ($x \in \bbR$) denote
\[
\overline{\vc{F} \ast \vc{G}}(x)
= \vc{F} \ast \vc{G}(\infty) - \vc{F} \ast \vc{G}(x)
= \vc{F}(\infty)\vc{G}(\infty) - \vc{F} \ast \vc{G}(x).
\]
When $\vc{F}$ is a square matrix-valued function (i.e., $m_0=m_1$), we
define $\vc{F}^{*n}$ ($n =1,2,\dots$) as the $n$-fold convolution of
$\vc{F}$ itself, i.e.,
\[
\vc{F}^{\ast n}(x) = \vc{F}^{\ast(n-1)} \ast \vc{F}(x),\quad x \in \bbR,
\]
and for convenience, define $\vc{F}^{*0}(x) = \vc{O}$ for $x < 0$ and
$\vc{F}^{*0}(x) = \vc{I}$ for $x \ge 0$.  Further for the $n$-fold
convolution $\vc{F}^{\ast n}$, let $\overline{\vc{F}^{\ast n}}(x)$ ($x
\in \bbR$) denote
\[
\overline{\vc{F}^{\ast n}}(x)
= \vc{F}^{\ast n}(\infty) - \vc{F}^{\ast n}(x)
= (\vc{F}(\infty))^n  - \vc{F}^{\ast n}(x).
\]

The following lemma is the upper-limit version of
Proposition A.3 in \cite{Masu11} and Lemma~6 in~\cite{Jele98}.
\begin{lem}\label{lem-convo-tail}
Suppose that for some r.v.\ $Y \in \calS$,
\begin{equation}
\limsup_{x\to\infty}{\overline{\vc{F}}(x) \over \PP(Y > x)} 
\le \widetilde{\vc{F}},
\qquad
\limsup_{x\to\infty}{\overline{\vc{G}}(x) \over \PP(Y > x)} 
\le \widetilde{\vc{G}},
\label{eqn-05}
\end{equation}
where $\widetilde{\vc{F}}=(\widetilde{F}_{i,j})$ and
$\widetilde{\vc{G}}=(\widetilde{G}_{i,j})$ are finite, and where
$\widetilde{\vc{F}} = \widetilde{\vc{G}} = \vc{O}$ is allowed.  We
then have
\begin{equation}
\limsup_{x\to\infty} {\overline{\vc{F} \ast \vc{G}}(x) \over \PP(Y > x)} 
\le \widetilde{\vc{F}} \vc{G}(\infty)
+ \vc{F}(\infty) \widetilde{\vc{G}}.
\label{eqn-05a}
\end{equation}
Further if $\vc{F}$ is a square matrix-valued function, then
\begin{equation}
\limsup_{x\to\infty} {\overline{\vc{F}^{\ast n}}(x) \over \PP(Y > x)} 
\le \sum_{\nu=0}^{n-1} (\vc{F}(\infty))^{\nu} \widetilde{\vc{F}}
(\vc{F}(\infty))^{n-\nu-1}.
\label{eqn-05b}
\end{equation}
In addition to the above conditions, assume that $\sum_{n=0}^{\infty}
(\vc{F}(\infty))^n=(\vc{I} - \vc{F}(\infty) )^{-1} < \infty$. We then
have
\begin{equation}
\limsup_{x\to\infty} \sum_{n=0}^{\infty}
{\overline{\vc{F}^{\ast n}}(x) \over \PP(Y > x)} 
\le (\vc{I} -  \vc{F}(\infty) )^{-1} \widetilde{\vc{F}}
(\vc{I} -  \vc{F}(\infty) )^{-1}.
\label{eqn-05c}
\end{equation}
\end{lem}

\proof See Appendix~\ref{proof-lem-convo-tail}. \qed

\section{Proofs of Technical Lemmas}\label{sec-proof-lemmas}

\subsection{Proof of Lemma \ref{lem-01}}\label{proof-lem-01}

We first
prove the statement (i). It follows from $X^{\theta} \in \calL$ that
$\lim_{y\to\infty}e^{\varepsilon y}\PP(X^{\theta}>y) = \infty$ for any
$\varepsilon > 0$. Thus letting $x=y^{1/\theta}$ for $y > 0$, we have
\[
\lim_{x\to\infty}e^{\varepsilon x^{\theta}}\PP(X > x)
= \lim_{x\to\infty}e^{\varepsilon x^{\theta}}\PP(X^{\theta} > x^{\theta})
= \lim_{y\to\infty}e^{\varepsilon y}\PP(X^{\theta}>y)
= \infty.
\]

Next we prove the statement (ii). For all $x,y \ge 0$, we have
\begin{equation}
1 \ge {\PP(X^{\eta} > x+y) \over \PP(X^{\eta} > x)} 
= {\PP(X^{\theta} > (x+y)^{\theta/\eta}) 
\over \PP(X^{\theta} > x^{\theta/\eta})}.
\label{lem-01-ineqn-00}
\end{equation}
It follows from $0 < \theta/\eta < 1$ that for all $x,y \ge 0$,
\[
(x+y)^{\theta/\eta} \le x^{\theta/\eta} + y^{\theta/\eta},
\]
from which and $X^{\theta} \in \calL$ we obtain
\begin{equation}
\lim_{x \to \infty}
{\PP(X^{\theta} > (x+y)^{\theta/\eta}) 
\over \PP(X^{\theta} > x^{\theta/\eta})}
\ge \lim_{x \to \infty}
{\PP(X^{\theta} > x^{\theta/\eta}  + y^{\theta/\eta})
\over \PP(X^{\theta} > x^{\theta/\eta})}
= 1.
\label{lem-01-ineqn-01}
\end{equation}
Combining (\ref{lem-01-ineqn-01}) with (\ref{lem-01-ineqn-00}) yields
$\PP(X^{\eta} > x+y) \simhm{x} \PP(X^{\eta} > x)$ for any $y > 0$,
i.e., $X^{\eta} \in \calL$.

\subsection{Proof of Lemma \ref{lem-02}}\label{proof-lem-02}

To prove Lemma~\ref{lem-02}, we use the following proposition, whose
proof is given in Appendix~\ref{proof-prop-00}.
\begin{prop}\label{prop-00}
For any $\gamma > 0$ and $x > y \ge 0$, 
\begin{eqnarray}
(x+y)^{\gamma} 
&\le& x^{\gamma} + C\left(1-{y \over x}\right)^{-1} y x^{\gamma - 1},
\label{add-lem00-01}
\\
(x-y)^{\gamma} 
&\ge& x^{\gamma} - C\left(1-{y \over x}\right)^{-1} y x^{\gamma - 1},
\label{add-lem00-02}
\end{eqnarray}
where $C$ is independent of $x$ and $y$.
\end{prop}

We first prove the ``if" part of
Lemma~\ref{lem-01}. Proposition~\ref{prop-00} implies that
$(x+y)^{1/\theta} \le x^{1/\theta} + C y x^{{1/\theta}-1}$ for any $x
> 0$ and $0 \le y < x/2$.  Thus for any $y > 0$,
\begin{eqnarray*}
1 \ge 
\lim_{x\to\infty}
{\PP(X^{\theta} > x+y) \over \PP(X^{\theta} > x)}
&\ge& \lim_{x\to\infty} 
{\PP(X > x^{1/\theta} + C y \cdot (x^{1/\theta})^{1-\theta} ) 
\over \PP(X > x^{1/\theta})}
\nonumber
\\
&=& \lim_{w\to\infty} 
{\PP(X > w + C y \cdot w^{1-\theta} ) 
\over \PP(X > w)}
= 1,
\end{eqnarray*}
which shows that $X^{\theta} \in \calL$.

Next we prove the ``only if" part. We fix $\xi$ such that $x^{\theta}
> 2|\xi|$. It then follows from Proposition~\ref{prop-00} that
\begin{align*}
&&
(x - \xi x^{1-\theta})^{\theta} 
&\ge x^{\theta} - C\left(1 - {\xi \over x^{\theta}} \right)^{-1}
\xi 
\ge x^{\theta} - 2C\xi, 
& \xi \ge 0, & &
\\
&&
(x - \xi x^{1-\theta})^{\theta} 
&\le x^{\theta} + C\left(1 - {-\xi \over x^{\theta}} \right)^{-1}
(-\xi) 
\le x^{\theta} + 2C(-\xi),
& \xi < 0. & &
\end{align*}
Thus for $\xi \ge 0$, 
\begin{eqnarray*}
1 
&\le& \lim_{x\to\infty}
{\PP(X > x - \xi x^{1-\theta}) \over \PP(X > x)}
\le \lim_{x\to\infty}
{\PP(X^{\theta} > x^{\theta} - C\xi) \over \PP(X^{\theta} > x^{\theta})} 
= 1,
\end{eqnarray*}
where the last equality follows from $X^{\theta} \in
\calL$. Similarly, for $\xi < 0$, 
\[
1 \ge 
\lim_{x\to\infty}{\PP(X > x - \xi x^{1-\theta}) \over \PP(X > x)}
\ge \lim_{x\to\infty}
{\PP(X^{\theta} > x^{\theta} + C(-\xi)) 
\over \PP(X^{\theta} > x^{\theta})}
= 1.
\]
As a result, $\PP(X > x - \xi x^{1-\theta}) \simhm{x} \PP(X
  > x)$ for any $\xi \in \bbR$. 

\subsection{Proof of Lemma \ref{add-lem-02}}\label{proof-add-lem-02}

It follows from Proposition~\ref{prop-00} that there exists some
$\widetilde{C} > 0$ such that for all $x > 2\sigma > 0$,
\begin{eqnarray}
(x-\sigma)^{1/\theta} 
&\ge& x^{1/\theta} - \sigma \widetilde{C} \cdot (x^{1/\theta})^{1-\theta},
\label{add-eqn-21a}
\\
(x+\sigma)^{1/\theta} 
&\le&
x^{1/\theta} + \sigma \widetilde{C} \cdot (x^{1/\theta})^{1-\theta}.
\label{add-eqn-21b}
\end{eqnarray}
We fix $\sigma = |\xi|/\widetilde{C}$, where $\xi \in \bbR \setminus
\{0\}$.  Using (\ref{add-eqn-21a}), we then have
\begin{eqnarray}
1 \le
\lim_{x\to\infty}{ \PP(X^{\theta} > x-\sigma) \over \PP(X^{\theta} > x)}
&\le& \lim_{x\to\infty}
{ \PP(X > x^{1/\theta} 
- \sigma \widetilde{C} \cdot (x^{1/\theta})^{1-\theta})
\over \PP(X > x^{1/\theta})}
\nonumber
\\
&=& \lim_{x\to\infty}
{ \PP(X > x - |\xi| x^{1-\theta}) \over \PP(X > x)}.
\label{add-3rd-04}
\end{eqnarray}
Similarly from (\ref{add-eqn-21b}), we have
\begin{eqnarray}
1 \ge
\lim_{x\to\infty}{ \PP(X^{\theta} > x+\sigma) \over \PP(X^{\theta} > x)}
&\ge& \lim_{x\to\infty}
{ \PP(X > x + |\xi| x^{1-\theta}) \over \PP(X > x)}.
\label{add-3rd-05}
\end{eqnarray}

We now suppose that $\PP(X > x -\xi x^{1-\theta}) \simhm{x} \PP(X >
x)$ for some $\xi \in \bbR \setminus \{0\}$, which implies
\[
\PP(X > x - |\xi| x^{1-\theta}) \simhm{x} \PP(X > x)
~~\mbox{or}~~
\PP(X > x + |\xi| x^{1-\theta}) \simhm{x} \PP(X > x).
\]
It thus follows from (\ref{add-3rd-04}) and (\ref{add-3rd-05}) that
$\PP(X^{\theta} > x-\sigma) \simhm{x} \PP(X^{\theta} > x)$ or
$\PP(X^{\theta} > x+\sigma) \simhm{x} \PP(X^{\theta} > x)$, which
shows $X^{\theta} \in \calL$, i.e., $X \in \calL^{1/\theta}$.

\subsection{Proof of Lemma~\ref{lem-C-in-L^{infty}}}\label{proof-lem-C-in-L^{infty}}

Suppose $X \in \calC$. It then follows from
Definition~\ref{defn-class-C} that for any $v > 1$ there exists some
$c(v) > 0$ such that $\lim_{v\downarrow1}c(v) = 1$ and
\[
\liminf_{x\to\infty}{\PP(X>vx) \over \PP(X>x)}=c(v).
\]
Since $x+1 \le vx$ for any fixed $v>1$ and all sufficiently large $x >
0$, we have for any $0 < \theta \le 1$,
\[
\liminf_{x\to\infty}{\PP(X^{\theta}>x+1) \over \PP(X^{\theta}>x)}
\ge
\liminf_{x\to\infty}{\PP(X>(vx)^{1/\theta}) \over \PP(X>x^{1/\theta})}
=c(v^{1/\theta}) \to 1 \quad\mbox{as}~v\downarrow1.
\]
On the other hand, it is clear that $\PP(X^{\theta}>x+1) \lesimhm{x}
\PP(X^{\theta}>x)$. Therefore we obtain $\PP(X^{\theta}>x+1) \simhm{x}
\PP(X^{\theta}>x)$, i.e., $X^{\theta} \in \calL$.

\subsection{Proof of Lemma \ref{lem-03}}\label{proof-lem-03}

The case of $u = 0$ is obvious. Therefore we focus on the case of $u >
0$.  For any $x \ge u$, we have
\begin{equation}
{\PP(X > x - u) \over \PP(X > x)}
= {\PP(X^{\theta} > (x - u)^{\theta}) 
\over \PP(X^{\theta} > x^{\theta})}
\le
{\PP(X^{\theta} > x^{\theta} - u^{\theta}) 
\over \PP(X^{\theta} > x^{\theta})},
\label{add-lem3-ineqn-01}
\end{equation}
where we use $(x - u)^{\theta} \ge x^{\theta} - u^{\theta}$ for $0 \le
u \le x$.  Let $y$ denote a nonnegative number such that $y =
x^{\theta} - u^{\theta}$. We then have
\begin{equation}
{\PP(X^{\theta} > x^{\theta} - u^{\theta}) 
\over \PP(X^{\theta} > x^{\theta})}
= {\PP(X^{\theta} > y) 
\over \PP(X^{\theta} > y+u^{\theta})}
=   
\prod_{i=0}^{\lceil u^{\theta} \rceil-1}
{ 
\PP\left(
X^{\theta} >  y + i \dm{u^{\theta} \over \lceil u^{\theta} \rceil}
\right)
\over 
\PP\left(
X^{\theta} >  y + (i+1) \dm{u^{\theta} \over \lceil u^{\theta} \rceil}
\right) 
}. 
\qquad
\label{lem3-ineqn-01}
\end{equation}
It follows from $X^{\theta} \in \calL$ that for any $\varepsilon > 0$
there exists some $\breve{y}_{\varepsilon} > 0$ such that for all $y >
\breve{y}_{\varepsilon}$,
\[
{\PP(X^{\theta} > y) \over \PP(X^{\theta} > y+\gamma)}
\le
{\PP(X^{\theta} > y) \over \PP(X^{\theta} > y+1)}
\le e^{\varepsilon},
\qquad 0 \le \forall\gamma \le 1.
\]
Thus since $0 <
u^{\theta} / \lceil u^{\theta} \rceil \le 1$, we have
\[
\prod_{i=0}^{\lceil u^{\theta} \rceil-1}
{ 
\PP\left(
X^{\theta} >  y + i \dm{u^{\theta} \over \lceil u^{\theta} \rceil}
\right)
\over 
\PP\left(
X^{\theta} >  y + (i+1) \dm{u^{\theta} \over \lceil u^{\theta} \rceil}
\right)
}
\le e^{\varepsilon\lceil u^{\theta} \rceil} 
\le e^{\varepsilon (u^{\theta}+1)},
\qquad y > \breve{y}_{\varepsilon},
\]
from which, (\ref{add-lem3-ineqn-01}) and (\ref{lem3-ineqn-01}) it
follows that
\begin{equation}
{\PP(X > x - u) \over \PP(X > x)}
\le e^{\varepsilon (u^{\theta}+1)}
\quad \mbox{for all $x,u \ge 0$ such that $x^{\theta} - u^{\theta} >
\breve{y}_{\varepsilon}$}.
\label{add-eqn-11}
\end{equation}

Note here that for all $0 < u \le g(x)$,
\[
\liminf_{x\to\infty}(x^{\theta} - u^{\theta})
\ge \liminf_{x\to\infty}\left[x^{\theta} - \{g(x)\}^{\theta} \right] 
= \liminf_{x\to\infty}x^{\theta}
\left[1 - \left( g(x) \over x \right)^{\theta} \right] = \infty,
\]
where the last equality is due to $\limsup_{x\to\infty}g(x)/x<1$.
As a result, there exists some $\breve{x}_{\varepsilon} > 0$ such that
for all $x > \breve{x}_{\varepsilon}$ and $0 < u \le g(x)$
\[
x^{\theta} - u^{\theta} > \breve{y}_{\varepsilon},
\]
and thus (\ref{add-eqn-11}) holds.

\subsection{Proof of Lemma~\ref{add-lem-SC}}\label{proof-add-lem-SC}

For any $0 < \beta \le 1$, it follows from (\ref{ineqn-Q_X-02}) that
for all sufficiently large $x > 0$,
\begin{eqnarray}
1 \le {\PP(X > x - x^{1-\beta}) \over \PP(X > x)}
&=& \exp\{Q_X(x) - Q_X( x - x^{1-\beta})\}
\le \exp\{\alpha Q_X(x)/x^{\beta}\}. \quad
\label{add-ineqn-05}
\end{eqnarray}
Further according to condition (iii) of Definition~\ref{defn-SC'},
there exists some $x_0 > 0$ such that
\[
Q_X(x) \le Cx^{\alpha}, \qquad \forall x \ge x_0.
\]
Thus for any $\beta \in (\alpha,1]$, we have
\[
1 \le {\PP(X > x - x^{1-\beta}) \over \PP(X > x)}
\le \exp\{Cx^{\alpha-\beta}\} \to 1 \qquad \mbox{as}~x\to \infty,
\]
which implies $X^{\beta} \in \calL$ due to Lemma~\ref{add-lem-02}. In
addition, if (\ref{add-eqn-18}) holds, then substituting
(\ref{add-eqn-18}) into (\ref{add-ineqn-05}) with $\beta=\alpha$
yields $\PP(X > x - x^{1-\alpha}) \simhm{x} \PP(X>x)$, i.e.,
$X^{\alpha} \in \calL$.

\subsection{Proof of Lemma~\ref{prop-LD-B(t)}}\label{proof-prop-LD-B(t)}

To prove Lemma~\ref{prop-LD-B(t)}, we consider another (possibly
delayed) cumulative process $\{\breve{B}(t);t \ge 0\}$ on $\bbR$,
which satisfies $|\breve{B}(0)| < \infty$ w.p.1 and has the same
regenerative points as $\{B(t);t \ge 0\}$, i.e.,
$\{\breve{B}(t+\tau_n)-\breve{B}(\tau_n); t \ge 0\}$ ($n=0,1,\dots$)
is stochastically equivalent to
$\{\breve{B}(t+\tau_0)-\breve{B}(\tau_0); t \ge 0\}$ and is
independent of $\{\breve{B}(u);0 \le u < \tau_n\}$.  Let
\begin{eqnarray*}
\Delta \breve{B}_n 
&=& \left\{
\begin{array}{ll}
\breve{B}(\tau_0), & n = 0,
\\
\breve{B}(\tau_n)-\breve{B}(\tau_{n-1}), & n=1,2,\dots,
\end{array}
\right.
\\
\Delta \breve{B}_n^{\ast} 
&=& \left\{
\begin{array}{ll}
\dm\sup_{0 \le t \le \tau_0} \max(\breve{B}(t),0), 
& n=0,
\\
\dm\sup_{\tau_{n-1} \le t \le \tau_n} \breve{B}(t) - \breve{B}(\tau_{n-1}), 
& n=1,2,\dots.
\end{array}
\right.
\end{eqnarray*}
The $\Delta \breve{B}_n$'s (resp.\ $\Delta \breve{B}_n^{\ast}$'s)
($n=1,2,\dots$) are i.i.d.~and independent of $\Delta \breve{B}_0$
(resp.\ $\Delta \breve{B}_0^{\ast}$). We assume that
\[
\PP(0 \le \breve{B}_n^{\ast} < \infty) = 1~~(n=0,1),\quad
\EE[|\Delta \breve{B}_1|] <\infty,
\quad \breve{b} :=
{ \EE[\Delta \breve{B}_1] \big/ \EE[\Delta \tau_1]} \neq 0.
\]
Note that $\breve{b}$ can be negative.

The following lemma is an extension of Proposition~1 in
\cite{Jele04}. Using the lemma, we can readily prove
Lemma~\ref{prop-LD-B(t)}.
\begin{lem}\label{lemma-LD-B(t)}
Let $\Delta \Theta_n = \Delta \breve{B}_n^{\ast} -
\min(\breve{b},0)\Delta \tau_n \ge 0$ for $n=0,1,\dots$.  If
$\EE[(\Delta \breve{B}_1)^2]<\infty$, $\EE[(\Delta \tau_1)^2]<\infty$
and $\EE[\exp\{Q(\Delta \Theta_n)\}] < \infty$ ($n=0,1$) for some $Q
\in \SC$, then for all $x,u \ge 0$,
\begin{equation}
\PP\left(\sup_{0\le t \le x}\{\breve{B}(t) - \breve{b}t\} > u \right)
\le C\left(e^{-cu^2/x} + e^{-cx} + xe^{-cQ(u)}\right),
\label{prop-bound-01}
\end{equation}
where $C$ and $c$ are independent of $x$ and $u$.
\end{lem}

Note that if $\breve{B}(t)=B(t)$, then $\breve{b}=b > 0$, $\Delta
\Theta_n = \Delta B_n^{\ast}$ and
\[
\PP\left(\sup_{0\le t \le x}\{\breve{B}(t) - \breve{b}t\} > u \right)
=
\PP\left(\sup_{0\le t \le x}\{B(t) - bt\} > u \right).
\]
On the other hand, suppose that $\breve{B}(t)=-B(t)$. We then have
$\breve{b}=-b<0$ and
\begin{eqnarray*}
\Delta \Theta_0
&=& \sup_{0 \le t \le \tau_0}\max(-B(t),0) + b\Delta \tau_0
=: \Delta B_0^{\ast\ast}+ b\Delta \tau_0 \ge b\Delta \tau_0,
\\
\Delta \Theta_n 
&=& \sup_{\tau_{n-1} \le t \le \tau_n}(B(\tau_{n-1})-B(t)) + b\Delta \tau_n
=: \Delta B_n^{\ast\ast}+ b\Delta \tau_n \ge b\Delta \tau_n,
\quad n=1,2,\dots.
\end{eqnarray*}
Thus $\EE[\exp\{Q(\Delta\Theta_1)\}] < \infty$ implies $\EE[(\Delta
  \tau_1)^2]<\infty$ (see Remark~\ref{remark-SC'-02}). We also have
\[
\PP\left(\sup_{0\le t \le x}\{\breve{B}(t) - \breve{b}t\} > u \right)
=
\PP\left(\inf_{0\le t \le x}\{B(t) - bt\} < -u \right).
\]
As a result, Lemma~\ref{prop-LD-B(t)} follows from
Lemma~\ref{lemma-LD-B(t)}.  As for the proof of
Lemma~\ref{lemma-LD-B(t)}, see Appendix~\ref{proof-lem-LD-B(t)}.

\subsection{Proof of Lemma \ref{lem-upper-calC}}\label{proof-lem-upper-calC}

It follows from that Lemma~4.4 in \cite{Fay06} that there exists some
nonincreasing slowly varying function $l_0$ such that $l_0(0) = 1$,
$\lim_{x\to \infty}l_0(x) = 0$ and $\overline{F}_U(x) =
o(l_0(x)\overline{F}_Y(x))$. Thus for any $\varepsilon > 0$ and $x_0
\ge 0$, there exists some $x_2:=x_2(\varepsilon,x_0) > x_0$ such that
\begin{equation}
\overline{F}_U(x) < \varepsilon l_0(x)\overline{F}_Y(x)
\le {\overline{F}}_U(0),
\qquad \forall x \ge x_2.
\label{ineqn-overline{F}_U(x)}
\end{equation}

Let $x_1:=x_1(\varepsilon)$ denote
\[
x_1 = \inf \left\{
x \in [x_0,x_2]; \overline{F}_U(x) 
\le \varepsilon l_0(x_2)\overline{F}_Y(x_2)
\right\},
\]
from which and (\ref{ineqn-overline{F}_U(x)}) we obtain
$\overline{F}_U(x) \le \varepsilon l_0(x_2)\overline{F}_Y(x_2)$ for
all $x_1 < x \le x_2$. Further since $\overline{F}_U$ is
right-continuous,
\[
\overline{F}_U(x) \le \varepsilon l_0(x_2)\overline{F}_Y(x_2),
\qquad x_1 \le \forall x \le x_2.
\]

\begin{figure}[tp]
\centering
\includegraphics[scale=0.4,clip]{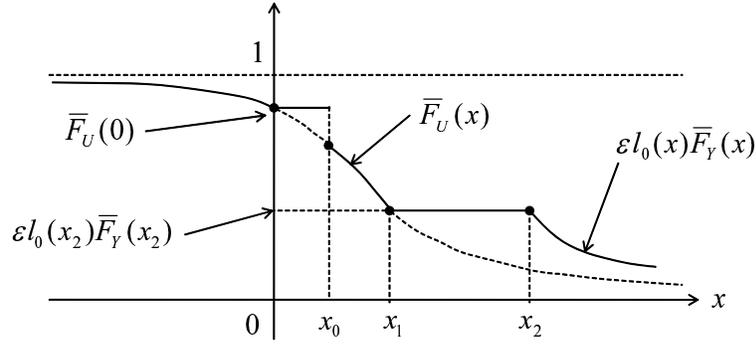}
\caption{Tail distribution of $Z(\varepsilon,x_0)$}
\label{fig-01}
\end{figure}

We now define $Z(\varepsilon,x_0)$ as a r.v.\ in $\bbR$ such that (see
Figure~\ref{fig-01})
\begin{eqnarray}
\overline{F}_{Z(\varepsilon,x_0)}(x)
&=& \left\{
\begin{array}{ll}
\overline{F}_U(x), & x < 0,
\\
\overline{F}_U(0), & 0 \le x < x_0,
\\
\overline{F}_U(x), & x_0 \le x < x_1,
\\
\varepsilon l_0(x_2) \overline{F}_Y(x_2), & x_1 \le x < x_2,
\\
\varepsilon l_0(x) \overline{F}_Y(x), & x \ge x_2.
\\
\end{array}
\right.
\label{defn-bar{F}_{X_varepsilon}}
\end{eqnarray}
Clearly, $\overline{F}_{Z(\varepsilon,x_0)}(x) \ge \overline{F}_U(x)$
for all $x \in \bbR$. Further it follows from
(\ref{defn-bar{F}_{X_varepsilon}}) that
\begin{eqnarray}
0 &\le& \EE[Z(\varepsilon,x_0)] - \EE[U]
=: S_1(\varepsilon,x_0) + S_2(\varepsilon,x_0) + S_3(\varepsilon,x_0),
\label{add-eqn-22}
\end{eqnarray}
where 
\begin{eqnarray}
S_1(\varepsilon,x_0)
&=& \int_{0}^{x_0}(\overline{F}_U(0) - \overline{F}_U(x)) \rd x,
\label{defn-S_1}
\\
S_2(\varepsilon,x_0)
&=&  \int_{x_1}^{x_2}
\left( \varepsilon l_0(x_2) \overline{F}_Y(x_2) - \overline{F}_U(x)\right) \rd x,
\label{defn-S_2}
\\
S_3(\varepsilon,x_0)
&=& \int_{x_2}^{\infty}
\left(\varepsilon l_0(x) \overline{F}_Y(x) - \overline{F}_U(x)\right)\rd x.
\label{defn-S_3}
\end{eqnarray}

From (\ref{add-eqn-22}), (\ref{defn-S_1}) and $\int_{0}^{\infty}
\overline{F}_U(x) \rd x = \EE[U^+] \le \EE[|U|]< \infty$, we have
\begin{equation}
\lim_{x_0 \to \infty}\sum_{j=1}^3 S_j(\varepsilon,x_0)
\ge \lim_{x_0 \to \infty} S_1(\varepsilon,x_0) = \infty.
\label{add-eqn-23}
\end{equation}
From (\ref{defn-S_2}) and (\ref{defn-S_3}), we also have
\begin{eqnarray*}
S_2(\varepsilon,x_0) + S_3(\varepsilon,x_0)
&\le& \int_{x_1}^{\infty}
\left( \varepsilon l_0(x) \overline{F}_Y(x)  - \overline{F}_U(x)\right)\rd x
\le \int_{x_1}^{\infty}\varepsilon\overline{F}_Y(x)\rd x,
\end{eqnarray*}
where the second inequality follows from $l_0(x) \le 1$ for $x \ge 0$.
Therefore since $\EE[Y] < \infty$,
\[
\lim_{\varepsilon \downarrow 0}
(S_2(\varepsilon,x_0) + S_3(\varepsilon,x_0)) = 0,
\]
which leads to
\begin{eqnarray}
\lim_{x_0 \downarrow 0}
\lim_{\varepsilon \downarrow 0} \sum_{j=1}^3 S_j(\varepsilon,x_0)
= \lim_{x_0 \downarrow 0} \lim_{\varepsilon \downarrow 0} S_1(\varepsilon,x_0) = 0.
\label{add-eqn-24}
\end{eqnarray}
According to (\ref{add-eqn-22}), (\ref{add-eqn-23}) and
(\ref{add-eqn-24}), we can fix $\EE[Z(\varepsilon,x_0)] - \EE[U] \in
(0,y)$ for any $y > 0$.  As a result, the statement of
Lemma~\ref{lem-upper-calC} holds for $Z = Z(\varepsilon,x_0)$.

\subsection{Proof of Lemma~\ref{add-prop-N(x)}}\label{proof-add-prop-N(x)}

Let $\widetilde{X}_n$'s ($n=1,2,\dots$) are independent copies of
$\widetilde{X} := X/\EE[X]$. We then have
\begin{eqnarray*}
\lefteqn{
\{N_X(x) > u + x/\EE[X]\} 
}
~~&&
\nonumber
\\
&\subseteq& \left\{
\sum_{n=1}^{\lfloor u + x/\EE[X] \rfloor} X_n \le x
\right\}
= \left\{
\sum_{n=1}^{\lfloor u + x/\EE[X] \rfloor} (1 - \widetilde{X}_n) \ge \lfloor u + x/\EE[X] \rfloor - x/\EE[X]
\right\}
\\
&\subseteq&
\left\{
\sum_{n=1}^{\lfloor u + x/\EE[X] \rfloor} (1 - \widetilde{X}_n) \ge u - 1
\right\},
\end{eqnarray*}
which leads to
\begin{eqnarray}
\PP\left(N_X(x) -  {x \over \EE[X]} > u \right) 
&\le& 
 \PP\left(\sum_{n=1}^{\lfloor u + x/\EE[X] \rfloor} 
(1 - \widetilde{X}_n) \ge u-1 \right).
\label{add-ineqn-01}
\end{eqnarray}
Using Markov's inequality (see, e.g.,
\cite{Will91}), we have for any $s > 0$,
\begin{eqnarray}
\PP\left(\sum_{n=1}^{\lfloor u + x/\EE[X] \rfloor} 
(1 - \widetilde{X}_n) \ge u-1 \right)
&\le& e^{-s(u-1)} 
\left(\EE[e^{s(1 - \widetilde{X})}] \right)^{\lfloor u + x/\EE[X] \rfloor}
\nonumber
\\
&\le& e^{-s(u-1)} 
\left(\EE[e^{s(1 - \widetilde{X})}] \right)^{u + x/\EE[X]}
\nonumber
\\
&=& e^{s(1+x/\EE[X])} 
\left(\EE[e^{-s\widetilde{X}}] \right)^{u + x/\EE[X]},
\label{add-ineqn-02}
\end{eqnarray}
where the second inequality follows from $\EE[e^{s(1 - \widetilde{X})}]
\ge \exp\{s(1- \EE[\widetilde{X}])\} = 1$ due to
$\EE[\widetilde{X}]=1$ and Jensen's inequality (see, e.g.,
\cite{Will91}). Further for any $s > 0$,
\begin{equation}
\EE[e^{-s\widetilde{X}}]
\le 1 - s\EE[\widetilde{X}] + s^2 \EE[\widetilde{X}^2]
= 1 - s + s^2 \EE[\widetilde{X}^2],
\label{add-ineqn-03}
\end{equation}
because $e^{-x} \le 1 - x + x^2$ for all $x \ge 0$ and
$\widetilde{X} \ge 0$ w.p.1. Substituting (\ref{add-ineqn-03}) into
(\ref{add-ineqn-02}), we obtain
\begin{eqnarray}
\PP\left(\sum_{n=1}^{\lfloor u + x/\EE[X] \rfloor} 
(1 - \widetilde{X}_n) \ge u-1 \right)
&\le& e^{s(1+x/\EE[X])} 
\left(1 - s + s^2\EE[\widetilde{X}^2] \right)^{u + x/\EE[X]}\nonumber
\\
&\le& e^{s(1+x/\EE[X])} e^{(- s + s^2\EE[\widetilde{X}^2])(u + x/\EE[X])}
\nonumber
\\
&\le& e^s
\exp\left\{- su + s^2 \cdot \EE[\widetilde{X}^2](\delta+1/\EE[X]) \cdot x\right\}
\nonumber
\\
&=:& e^s
\exp\left\{ - su + s^2 \widetilde{K}(\delta)x \right\},
\label{add-ineqn-04}
\end{eqnarray}
where we use $1+x \le e^x$ ($x\in\bbR$) and $u \le \delta x$ in the
second and third inequalities.

Finally, letting $s=(u/x)\{2\widetilde{K}(\delta)\}^{-1}$ in
(\ref{add-ineqn-04}) and using $u/x \le \delta$, we obtain
\begin{eqnarray}
\PP\left(\sum_{n=1}^{\lfloor u + x/\EE[X] \rfloor} 
(1 - \widetilde{X}_n) \ge u-1 \right)
&\le&  \exp\left\{ {u \over x}{1 \over 2\widetilde{K}(\delta)} \right\}
\cdot \exp\left\{ -{1 \over 4\widetilde{K}(\delta)} {u^2 \over x} \right\}
\nonumber
\\
&\le&  \exp\left\{ {\delta \over 2\widetilde{K}(\delta)} \right\}
\cdot \exp\left\{ -{1 \over 4\widetilde{K}(\delta)} {u^2 \over x} \right\}.
\label{add-3rd-03}
\end{eqnarray}
Note here that
$\widetilde{K}(\delta)=\EE[\widetilde{X}^2](\delta+1/\EE[X])$ is
finite and positive for any fixed $\delta > 0$. As a result,
substituting (\ref{add-3rd-03}) into (\ref{add-ineqn-01}) yields
(\ref{eqn-prob-N(x)}).

\subsection{Proof of Lemma~\ref{lem-LD-maxima}}\label{proof-lem-LD-maxima}

For all $n=1,2,\dots$ and $k=1,2,\dots,n$,
\begin{eqnarray*}
\left\{ \max_{1\le k\le n}\sum_{i=1}^k U_i \ge x \right\}
&=& \bigcup_{1 \le k \le n} \left\{ \sum_{i=1}^k U_i \ge x \right\},
\quad x > 0,
\\
\left\{ \sum_{i=1}^k U_i \ge x \right\} 
&\subseteq& \bigcup_{1 \le i \le k} \left\{ U_i \ge x/k \right\} 
\subseteq \bigcup_{1 \le i \le k} \left\{ U_i \ge x/n \right\},
\quad x > 0.
\end{eqnarray*}
Thus for any fixed positive integer $n_0$ and all $n=1,2,\dots,n_0$,
we have
\[
\left\{ \max_{1\le k\le n}\sum_{i=1}^k U_i \ge x \right\}
\subseteq \bigcup_{1 \le i \le n} \left\{ U_i \ge x/n \right\}
\subseteq \bigcup_{1 \le i \le n} \left\{ U_i \ge x/n_0 \right\},
\quad x > 0,
\]
which leads to
\[
\PP\left(\max_{1\le k\le n}\sum_{i=1}^k U_i \ge x \right)
\le n \PP(U \ge x / n_0 ),
\quad n=1,2,\dots,n_0,~x>0.
\]
Therefore it suffices to prove that (\ref{add-eqn-35}) holds for all
sufficiently large $n$.

Let $\widetilde{U}_i = \min(U_i,vx)$ for $i=1,2,\dots$, where $0 < v <
1/n_0$ is a constant. Since $\EE[U]=0$, we have $\EE[\widetilde{U}_1]
\le 0$.  Thus for all $x > 0$,
\begin{eqnarray}
\PP\left(\max_{1\le k\le n}\sum_{i=1}^k U_i \ge x \right)
&\le& \PP\left(\max_{1\le i \le n} U_i > vx \right)
+ \PP\left(\max_{1\le k\le n}\sum_{i=1}^k U_i \ge x, \max_{1\le i \le n} U_i \le vx  \right)
\nonumber
\\
&\le&  n\PP(U > vx)
+ \PP\left(\max_{1\le k\le n}\sum_{i=1}^k \widetilde{U}_i \ge x \right)
\nonumber
\\
&\le&  n\PP(U > vx)
+ \PP\left(\max_{1\le k\le n}\sum_{i=1}^k W_i \ge x \right),
\label{add-eqn-34}
\end{eqnarray}
where $W_i = \widetilde{U}_i - \EE[\widetilde{U}_1]$ for
$i=1,2,\dots$.  In what follows, we estimate the second term on the
right hand side of (\ref{add-eqn-34}).

Since $\{\sum_{i=1}^k W_i\};k=1,2,\dots\}$ is martingale,
$\{\exp\{s\sum_{i=1}^k W_i\};k=1,2,\dots\}$ is submartingale for any
$s > 0$ (see, e.g., \cite[Section 14.6, Lemma (b)]{Will91}). It thus
follows from Doob's submartingale inequality (see, e.g., \cite[Section
  14.6, Theorem (a)]{Will91}) that for any $s > 0$,
\begin{eqnarray}
\PP\left(\max_{1\le k\le n}\sum_{i=1}^k W_i \ge x \right)
&=& \PP\left(\max_{1\le k\le n}\exp\left\{s\sum_{i=1}^k W_i\right\} \ge e^{sx} 
\right)
\nonumber
\\
&\le& e^{-sx} \EE\left[\exp\left\{s\sum_{i=1}^n W_i\right\} \right]
= e^{-sx} \left(\EE[e^{s \widetilde{U}_1}] \right)^n  
e^{-sn\EE[\widetilde{U}_1]}. \quad
\label{add-eqn-20}
\end{eqnarray}

We first estimate $e^{-sx} (\EE[e^{s \widetilde{U}_1}])^n$ on the
right hand side of (\ref{add-eqn-20}).  Let $1 < q < \min(r,2)$ and
\begin{equation}
s = {1 \over vx} \log\left( {v^{q-1}x^q \over n\EE[(U^+)^q]} + 1\right).
\label{defn-h}
\end{equation}
Following the estimation of the right hand side of (2.4)
in \cite{Tang06}, we can prove that there exist some positive constant
$\widetilde{C}_1:=\widetilde{C}_1(v,\gamma)$ and some positive integer
$n_1$ such that
\[
e^{-sx} \left(\EE[e^{s \widetilde{U}_1}] \right)^n
\le \widetilde{C}_1x^{-(q-1)/(2v)},
\qquad \forall x \ge \gamma n,~\forall n \ge n_1.
\]
Fix $n_0 = n_1$ and $v:=v(r,p)$ such that $0 < v < 1/n_0$ and
$(q-1)/(2v) > p$. We then have
\begin{equation}
e^{-sx} \left(\EE[e^{s \widetilde{U}_1}] \right)^n
\le \widetilde{C}_1x^{-p},
\qquad \forall x \ge \gamma n,~\forall n \ge n_0.
\label{add-eqn-12}
\end{equation}

Next we estimate $e^{-sn\EE[\widetilde{U}_1]}$ on the right hand side of
(\ref{add-eqn-20}). From (\ref{defn-h}), $\EE[\widetilde{U}_1] \le 0$,
$x \ge \gamma n$ and $n \ge 1$, we have
\begin{eqnarray}
-sn\EE[\widetilde{U}_1]
&\le& {1 \over v\gamma} \log\left( {v^{q-1}x^q \over \EE[(U^+)^q]} + 1\right)
(-\EE[\widetilde{U}_1]).
\label{ineqn--hn-EE[widetilde{U}_1]}
\end{eqnarray}
Note here that 
\begin{eqnarray*}
\EE[\widetilde{U}_1] = \EE[U \cdot \dd{1}(U \le vx)] &+& vx\PP(U > vx),
\\
\EE[U \cdot \dd{1}(U \le vx)] &+& \EE[U \cdot \dd{1}(U > vx)] = \EE[U] = 0
\end{eqnarray*}
It thus follows from $\PP(U^+ > x) =
o(x^{-r})$ (due to $\EE[(U^+)^r] < \infty$) that for all $x > 0$,
\begin{eqnarray*}
-\EE[\widetilde{U}_1] 
&\le& -\EE[U \cdot \dd{1}(U \le vx)] = \EE[U \cdot \dd{1}(U > vx)]
= \EE[U^+ \cdot \dd{1}(U^+ > vx)]
\nonumber
\\
&=& vx\PP(U^+ > vx) + \int_{vx}^{\infty} \PP(U^+ > y) \rd y
= o(x^{-r+1}).
\end{eqnarray*}
This equation and (\ref{ineqn--hn-EE[widetilde{U}_1]}) imply that for
all $x \ge \gamma n$ and $n=1,2,\dots$,
\begin{equation}
e^{-sn\EE[\widetilde{U}_1]} \le \widetilde{C}_2 < \infty,
\label{add-eqn-13}
\end{equation}
where
\[
\widetilde{C}_2
:=\widetilde{C}_2(v,\gamma)
= 
\sup_{x \ge \gamma}
\exp\left\{
{1 \over v\gamma} \log\left( {v^{q-1}x^q \over \EE[(U^+)^q]} + 1\right)
Cx^{-r+1}
\right\}.
\]

Substituting (\ref{add-eqn-12}) and (\ref{add-eqn-13}) into
(\ref{add-eqn-20}) and letting
$\widetilde{C}:=\widetilde{C}(v,\gamma)=\widetilde{C}_1(v,\gamma)\widetilde{C}_2(v,\gamma)$
yield
\[
\PP\left(\max_{1\le k\le n}\sum_{i=1}^k W_i \ge x \right)
\le \widetilde{C} x^{-p},
\qquad \forall x \ge \gamma n,~\forall n \ge n_0.
\]
This inequality and (\ref{add-eqn-34}) show that (\ref{add-eqn-35})
holds for all $n \ge n_0$.

\subsection{Proof of Lemma~\ref{prop-Lin11}}\label{proof-prop-Lin11}

Let $V_i$'s ($i=1,2,\dots$) denote independent copies of $V := U -
\EE[U] - \varepsilon$, where $\varepsilon > 0$. Clearly, $V \le U$ and
$\EE[V] = -\varepsilon < 0$. Further since $U^+ \in \calC$, we have
$V^+ \in \calC \subset \calS^{\ast}$ (see
Remark~\ref{rem-class-D}). It thus follows from the theorem in
\cite{Kors02} that for all $x \ge (\EE[U] + \varepsilon)n$ and
$n=1,2,\dots$,
\begin{eqnarray}
\PP\left( \max_{1\le k\le n}\sum_{i=1}^k U_i \ge x \right)
&\le& \PP\left( 
\max_{1\le k\le n}\sum_{i=1}^k V_i \ge x - (\EE[U] + \varepsilon)n
\right)
\nonumber
\\
&\le& {C \over \varepsilon}\int_{x- (\EE[U] + \varepsilon)n}^{x- \EE[U]n}
\PP(V > y) \rd y
\nonumber
\\
&\le& Cn \PP(V > x- (\EE[U] + \varepsilon)n)
\nonumber
\\
&\le& Cn \PP(U > x- (\EE[U] + \varepsilon)n).
\label{add-ean-36}
\end{eqnarray}
It also follows from $U^+ \in \calC \subset \calD$ that for all $x \ge
(1+\varepsilon)(\EE[U] + \varepsilon) n$, $n=1,2,\dots$ and
$\varepsilon > 0$,
\begin{eqnarray}
\PP(U > x- (\EE[U] + \varepsilon)n)
\le \PP(U > \varepsilon x/(1+\varepsilon)) \le \widehat{C}\PP(U > x),
\label{add-ean-37}
\end{eqnarray}
where $\widehat{C}:=\widehat{C}(\gamma) \in
(0,\infty)$ is given by
\[
\widehat{C} 
= \sup_{x \ge \gamma}
{\PP(U > \varepsilon x/(1+\varepsilon)) \over \PP(U > x)}~~\mbox{with}~
\gamma:=\gamma(\varepsilon)=(1+\varepsilon)(\EE[U] + \varepsilon).
\]
According to (\ref{add-ean-36}) and (\ref{add-ean-37}), there exists
some $\widetilde{C}:=\widetilde{C}(\gamma) \in (0,\infty)$ such that
\[
\PP\left( \max_{1\le k\le n}\sum_{i=1}^k U_i \ge x \right)
\le \widetilde{C}n\PP(U > x),
~~x \ge \gamma n,~n=1,2,\dots.
\]

\subsection{Proof of Lemma~\ref{lem-convo-tail}}\label{proof-lem-convo-tail}
 It follows from (\ref{eqn-05}) that for any $\varepsilon > 0$ there
 exists some $x_0:=x_0(\varepsilon) > 0$ such that for all $x \ge x_0$,
\begin{align*}
&&
\overline{F}_{i,j}(x) &\le (\widetilde{F}_{i,j} + \varepsilon)\PP(Y > x),
& 1 &\le i \le m_0,~1 \le j \le m_1,&&
\\
&&
\overline{G}_{i,j}(x) &\le (\widetilde{G}_{i,j} + \varepsilon)\PP(Y > x),
& 1 &\le i \le m_1,~1 \le j \le m_2.&&
\end{align*}
Without loss of generality, we assume that $\lim_{\varepsilon
  \downarrow 0}x_0(\varepsilon) = \infty$.

We now define $\vc{P}=(P_{i,j})$ and $\vc{Q}=(Q_{i,j})$ as $m_0 \times
m_1$ and $m_1 \times m_2$ matrix-valued functions on $\bbR$ such that
$\overline{P}_{i,j}(x) := P_{i,j}(\infty) - P_{i,j}(x) $ and
$\overline{Q}_{i,j}(x) := Q_{i,j}(\infty) - Q_{i,j}(x)$ are given by
\begin{eqnarray*}
\overline{P}_{i,j}(x) 
&=& 
\left\{
\begin{array}{ll}
\max\left(\overline{F}_{i,j}(x),(\widetilde{F}_{i,j} + \varepsilon)\PP(Y > x_0) \right), & x < x_0,
\\
(\widetilde{F}_{i,j} + \varepsilon)\PP(Y > x), & x \ge x_0,
\end{array}
\right.
\label{eqn-overline{P}_{i,j}(x)}
\\
\overline{Q}_{i,j}(x) 
&=& \left\{
\begin{array}{ll}
\max\left(\overline{G}_{i,j}(x),(\widetilde{G}_{i,j} + \varepsilon)\PP(Y > x_0) \right), & x < x_0,
\\
(\widetilde{G}_{i,j} + \varepsilon)\PP(Y > x), & x \ge x_0,
\end{array}
\right..
\label{eqn-overline{Q}_{i,j}(x)}
\end{eqnarray*}
Clearly, $\overline{F}_{i,j}(x) \le \overline{P}_{i,j}(x)$ and
$\overline{G}_{i,j}(x) \le \overline{Q}_{i,j}(x)$ for all $x \in
\bbR$. Further,
\[
\lim_{x\to\infty}{\overline{P}_{i,j}(x) \over \PP(Y>x)}
= \widetilde{F}_{i,j} + \varepsilon =: \widetilde{P}_{i,j},
\quad
\lim_{x\to\infty}{\overline{Q}_{i,j}(x) \over \PP(Y>x)}
= \widetilde{G}_{i,j} + \varepsilon =: \widetilde{Q}_{i,j}.
\]
Therefore using Proposition~A.3 in \cite{Masu11} yields
\begin{equation}
\limsup_{x\to\infty} {\overline{\vc{F} \ast \vc{G}}(x) \over \PP(Y > x)} 
\le \lim_{x\to\infty} {\overline{\vc{P} \ast \vc{Q}}(x) \over \PP(Y > x)} 
= \widetilde{\vc{P}} \vc{Q}(\infty)
+ \vc{P}(\infty) \widetilde{\vc{Q}},
\label{add-eqn-25}
\end{equation}
where $\widetilde{\vc{P}} = (\widetilde{P}_{i,j})$ and
$\widetilde{\vc{Q}} = (\widetilde{Q}_{i,j})$. Note here that
\[
\lim_{\varepsilon \downarrow 0} 
\left(
\widetilde{\vc{P}} \vc{Q}(\infty) + \vc{P}(\infty) \widetilde{\vc{Q}} \right)
= \widetilde{\vc{F}} \vc{G}(\infty)
+ \vc{F}(\infty) \widetilde{\vc{G}}.
\]
Combining this with (\ref{add-eqn-25}), we have the first statement
(\ref{eqn-05a}). The second statement (\ref{eqn-05b}) can be proved by
induction using the first statement.

Finally we prove the third statement (\ref{eqn-05c}). Since
$\sum_{n=0}^{\infty} (\vc{F}(\infty))^n = (\vc{I} - \vc{F}(\infty)
)^{-1}$,
\[
\sum_{n=0}^{\infty} (\vc{P}(\infty))^n = (\vc{I} - \vc{P}(\infty))^{-1} 
\quad \mbox{for any sufficiently small $\varepsilon > 0$}.
\]
It thus follows from $\overline{\vc{F}^{\ast n}}(x) \le
\overline{\vc{P}^{\ast n}}(x)$ ($x \in \bbR$, $n=0,1,\dots$) and
Lemma~6 in \cite{Jele98} that
\begin{equation}
\limsup_{x\to\infty} \sum_{n=0}^{\infty}
{\overline{\vc{F}^{\ast n}}(x) \over \PP(Y > x)} 
\le \lim_{x\to\infty} \sum_{n=0}^{\infty}
{\overline{\vc{P}^{\ast n}}(x) \over \PP(Y > x)} 
= (\vc{I} -  \vc{P}(\infty) )^{-1} \widetilde{\vc{P}}
(\vc{I} -  \vc{P}(\infty) )^{-1}.
\label{add-eqn-26}
\end{equation}
Letting $\varepsilon \downarrow 0$ in (\ref{add-eqn-26}) and using the
dominated convergence theorem, we obtain the third statement
(\ref{eqn-05c}).

\subsection{Proof of Lemma \ref{lemma-LD-B(t)}}\label{proof-lem-LD-B(t)}

It follows from condition (iii) of Definition~\ref{defn-SC'} that
there exists some $x_{\ast} > 0$ such that
\begin{equation}
Q(x/3) \ge Q(x)/3^{\alpha} \ge Q(x)/3,  \qquad \forall x \ge x_{\ast}.
\label{add-proof-lem-LD-B(t)-01}
\end{equation}
Let $\eta$ denote any fixed positive number such that $\eta x_{\ast}^2
\ge 1$. We then discuss three cases: (a) $0 \le x < \eta x_{\ast}^2$,
(b) $x > \eta u^2$ and (c) $\eta x_{\ast}^2 \le x \le \eta u^2$
separately. In case (a), (\ref{prop-bound-01}) holds for $C \ge
e^{(\eta x_{\ast})^2}$ because $Ce^{-\eta x} > Ce^{-(\eta x_{\ast})^2}
\ge 1$. In case (b), (\ref{prop-bound-01}) also holds for $C \ge e$
because $Ce^{-\eta u^2/x} > Ce^{-1} \ge 1$. Therefore in what follows,
we consider case (c).

For all $t \ge 0$, we have
\begin{eqnarray*}
\breve{B}(t) - \breve{b}t 
&\le& 
 \Delta \breve{B}_0^{\ast} +  \Delta \breve{B}_{N(t-\tau_0)+1}^{\ast} 
- \min(\breve{b},0)(\Delta \tau_0 + \Delta \tau_{N(t-\tau_0)+1})
+ \sum_{i=1}^{N(t-\tau_0)} (\Delta \breve{B}_i - \breve{b}\Delta \tau_i)
\nonumber
\\
&=& \Delta \Theta_0 + \Delta \Theta_{N(t-\tau_0)+1}
+ \sum_{i=1}^{N(t-\tau_0)} (\Delta \breve{B}_i - \breve{b}\Delta \tau_i),
\end{eqnarray*}
where $N(t) = \max\{n\ge0; \sum_{i=1}^n \Delta \tau_i \le t\}$ for $t
\in \bbR$.  Thus we obtain
\begin{eqnarray}
\lefteqn{
\PP\left(\sup_{0 \le t \le x}\{\breve{B}(t) - \breve{b}t\} > u\right)
}
\quad && 
\nonumber
\\
&\le& \PP\left(\Delta \Theta_0 > {u \over 3} \right)
+ \PP\left(\Delta \Theta_1 > {u \over 3} \right)
+ \PP\left(\max_{1 \le n \le N(x-\tau_0)}
\sum_{i=1}^n (\Delta \breve{B}_i - \breve{b}\Delta \tau_i) 
> {u \over 3}
\right)
\nonumber
\\
&\le& \PP\left(\Delta \Theta_0 > {u \over 3} \right)
+ \PP\left(\Delta \Theta_1 > {u \over 3} \right)
 +
\PP\left(\max_{1 \le n \le N(x)}
\sum_{i=1}^n (\Delta \breve{B}_i - \breve{b}\Delta \tau_i) 
> {u \over 3}
\right),
\label{ineqn-prop1-01}
\end{eqnarray}
where we use the inequality $\PP(X^{(1)} + X^{(2)} + X^{(3)} > u) \le
\sum_{m=1}^3\PP(X^{(m)} > u/3)$ for any triple of r.v.s $X^{(m)}$'s
($m=1,2,3$).  Note here that $\PP\left(\Delta \Theta_n > x \right) \le
C e^{-Q(x)}$ for all $x \ge 0$ due to $\EE[\exp\{Q(\Delta \Theta_n)\}]
< \infty$ ($n=0,1$). It then follows from $\eta x_{\ast}^2 \ge 1$ that
for all $x$ and $u$ such that $\eta x_{\ast}^2 \le x \le \eta u^2$,
\[
\PP\left(\Delta \Theta_n > {u \over 3} \right)
\le C e^{-Q(u/3)} \le C \eta x_{\ast}^2 e^{-Q(u/3)}
\le C x e^{-Q(u/3)},
\qquad n= 0,1,
\]
from which and (\ref{ineqn-prop1-01}) we have
\begin{equation}
\PP\left(\sup_{0 \le t \le x}\{\breve{B}(t) - \breve{b}t\} > u\right)
\le Cxe^{-Q(u/3)} 
+
\PP\left(\max_{1 \le n \le N(x)}
\sum_{i=1}^n (\Delta \breve{B}_i - \breve{b}\Delta \tau_i) > {u \over 3}
\right).
\label{add-ineqn-prop1-01}
\end{equation}

We now fix $\delta > 0$ arbitrarily and then have
\begin{eqnarray}
\lefteqn{\PP\left(\max_{1 \le n \le N(x)}
\sum_{i=1}^n (\Delta \breve{B}_i - \breve{b}\Delta \tau_i) > {u \over 3}
\right)
}
\quad &&
\nonumber
\\
&\le& 
\PP\left(N(x) - {x \over \EE[\Delta \tau_1]} > \delta x \right)
+
\PP\left( 
\max_{1 \le n \le (\delta + 1 /\EE[\Delta \tau_1])x}
\sum_{i=1}^n (\Delta \breve{B}_i - \breve{b}\Delta \tau_i) > {u \over 3}
\right). \qquad
\label{ineqn-prop1-02}
\end{eqnarray}
Applying Lemma~\ref{add-prop-N(x)} (which requires $\EE[(\Delta
  \tau_1)^2] < \infty$) to the first term on the right hand side of
(\ref{ineqn-prop1-02}), we have
\begin{equation}
\PP\left(N(x) - {x \over \EE[\Delta \tau_1]} > \delta x \right)
\le Ce^{-cx}, \qquad x \ge 0.
\label{add-eqn-01}
\end{equation}
Note that $\Delta \Theta_1 \ge 0$ and $\Delta \Theta_1 \ge \Delta
\breve{B}_1-\breve{b}\Delta \tau_1$, which lead to $\Delta \Theta_1
\ge (\Delta \breve{B}_1-\breve{b}\Delta \tau_1)^+$. Thus $\EE[\exp\{
  Q(\Delta \Theta_1) \}] < \infty$ yields
\[
\EE[\exp\{Q((\Delta \breve{B}_1-\breve{b}\Delta \tau_1)^+)\}]< \infty.
\] 
Further it follows from $\EE[(\Delta \tau_1)^2] < \infty$,
$\EE[(\Delta \breve{B}_1)^2] < \infty$ and H\"older's inequality (see,
e.g., \cite{Will91}) that
\[
\left|\EE[\Delta \breve{B}_1 \Delta \tau_1]\right|
\le \sqrt{\EE[(\Delta \breve{B}_1)^2]} \sqrt{\EE[(\Delta \tau_1)^2]} < \infty,
\]
which implies $\EE[(\Delta \breve{B}_1-\breve{b}\Delta \tau_1)^2] <
\infty$. 

We now need the following result:
\begin{prop}[Lemma~5 in \cite{Jele04}]\label{lem-LD}
Suppose that $U,U_1,U_2,\dots$ are i.i.d.\ r.v.s in $\bbR$. If
$\EE[U^2] < \infty$ and $\EE[e^{Q(U^+)}] < \infty$ for some $Q \in
\SC$, then for all $x,u \ge 0$,
\[
\PP\left(\max_{1 \le n \le x} 
\left\{ \sum_{i=1}^n U_i -n\EE[U]\right\} > u\right)
\le C\left(e^{-cu^2/x} + xe^{-(1/2)Q(u)} \right),
\]
where $C$ and $c$ are independent of $x$ and $u$.
\end{prop}

\begin{rem}
Although $\EE[U^2] < \infty$ is not
  explicitly assumed in Lemma~5 in \cite{Jele04}, this condition is
  required to prove the lemma (see p.~110 therein).
\end{rem}

Applying Proposition~\ref{lem-LD} to the second term on the right hand
side of (\ref{ineqn-prop1-02}) and using $\EE[\Delta
  \breve{B}_1-\breve{b}\Delta \tau_1] = 0$, we obtain
\begin{equation}
\PP\left( 
\max_{1 \le n \le (\delta + 1 /\EE[\Delta \tau_1])x}
\sum_{i=1}^n (\Delta \breve{B}_i - \breve{b}\Delta \tau_i) > {u \over 3}
\right)
\le C\left(e^{-cu^2/x} + xe^{-(1/2)Q(u/3)} \right).
\label{add-proof-lem-LD-B(t)-04}
\end{equation}
Substituting (\ref{add-eqn-01}) and (\ref{add-proof-lem-LD-B(t)-04})
into (\ref{ineqn-prop1-02}) yields
\[
\PP\left(\max_{1 \le n \le N(x)}
\sum_{i=1}^n (\Delta \breve{B}_i - \breve{b}\Delta \tau_i) > {u \over 3}
\right)
\le C\left(e^{-cx} + e^{-cu^2/x} + xe^{-(1/2)Q(u/3)} \right),
\]
from which and (\ref{add-ineqn-prop1-01}), we have
\begin{equation}
\PP\left(\sup_{0 \le t \le x}\{\breve{B}(t) - \breve{b}t\} > u\right)
\le C\left( e^{-cx} + e^{-cu^2/x} + xe^{-(1/2)Q(u/3)} \right).
\label{ineqn-prop1-04}
\end{equation}
Recall that in case (c), we have $u \ge x_{\ast}$ and thus $Q(u/3) \ge
(1/ 3)Q(u)$ due to (\ref{add-proof-lem-LD-B(t)-01}). Finally,
(\ref{ineqn-prop1-04}) yields (\ref{prop-bound-01}).

\subsection{Proof of Proposition~\ref{prop-00}}\label{proof-prop-00}

We first prove (\ref{add-lem00-01}). The Taylor expansion of
$(x+y)^{\gamma}$ is given by
\[
(x+y)^{\gamma}
= \sum_{n=0}^{\infty} 
{\gamma(\gamma-1)\cdots(\gamma-n+1) \over n!} y^n x^{\gamma-n},
\quad x > y \ge 0,
\]
from which we have
\begin{equation}
(x+y)^{\gamma}
\le x^{\gamma} 
+ \sum_{n=1}^{\infty} 
{|\gamma(\gamma-1)\cdots(\gamma-n+1)| \over n!} 
\left({y \over x}\right)^{n-1}
y x^{\gamma-1},\quad x > y \ge 0.
\label{eqn-lem00-01}
\end{equation}
Note here that
\begin{eqnarray*}
|\gamma(\gamma-1)\cdots(\gamma-n+1)|
&=& \gamma(\gamma-1) \cdots (\gamma - \lfloor \gamma \rfloor)
\cdot \prod_{i=1}^{n-1- \lfloor \gamma \rfloor}
(i+\lfloor \gamma \rfloor - \gamma)
\nonumber
\\
&\le& \gamma(\gamma-1) \cdots (\gamma - \lfloor \gamma \rfloor)
\cdot n!.
\end{eqnarray*}
Thus since $y/x < 1$, 
\begin{eqnarray}
\sum_{n=1}^{\infty} 
{|\gamma(\gamma-1)\cdots(\gamma-n+1)| \over n!} 
\left({y \over x}\right)^{n-1}
&\le& \gamma(\gamma-1) \cdots 
(\gamma - \lfloor \gamma \rfloor)
\cdot \sum_{n=1}^{\infty}
\left({y \over x}\right)^{n-1}
\nonumber
\\
&=& \gamma(\gamma-1) \cdots (\gamma - \lfloor \gamma \rfloor)
\cdot \left(1 - {y \over x} \right)^{-1}. \quad
\label{eqn-lem00-03}
\end{eqnarray}
Substituting (\ref{eqn-lem00-03}) into (\ref{eqn-lem00-01}) yields
(\ref{add-lem00-01}).

In the same way as the proof of (\ref{add-lem00-01}), we have
\[
(x-y)^{\gamma}
\ge x^{\gamma} 
- \sum_{n=1}^{\infty} 
{|\gamma(\gamma-1)\cdots(\gamma-n+1)| \over n!} 
\left({y \over x}\right)^{n-1}
y x^{\gamma-1},
\]
from which and (\ref{eqn-lem00-03}) we obtain (\ref{add-lem00-02}).

\section{Proofs of Main Results}\label{sec-proof-main-results}

\subsection{Proof of Theorem \ref{thm-dependent-01}}\label{proof-thm-dependent-01} 

Since $\{B(t)\}$ is nondecreasing with $t$, we have for $x \ge 1$,
\begin{eqnarray*}
\PP(B(T)>x)
&\le& \PP(T > x - x^{2/3}) 
+ \PP(B(T)>x, T \le x - x^{2/3})
\nonumber
\\
&\le& \PP(T > x - x^{2/3}) 
+ \PP(B(x - x^{2/3}) > x),
\\
\PP(B(T) > x)
&\ge& \PP(B(T) > x, T > x + x^{2/3}) 
\nonumber
\\
&=& \PP(T > x + x^{2/3})
- \PP(B(T) \le x, T > x + x^{2/3})
\nonumber
\\
&\ge& \PP(T > x + x^{2/3})
- \PP(B(x + x^{2/3}) \le x).
\end{eqnarray*}
Note here that condition (i) is equivalent to $T \in \calL^3$ (see
Lemma~\ref{lem-01}~(ii)). It thus follows from Lemma~\ref{lem-02} that
\[
\PP(T > x + x^{2/3})
\simhm{x}
\PP(T > x - x^{2/3})
\simhm{x}
\PP(T > x).
\] 
Therefore it suffices to show that 
\begin{eqnarray*}
\PP(B(x - x^{2/3}) > x) = o(\PP(T>x - x^{2/3})),
\\
\PP(B(x + x^{2/3}) \le x) = o(\PP(T>x + x^{2/3})).
\end{eqnarray*}

For $x \ge 1$, we have
\begin{eqnarray}
\PP(B(x - x^{2/3})>x)
&=& 
\PP\left( B(x - x^{2/3}) - (x - x^{2/3}) > x^{2/3} 
\right)
\nonumber
\\
&\le& 
\PP\left( \sup_{0 \le t \le x - x^{2/3}}(B(t)-t) >  x^{2/3} 
\right)
\nonumber
\\
&\le& 
\PP\left( \sup_{0 \le t \le x}(B(t) - t) > x^{2/3} \right),
\label{thm02-eqn-02}
\end{eqnarray}
and 
\begin{eqnarray}
\PP(B(x + x^{2/3}) \le x)
&\le&
\PP(B(x + x^{2/3}) < x + (1/2)x^{2/3})
\nonumber
\\
&=& \PP\left( B(x + x^{2/3}) - (x + x^{2/3}) < -(1/2)x^{2/3}
\right)
\nonumber
\\
&\le& \PP\left(
\inf_{0 \le t \le x + x^{2/3}}(B(t) - t) < -(1/2)x^{2/3}
\right)
\nonumber
\\
&\le&
\PP\left(
\inf_{0 \le t \le 2x}(B(t) - t) < -(1/2)x^{2/3}
\right).
\qquad
\label{thm02-eqn-06}
\end{eqnarray}
Applying Lemma~\ref{prop-LD-B(t)} (i) to the right hand side of
(\ref{thm02-eqn-02}), we obtain
\begin{equation}
\PP(B(x - x^{2/3})>x)
\le
C\left(e^{-cx^{1/3} }
+ e^{-cx}
+ xe^{-cQ(x^{2/3})}\right),
\qquad x \ge 1.
\label{add-eqn-40}
\end{equation}
Since $T^{\theta} \in \calL$, we have $\PP(T > x)=e^{-o(x^{\theta})}$
(see Lemma~\ref{lem-01}~(i)) and thus for any $0 < \theta \le 1/3$,
\begin{eqnarray*}
\limsup_{x\to\infty}
{e^{-cx^{1/3}} \over \PP(T > x - x^{2/3})} 
&=& \limsup_{x\to\infty}e^{-cx^{1/3} + o(x^{\theta})}=0,
\\
\limsup_{x\to\infty}
{e^{-cx} \over \PP(T > x - x^{2/3})} 
&=& \limsup_{x\to\infty}e^{-cx + o(x^{\theta})}=0.
\end{eqnarray*}
Further it follows from (\ref{add-eqn-28}) and 
$\PP(T > x)=e^{-o(x^{\theta})}$ that
\[
\limsup_{x\to\infty}{xe^{-cQ(x^{2/3})} \over \PP(T>x-x^{2/3})}
\le \limsup_{x\to\infty}
e^{-cx^{\theta}+\log x + o(x^{\theta})} = 0.
\]
As a result, we have $\PP(B(x - x^{2/3}) > x)=o(\PP(T > x-x^{2/3}))$.

Next we estimate $\PP(B(x + x^{2/3}) \le x)$. Applying
Lemma~\ref{prop-LD-B(t)}~(ii) to the right hand side of
(\ref{thm02-eqn-06}) yields
\[
\PP(B(x + x^{2/3}) \le x)
\le
C\left( e^{-cx^{1/3}}
+ e^{-cx}
+ x e^{-cQ((1/2)x^{2/3})}\right).
\]
Therefore similarly to the estimation of (\ref{add-eqn-40}), we can
readily show $\PP(B(x + x^{2/3}) \le x)=o(\PP(T > x+x^{2/3}))$.  \qed

\subsection{Proof of Theorem \ref{thm-dependent-02}}\label{proof-thm-dependent-02}

We fix $\varepsilon \in (0,1)$ arbitrarily. Since $\{B(t)\}$ is
nondecreasing with $t$, we have for $x > 0$,
\begin{eqnarray*}
\PP(B(T) > x) 
&\le& \PP(T > (1-\varepsilon)x) + \PP(B(T) > x, T \le (1-\varepsilon)x)
\nonumber
\\
&\le& \PP(T > (1-\varepsilon)x) + \PP(B((1-\varepsilon)x) > x),
\\
\PP(B(T) > x) 
&\ge& \PP(B(T) > x, T > (1+\varepsilon)x)
\nonumber
\\
&=& \PP(T > (1+\varepsilon)x) - \PP(B(T) \le x, T > (1+\varepsilon)x)
\nonumber
\\
&\ge& \PP(T > (1+\varepsilon)x) - \PP(B((1+\varepsilon)x) \le x).
\end{eqnarray*}
Since $T \in \calC$ (see Definition~\ref{defn-class-C}),
\begin{eqnarray}
\lim_{\varepsilon\downarrow0}
\liminf_{x\to\infty}{\PP(T > (1+\varepsilon)x) \over \PP(T > x)} &=& 1,
\label{property-T-in-C-01}
\\
\lim_{\varepsilon\downarrow0}
\limsup_{x\to\infty}{\PP(T > (1-\varepsilon)x) \over \PP(T > x)} &=& 1.
\label{property-T-in-C-02}
\end{eqnarray}
Therefore it suffices to show that
\begin{eqnarray}
\PP(B((1-\varepsilon)x) > x) &=& o(\PP(T > x)),
\label{estimate-B(T)-01}
\\
\PP(B((1+\varepsilon)x) \le x) &=& o(\PP(T > x)).
\label{estimate-B(T)-02}
\end{eqnarray}

For $x > 0$, we have
\begin{eqnarray}
\PP(B((1-\varepsilon)x) > x)
&=& \PP(B((1-\varepsilon)x) - (1-\varepsilon)x > \varepsilon x)
\nonumber
\\
&\le& \PP\left(\sup_{0 \le t \le (1-\varepsilon)x} (B(t) - t) > \varepsilon x \right)
\nonumber
\\
&\le& \PP\left(\sup_{0 \le t \le (1+\varepsilon)x} (B(t) - t) > \varepsilon x \right),
\label{add-eqn-02a}
\\
\PP(B((1+\varepsilon)x) \le x)
&=& \PP(B((1+\varepsilon)x) - (1+\varepsilon)x \le -\varepsilon x)
\nonumber
\\
&\le& \PP(B((1+\varepsilon)x) - (1+\varepsilon)x < -\varepsilon x/2)
\nonumber
\\
&\le& \PP\left(\sup_{0 \le t \le (1+\varepsilon)x} (t - B(t)) > \varepsilon x/2 \right).
\label{add-eqn-02b}
\end{eqnarray}
Note here that 
\begin{eqnarray*}
\sum_{i=1}^{N(t-\tau_0)}\Delta \tau_i
&\le& t \le \sum_{i=0}^{N(t-\tau_0)+1}\Delta \tau_i,
\\
\sum_{i=1}^{N(t-\tau_0)}\Delta B_i
&\le& B(t) - B(0) \le \sum_{i=0}^{N(t-\tau_0)+1}\Delta B_i - B(0),
\end{eqnarray*}
where $N(t) = \max\{n\ge0; \sum_{i=1}^n \Delta \tau_i \le t\}$ for $t
\in \bbR$.  We thus have
\begin{eqnarray}
B(t) - t
&\le& \Delta B_0 + \Delta B_{N(t-\tau_0)+1} 
+ \sum_{i=1}^{N(t-\tau_0)}(\Delta B_i - \Delta \tau_i),
\label{ineqn-t-B(t)a}
\\
t - B(t) 
&\le& \Delta \tau_0 + \Delta \tau_{N(t-\tau_0)+1} 
+ \sum_{i=1}^{N(t-\tau_0)}(\Delta \tau_i - \Delta B_i)- B(0).
\label{ineqn-t-B(t)b}
\end{eqnarray}
Therefore similarly to (\ref{ineqn-prop1-01}), it follows from
(\ref{add-eqn-02a}) and (\ref{ineqn-t-B(t)a}) that
\begin{eqnarray}
\PP(B((1-\varepsilon)x) > x)
&\le& \PP(\Delta B_0 > \varepsilon x/3) + \PP(\Delta B_1 > \varepsilon x/3) 
\nonumber
\\
&& {} +  \PP\left(\max_{1 \le k \le N((1+\varepsilon)x)} 
\sum_{i=1}^{k} (\Delta B_i - \Delta \tau_i) > {\varepsilon \over 3}x \right),
\quad x > 0, \quad
\label{add-eqn-03a}
\end{eqnarray}
and it follows from (\ref{add-eqn-02b}) and (\ref{ineqn-t-B(t)b}) that
\begin{eqnarray}
\PP(B((1+\varepsilon)x) \le x)
&\le& 
\PP(-B(0) > \varepsilon x/8)
+\PP(\Delta \tau_0 > \varepsilon x/8) + \PP(\Delta \tau_1 > \varepsilon x/8) 
\nonumber
\\
&& {} +  \PP\left(\max_{1 \le k \le N((1+\varepsilon)x)} 
\sum_{i=1}^{k} (\Delta \tau_i - \Delta B_i) > {\varepsilon \over 8}x \right),
\quad x > 0. \quad
\label{add-eqn-03b}
\end{eqnarray}
Since $T \in \calC \subset \calD$, condition
(iii) yields
\[
\lim_{x\to\infty}{\PP(-B(0) > \varepsilon x/8) \over \PP(T > x)}
\le \limsup_{x\to\infty}
{\PP(-B(0) > \varepsilon x/8) \over \PP(T > \varepsilon x/8)}
\limsup_{x\to\infty}{\PP(T > \varepsilon x/8) \over \PP(T > x)} = 0,
\]
which shows that $\PP(-B(0) > \varepsilon x/8) = o(\PP(T > x))$. In
addition, $\PP(\Delta \tau_n > \varepsilon x/8) = o(\PP(T > x))$ and
$\PP(\Delta B_n > \varepsilon x/3) = o(\PP(T > x))$ ($n=0,1$).

As for the last terms in (\ref{add-eqn-03a}) and (\ref{add-eqn-03b}),
we have for any $\delta > 0$,
\begin{eqnarray}
\lefteqn{
\PP\left(\max_{1 \le k \le N((1+\varepsilon)x)} 
\sum_{i=1}^{k} (\Delta B_i - \Delta \tau_i)
> {\varepsilon \over 3}x \right)
}
\quad &&
\nonumber
\\
&\le& 
\PP\left(
N((1+\varepsilon)x) - {(1+\varepsilon)x \over \EE[\Delta \tau_1]} > \delta x 
\right)
\nonumber
\\
&& {}
+
\PP\left(N((1+\varepsilon)x) - {(1+\varepsilon)x \over \EE[\Delta \tau_1]} \le \delta x,
\max_{1 \le k \le N((1+\varepsilon)x)}
\sum_{i=1}^k (\Delta B_i - \Delta \tau_i) > {\varepsilon \over 3}x
\right)
\nonumber
\\
&\le& 
\PP\left(
N((1+\varepsilon)x) - {(1+\varepsilon)x \over \EE[\Delta \tau_1]} > \delta x 
\right)
\nonumber
\\
&& {}
+
\PP\left(
\max_{1 \le k \le \{\delta + (1+\varepsilon) /\EE[\Delta \tau_1]\}x}
\sum_{i=1}^k (\Delta B_i - \Delta \tau_i) > {\varepsilon \over 3}x
\right),\qquad x > 0,
\label{add-eqn-04a}
\end{eqnarray}
and
\begin{eqnarray}
\lefteqn{
\PP\left(\max_{1 \le k \le N((1+\varepsilon)x)} 
\sum_{i=1}^{k} (\Delta \tau_i - \Delta B_i)
> {\varepsilon \over 8}x \right)
}
\quad &&
\nonumber
\\
&\le& 
\PP\left(
N((1+\varepsilon)x) - {(1+\varepsilon)x \over \EE[\Delta \tau_1]} > \delta x 
\right)
\nonumber
\\
&& {}
+
\PP\left(
\max_{1 \le k \le \{\delta + (1+\varepsilon) /\EE[\Delta \tau_1]\}x}
\sum_{i=1}^k (\Delta \tau_i - \Delta B_i) > {\varepsilon \over 8}x
\right),\qquad x > 0.
\label{add-eqn-04b}
\end{eqnarray}
According to Lemma~\ref{add-prop-N(x)}, the first terms in
(\ref{add-eqn-04a}) and (\ref{add-eqn-04b}) are bounded from above by
$Ce^{-cx} = o(\PP(T > x))$. Let 
\begin{eqnarray}
\gamma 
={\varepsilon \over 8} 
\cdot {1 \over \delta + (1+\varepsilon)/\EE[\Delta \tau_1]}.
\label{eqn-gamma-01}
\end{eqnarray}
We then have $\varepsilon x/(8\gamma)=\{\delta + (1+\varepsilon)
/\EE[\Delta \tau_1]\}x$ and thus for $x>0$,
\begin{eqnarray*}
\PP\left(
\max_{1 \le k \le \{\delta + (1+\varepsilon) /\EE[\Delta \tau_1]\}x}
\sum_{i=1}^k (\Delta B_i - \Delta \tau_i) > {\varepsilon \over 3}x
\right)
\le
\PP\left(
\max_{1 \le k \le \varepsilon x/(8\gamma)}
\sum_{i=1}^k (\Delta B_i - \Delta \tau_i) > {\varepsilon \over 8}x
\right).
\end{eqnarray*}
As a result, to prove
(\ref{estimate-B(T)-01}) and (\ref{estimate-B(T)-02}), it suffices to
show that
\begin{eqnarray}
\PP\left(
\max_{1 \le k \le \varepsilon x/(8\gamma)}
\sum_{i=1}^k (\Delta B_i - \Delta \tau_i) > {\varepsilon \over 8}x
\right)
&=& o(\PP(T > x)),
\label{2nd-term-A}
\\
\PP\left(
\max_{1 \le k \le \varepsilon x/(8\gamma)}
\sum_{i=1}^k (\Delta \tau_i - \Delta B_i) > {\varepsilon \over 8}x
\right)
&=& o(\PP(T > x)).
\label{2nd-term-B}
\end{eqnarray}
In what follows, we prove (\ref{2nd-term-A}) under condition (v.a) and
condition (v.b), separately. We omit the proof of (\ref{2nd-term-B}),
which is almost the same as that of (\ref{2nd-term-A}).

\subsubsection{Condition (v.a)}\label{subsubsec-v.a}

Suppose condition (v.a) holds. 
It then follows from Lemma~\ref{lem-LD-maxima} that for any fixed $p >
0$,
\begin{eqnarray*}
\PP\left(
\max_{1 \le k \le \varepsilon x/(8\gamma)}
\sum_{i=1}^k (\Delta B_i - \Delta \tau_i) > {\varepsilon \over 8}x
\right)
&\le& 
Cx\PP(\Delta B_1 - \Delta \tau_1 > vx) + Cx^{-p}
\\
&&{}\mbox{for all sufficiently large $x>0$},
\end{eqnarray*}
where $v:=v(r,p)$ is some finite positive constant.  Further
from $T \in \calC \subset \calD$ and condition (iv), we have
\begin{eqnarray*}
\lim_{x\to\infty}{x\PP(\Delta B_1 - \Delta \tau_1 > vx) \over \PP(T > x)}
&\le& \limsup_{x\to\infty}
{x\PP(\Delta B_1 - \Delta \tau_1 > vx) \over \PP(T > vx)}
\limsup_{x\to\infty}{\PP(T > vx) \over \PP(T > x)} = 0.
\end{eqnarray*}

We now fix $p > r_+$, where $r_+$ denotes the upper Matuszewska index
of the d.f.\ of $\Delta B_1 - \Delta \tau_1$ (see
subsection~\ref{subsec-class-D}). It then follows from
Proposition~\ref{prop-Bing89} and condition (iv) that $x^{-p} =
o(\PP(\Delta B_1 - \Delta \tau_1 > x))$ and thus $x^{-p} = o(\PP(T >
x))$. Therefore (\ref{2nd-term-A}) holds.

\subsubsection{Condition (v.b)}\label{subsubsec-v.b}

Suppose that condition (v.b) holds and define $Y$ as a nonnegative r.v.\ 
such that
\[
\PP(Y > x) = \min(1,c\PP(T > x)/x),\qquad x > 0.
\]
It then follows from $T \in \calC$ and conditions (iv) and (v.b) that
\[
Y \in \calC,\quad \EE[Y] < \infty,\quad
\PP(\Delta B_1 - \Delta \tau_1 > x)= o(\PP(T>x)/x) = o(\PP(Y>x)).
\]
Therefore Lemma~\ref{lem-upper-calC} implies that there exists a
r.v.\ $Z$ in $\bbR$ such that $0 < \EE[Z] < \gamma$,
\begin{eqnarray}
\PP(Z > x) &\ge& \PP(\Delta B_1 - \Delta \tau_1 > x)
\quad \mbox{for all $x \in \bbR$ and},
\label{defn-Z-01}
\\
\PP(Z > x) &=& \tilde{l}(x)\PP(Y>x)
\quad \mbox{for all sufficiently large $x > 0$},
\label{defn-Z-02}
\end{eqnarray}
where $\gamma$ is given in (\ref{eqn-gamma-01}) and $\tilde{l}$ is
some slowly varying function such that
$\lim_{x\to\infty}\tilde{l}(x)=0$.

The inequality (\ref{defn-Z-01}) enables us to assume that $Z$ and
$\Delta B_1 - \Delta \tau_1$ are on the same probability space and $Z
\ge \Delta B_1 - \Delta \tau_1$, without loss of generality (see,
e.g., Theorem 1.2.4 in \cite{Mull02}). We thus have
\begin{eqnarray}
\PP\left(
\max_{1 \le k \le \varepsilon x/(8\gamma)}
\sum_{i=1}^k (\Delta B_i - \Delta \tau_i) > {\varepsilon \over 8}x
\right)
&\le& \PP\left(
\max_{1 \le k \le \varepsilon x/(8\gamma)}
\sum_{i=1}^k Z_i > {\varepsilon \over 8}x
\right),
\label{add-eqn-27}
\end{eqnarray}
where $Z_i$'s ($i=1,2,\dots$) are independent copies of $Z$. Note here
that $Z \in \calC$ due to $Y \in \calC$ and
(\ref{defn-Z-02}). Therefore applying Lemma~\ref{prop-Lin11} to the
right hand side of (\ref{add-eqn-27}) yields
\begin{eqnarray*}
\PP\left(
\max_{1 \le k \le \varepsilon x/(8\gamma)}
\sum_{i=1}^k \Delta B_i - \Delta \tau_i > {\varepsilon \over 8}x
\right)
\le Cx\PP(Z > \varepsilon x/8)~~\mbox{for all sufficiently large $x>0$}.
\end{eqnarray*}
In addition, it follows from condition (iv) and the definitions of $Y$
and $Z$ that $\PP(Z > x) = o(\PP(T>x)/x)$. Using this and $T \in \calC
\subset \calD$, we have
\[
\lim_{x\to\infty}{x\PP(Z > \varepsilon x/8) \over \PP(T > x)} 
\le
\limsup_{x\to\infty}{x\PP(Z > \varepsilon x/8) \over \PP(T > \varepsilon x/8)} 
\limsup_{x\to\infty}{\PP(T > \varepsilon x/8) \over \PP(T > x)} = 0.
\]
As a result, we obtain (\ref{2nd-term-A}).

\subsection{Proof of Theorem \ref{thm-independ-01}}\label{proof-thm-main}

The asymptotic lower bound $\PP(B(T) > x) \gesimhm{x} \PP(T > x)$ can
be proved in the same way as the proof of Theorem~3 in
\cite{Jele04}. Thus we here prove only the asymptotic upper bound
$\PP(M(T) > x) \lesimhm{x} \PP(T > x)$.

We fix $\delta$ $(0 < \delta < 1)$ arbitrarily and also fix $x$ such
that $0 < \delta x \le x - \xi \sqrt{x}$ and $\xi \ge 1$, which leads
to $\sqrt{x} \ge \xi /(1-\delta) > 1$. We then have
\begin{eqnarray}
\PP(M(T) > x)
&=& \PP(M(T) > x, T > x - \xi \sqrt{x})
\nonumber
\\
&& {} + \PP(M(T) > x, \delta x < T \le x - \xi \sqrt{x})
 + \PP(M(T) > x, T \le \delta x) 
\nonumber
\\
&\le& \PP(T > x - \xi \sqrt{x})
\nonumber
\\
&&{} + 
\PP(M(T) > x, \delta x < T \le x - \xi \sqrt{x})
+ \PP(M(\delta x) > x).
\qquad
\label{thm3-ineqn-00}
\end{eqnarray}
Since $\PP(T > x - \xi \sqrt{x}) \simhm{x} \PP(T > x)$ (due to
condition~(i); see Lemmas~\ref{lem-01} and \ref{lem-02}), it suffices
to show that the second and third terms in (\ref{thm3-ineqn-00}) are
$o(\PP(T>x))$.

Note first that $\sup_{0 \le t \le \delta x}B(t) - \delta x \le
\sup_{0 \le t \le \delta x}(B(t) - t)$ and thus
\begin{eqnarray}
\PP(M(\delta x) > x)
&\le& \PP\left(\sup_{0 \le t \le \delta x}(B(t) - t) > (1 - \delta) x 
\right).
\label{thm3-ineqn-01}
\end{eqnarray}
Applying Lemma~\ref{prop-LD-B(t)}~(i) to (\ref{thm3-ineqn-01})
yields
\begin{eqnarray*}
\PP(M(\delta x) > x) 
&\le& C\left(e^{-cx} 
+ xe^{-cQ((1-\delta)x)} \right)
= o(\PP(T > x)) + Cxe^{-cQ((1-\delta)x)}.
\end{eqnarray*}
Further since $x^{\theta}=O(Q(x))$ and $\PP(T > x) =
e^{-o(x^{\theta})}$ (due to $T \in \calL^{1/\theta}$; see Lemma~\ref{lem-01}~(i)),
\begin{eqnarray}
\limsup_{x\to\infty}{xe^{-cQ((1-\delta)x)} \over \PP(T > x)}
\le \limsup_{x\to\infty}
\exp\left\{-cx^{\theta}/C + \log x + o(x^{\theta})\right\}
=0.
\label{add-eqn-30}
\end{eqnarray}
Consequently, we have $\PP(M(\delta x) > x) = o(\PP(T > x))$.

Next we consider the second term on the right hand side of
(\ref{thm3-ineqn-00}). Note that
\begin{eqnarray}
\lefteqn{
\PP(M(T) > x, \delta x < T \le x - \xi \sqrt{x})
}
\quad &&
\nonumber
\\
&=& \int_{\delta x}^{x - \xi \sqrt{x}}\PP(M(u) > x) \rd\PP(T \le u)
\nonumber
\\
&\le& \int_{\delta x}^{x - \xi \sqrt{x}}
\PP\left(\sup_{0 \le t \le u}(B(t) - t) > x - u\right) 
\rd\PP(T \le u).
\qquad
\label{thm3-ineqn-06}
\end{eqnarray}
Applying Lemma~\ref{prop-LD-B(t)}~(i) to the right hand side of
(\ref{thm3-ineqn-06}) and using $\delta x \le u \le x$, we obtain
\begin{eqnarray*}
\lefteqn{
\PP(M(T) > x, \delta x < T \le x - \xi \sqrt{x})
}
\quad &&
\nonumber
\\
&\le& \int_{\delta x}^{x - \xi \sqrt{x}}
C\left(e^{-c(x-u)^2/u} + e^{-cu} 
+ ue^{-cQ(x-u)}\right) 
\rd\PP(T \le u)
\\
&\le& C e^{-c\delta x} + C\int_{\delta x}^{x - \xi \sqrt{x}}
\left(e^{-c(x-u)^2/x}  
+ xe^{-cQ(x-u)}\right) 
\rd\PP(T \le u)
\\
&=& o(\PP(T > x)) + Cf_1(x) + Cf_2(x),
\end{eqnarray*}
where
\begin{eqnarray}
f_1(x) 
&=& \int_{\delta x}^{x - \xi \sqrt{x}} e^{-c(x-u)^2/x}\rd\PP(T \le u),
\label{def-f_2(x)}
\\
f_2(x) 
&=& \int_{\delta x}^{x - \xi \sqrt{x}} xe^{-cQ(x-u)}\rd\PP(T \le u).
\label{def-f_1(x)}
\end{eqnarray}
In what follows, we prove $f_1(x) = o(\PP(T > x))$ and $f_2(x) =
o(\PP(T > x))$. 

Note that $e^{-c(x-u)^2/x}$ is differentiable with respect
to $u$. Thus integrating the right hand side of (\ref{def-f_2(x)}) by
parts (see, e.g., Theorems~6.1.7 and 6.2.2 in \cite{Cart00}) and
letting $y=(x-u)/\sqrt{x}$ yield
\begin{eqnarray}
f_1(x)
&\le& e^{-c(1-\delta)^2x} 
+ \int_{\delta x}^{x - \xi \sqrt{x}} \PP(T > u) \rd_u(e^{-c(x-u)^2/x}) 
\nonumber
\\
&=& e^{-c(1-\delta)^2x} 
+ \int_{\delta x}^{x - \xi \sqrt{x}} \PP(T > u) {2c(x-u) \over x}e^{-c(x-u)^2/x} \rd u  
\nonumber
\\
&=& o(\PP(T > x))
+  \int_{\xi}^{(1-\delta)\sqrt{x}} 
\PP(T > x - y \sqrt{x}) 
2cye^{-cy^2}\rd y
\nonumber
\\
&\le& o(\PP(T > x))
 +
\int_{\xi}^{(1-\delta)\sqrt{x}} 
\PP( \sqrt{T} > \sqrt{x} - y)
2cye^{-cy^2}\rd y,
\label{add-proposition-4-proof-pre}
\end{eqnarray}
where the last inequality holds because $(x - y\sqrt{x})^{1/2} \ge
\sqrt{x} - y$ for $0 \le y \le \sqrt{x}$. It thus follows from
$\sqrt{T} \in \calL$ and Lemma~\ref{lem-03} that for any $\varepsilon
> 0$,
\begin{eqnarray}
&& \lim_{\xi\to\infty} \limsup_{x\to\infty}
\int_{\xi}^{(1-\delta)\sqrt{x}} 
{\PP( \sqrt{T} > \sqrt{x} - y)  \over \PP(T > x)}
2cy e^{-cy^2} \rd y
\nonumber
\\
&& \qquad =  \lim_{\xi\to\infty} \limsup_{x\to\infty}
\int_{\xi}^{(1-\delta)\sqrt{x}} 
{\dm \PP( \sqrt{T} > \sqrt{x} - y)
\over \PP( \sqrt{T} > \sqrt{x}) }
2cy e^{-cy^2}\rd y
\nonumber
\\
&& \qquad \le e^{\varepsilon} \lim_{\xi\to\infty} \limsup_{x\to\infty}
\int_{\xi}^{(1-\delta)\sqrt{x}} 
2cy \exp\{-cy^2 + \varepsilon y\}\rd y
\nonumber
\\
&& \qquad \le e^{\varepsilon}
 \lim_{\xi\to\infty} \int_{\xi}^{\infty} 
2cy \exp\{-cy^2 + \varepsilon y\}\rd y = 0.
\label{add-eqn-17a}
\end{eqnarray}
Combining (\ref{add-proposition-4-proof-pre}) with
(\ref{add-eqn-17a}) yields $f_1(x) = o(\PP(T > x))$.

We proceed to the proof of $f_2(x) = o(\PP(T > x))$. Since $Q$ is
eventually concave (see Definition~\ref{defn-SC'}), $Q$ is continuous
for all sufficiently large $x > 0$. Therefore without loss of
generality, we fix $x$ to be sufficiently large such that $Q(x-u)$ is
continuous for all $\delta x \le u \le x - \xi \sqrt{x}$.

For $\delta x \le u \le x - \xi \sqrt{x}$, we have
\[
e^{-c Q(x-u)}
= e^{-(c/2)Q(x-u)}e^{-(c/2)Q(x-u)} 
\le e^{-(c/2)Q(\xi \sqrt{x})} e^{-(c/2)Q(x-u)}.
\]
Substituting this into the right hand side of (\ref{def-f_1(x)}) and
integrating it by parts yield
\begin{eqnarray*}
f_2(x)
&\le&
x e^{-cQ(\xi \sqrt{x})}
\int_{\delta x}^{x - \xi \sqrt{x}} \{- e^{-cQ(x-u)} \}
\rd\PP(T > u)
\nonumber
\\
&\le&
x e^{-cQ(\xi \sqrt{x})}
\left[ e^{-cQ((1-\delta)x)} 
+ \int_{\delta x}^{x - \xi \sqrt{x}} \PP(T > u) \rd_u (e^{-cQ(x-u)})\right]
\nonumber
\\
&=&
x e^{-cQ(\xi \sqrt{x})}
\left[ o(\PP(T > x))
+ \int_{\delta x}^{x - \xi \sqrt{x}} \PP(T > u) \rd_u (e^{-cQ(x-u)})\right],
\end{eqnarray*}
where the last equality follows from $e^{-cQ((1-\delta)x)} = o(\PP(T >
x))$ due to (\ref{add-eqn-30}). Further using $\log x = o(Q(x))$ and
$x^{\theta}=O(Q(x))$, we have
\begin{eqnarray*}
\lim_{x\to\infty}x e^{-cQ(\xi \sqrt{x})}
&=& \lim_{x\to\infty}e^{-cQ(\xi \sqrt{x}) + 2\log \sqrt{x}}
= \lim_{x\to\infty}e^{-cQ(\xi \sqrt{x}) + o(Q(\xi \sqrt{x}))}
= 0.
\end{eqnarray*}
Finally, it follows from $T^{\theta} \in \calL$ and Lemma~\ref{lem-03}
that for sufficiently small $\varepsilon > 0$,
\begin{eqnarray*}
\lefteqn{
\limsup_{x\to\infty}
\int_{\delta x}^{x - \xi \sqrt{x}} {\PP(T > u) \over \PP(T>x)}
\rd_u (e^{-cQ(x-u)})
}
\quad &&
\nonumber
\\
&\le&
e^{\varepsilon}
\limsup_{x\to\infty}
\int_{\delta x}^{x - \xi \sqrt{x}} e^{\varepsilon (x-u)^{\theta}}
\rd_u (e^{-cQ(x-u)})
\nonumber
\\
&\le&
e^{\varepsilon}\limsup_{x\to\infty}
\left[ e^{\varepsilon(\xi \sqrt{x})^{\theta} - cQ(\xi \sqrt{x})}
+ \int_{\delta x}^{x - \xi \sqrt{x}} 
\varepsilon \theta (x-u)^{\theta-1} e^{\varepsilon(x-u)^{\theta} - cQ(x-u)} \rd u
\right]
\nonumber
\\
&\le&
e^{\varepsilon}\limsup_{x\to\infty}
\left[ e^{\varepsilon(\xi \sqrt{x})^{\theta} - cQ(\xi \sqrt{x})}
+ \varepsilon \theta \int_{\xi \sqrt{x}}^{(1-\delta) x}
 e^{\varepsilon y^{\theta} - cQ(y)} \rd y
\right]
= 0,
\end{eqnarray*}
where the last equality is due to $x^{\theta}=O(Q(x))$.
As a result, we have $f_2(x) = o(\PP(T>x))$.

\subsection{Proof of Theorem~\ref{thm-independ-02}}\label{proof-thm-independ-02}

For any $\varepsilon > 0$, we have
\begin{eqnarray}
\PP(B(T) > x) 
&\ge& \int_{(1+\varepsilon)x}^{\infty}\PP(B(u) > x) \rd \PP(T \le u)
\nonumber
\\
&\ge& \inf_{u>(1+\varepsilon)x}
\PP(B(u) > x) \PP(T > (1+\varepsilon)x)
\nonumber
\\
&=& \inf_{u>(1+\varepsilon)x}
\PP\left( {B(u) - u \over u} > {x - u \over u} \right)
\PP(T > (1+\varepsilon)x )
\nonumber
\\
&\ge& \inf_{u>(1+\varepsilon)x}
\PP\left( {B(u) - u \over u} > {-\varepsilon \over 1+\varepsilon} \right)
\PP(T > (1+\varepsilon)x ).
\label{add-eqn-06}
\end{eqnarray}
It follows from the SLLN for $\{B(t)\}$ (see \cite[Chapter~VI,
  Theorem~3.1]{Asmu03}) that for any $\varepsilon>0$,
\[
\lim_{x\to\infty}\inf_{u>(1+\varepsilon)x}
\PP\left( {B(u) - u \over u} > {-\varepsilon \over 1+\varepsilon} \right)
\ge \lim_{x\to\infty}\inf_{u>(1+\varepsilon)x}
\PP\left( \left|{B(u) - u \over u} \right| < {\varepsilon \over 1+\varepsilon} \right) 
 = 1.
\]
Note here that (\ref{property-T-in-C-01}) holds due to $T \in \calC$.
Thus from (\ref{add-eqn-06}), we have $\PP(B(T) > x) \gesimhm{x}
\PP(T>x)$.

In what follows, we prove $\PP(M(T) > x) \lesimhm{x} \PP(T>x)$. 
For any $\varepsilon \in (0,1)$, 
\[
\PP(M(T) > x) 
\le \PP(T > (1-\varepsilon)x) + \PP(M(T) > x, T \le (1-\varepsilon)x).
\]
Since (\ref{property-T-in-C-02}) holds, it suffices to show $\PP(M(T)
> x, T \le (1-\varepsilon)x)=o(\PP(T>x))$.

It follows from $M(u) - u \le \sup_{0\le t \le u}(B(t) - t)$ ($u \ge
0$) that for $x>0$,
\begin{eqnarray*}
\PP(M(T) > x, T \le (1-\varepsilon)x)
&\le& \int_0^{(1-\varepsilon)x}
\PP\left(\sup_{0 \le t \le u}\{B(t) - t\} > x - u\right)\rd \PP(T \le u)
\nonumber
\\
&\le& \int_0^{(1-\varepsilon)x}
\PP\left(\sup_{0 \le t \le u}\{B(t) - t\} > \varepsilon x\right)
\rd \PP(T \le u).
\end{eqnarray*}
Similarly to (\ref{ineqn-prop1-01}), we estimate the integrand on the
right hand side of the above inequality as follows:
\begin{eqnarray*}
\lefteqn{
\PP\left(\sup_{0 \le t \le u}\{B(t) - t\} > \varepsilon x\right)
}
\quad && \nonumber
\\
&\le& \PP\left(\Delta B_0^{\ast} > \varepsilon x / 3\right)
+ \PP\left(\Delta B_1^{\ast} > \varepsilon x / 3\right)
 + \PP\left(\max_{1\le k \le N(u)}\sum_{i=1}^k
 (\Delta B_i - \Delta \tau_i) > {\varepsilon x \over 3}\right).
\end{eqnarray*}
From conditions (i) and (iii), we have $\PP\left(\Delta B_n^{\ast} >
\varepsilon x / 3\right) = o(\PP(T>x))$ ($n=0,1$).  Therefore it
remains to show that
\begin{eqnarray}
\int_0^{(1-\varepsilon)x}
\PP\left(\max_{1\le k \le N(u)}\sum_{i=1}^k
 (\Delta B_i - \Delta \tau_i) > {\varepsilon x \over 3}\right)
\rd \PP(T \le u) = o(\PP(T>x)).
\label{add-eqn-19}
\end{eqnarray}

Fix a positive number $\gamma$ such that
\begin{equation}
{\varepsilon \over 3\gamma} > 
{1-\varepsilon \over \EE[\Delta \tau_1]}.
\label{add-eqn-31}
\end{equation}
We then decompose the left hand side of
(\ref{add-eqn-19}) into $R_1(x) + R_2(x)$ in the following way:
\begin{eqnarray}
R_1(x)
&=& \int_0^{(1-\varepsilon)x}\rd \PP(T \le u)
\PP\left(
\max_{1\le k \le N(u)}\sum_{i=1}^k
 (\Delta B_i - \Delta \tau_i) > {\varepsilon x \over 3},~ 
N(u) > {\varepsilon x \over 3\gamma}
\right),
\nonumber
\\
R_2(x)
&=& \int_0^{(1-\varepsilon)x}\rd \PP(T \le u)
\PP\left(
\max_{1\le k \le N(u)}\sum_{i=1}^k
 (\Delta B_i - \Delta \tau_i) > {\varepsilon x \over 3},~
N(u) \le {\varepsilon x \over 3\gamma}
\right). \qquad
\label{eqn-R_2(x)}
\end{eqnarray}
For $x > 0$, we have
\begin{eqnarray}
R_1(x) 
&\le& \int_0^{(1-\varepsilon)x}
\PP\left(N(u) > {\varepsilon x \over 3\gamma} \right) \rd \PP(T \le u) 
\le \PP\left( N((1-\varepsilon)x) > {\varepsilon x \over 3\gamma} \right).
\label{add-eqn-10}
\end{eqnarray}
Note here that $\varepsilon / (3\gamma) -
(1-\varepsilon) / \EE[\Delta \tau_1] > 0$ due to (\ref{add-eqn-31}). Thus
Lemma~\ref{add-prop-N(x)} yields
\begin{eqnarray*}
\PP\left( N((1-\varepsilon)x) > {\varepsilon x \over 3\gamma} \right)
&=& \PP\left(
N((1-\varepsilon)x) - {(1-\varepsilon) x \over \EE[\Delta \tau_1]}
> \left( {\varepsilon \over 3\gamma} - {1-\varepsilon  \over \EE[\Delta \tau_1]}
 \right)x
\right)
\nonumber
\\
&\le& Ce^{-cx} = o(\PP(T>x)).
\end{eqnarray*}
Combining this with (\ref{add-eqn-10}), we have $R_1(x)= o(\PP(T >
x))$.  

Next we consider $R_2(x)$. From (\ref{eqn-R_2(x)}), we have
\begin{eqnarray*}
R_2(x)
&\le& \int_0^{(1-\varepsilon)x}\rd \PP(T \le u)
\PP\left(
\max_{1\le k \le \varepsilon x/(3\gamma)}\sum_{i=1}^k
 (\Delta B_i - \Delta \tau_i) > {\varepsilon x \over 3}
\right).
\label{eqn-R_2(x)-02}
\end{eqnarray*}
Following the proof of (\ref{2nd-term-A}), we can show that
\[
\PP\left(
\max_{1\le k \le \varepsilon x/(3\gamma)}\sum_{i=1}^k
 (\Delta B_i - \Delta \tau_i) > {\varepsilon x \over 3}
\right)
= o(\PP(T > x)),
\]
which leads to $R_2(x) = o(\PP(T > x))$ (see
subsection~\ref{subsubsec-v.a} and \ref{subsubsec-v.b} in the proof of
Theorem~\ref{thm-dependent-02}).

\begin{rem}\label{rem-R_2(x)}
Except for the estimation of $R_2(x)$, conditions (i)--(iii) and the
independence between $\{B(t)\}$ and $T$ are sufficient for the proof
of Theorem~\ref{thm-independ-02}. Conditions (iv) and (v) are required
by the estimation of $R_2(x)$.
\end{rem}

\subsection{Proof of Lemma~\ref{add-lem-MAdP-01}}\label{proof-add-lem-MAdP-01}

We first partition $\widetilde{\vc{\beta}}$ and
$\widetilde{\vc{H}}$ as
\[
\widetilde{\vc{\beta}}
= 
\bordermatrix{
               & \{0\} &   \bbD\setminus\{0\}       
\cr
 & \widetilde{\beta}_{0} & \widetilde{\vc{\beta}}_{+}
},
\quad
\widetilde{\vc{H}}
= 
\bordermatrix{
               & \{0\} &   \bbD\setminus\{0\}       
\cr
\{0\} & \widetilde{H}_{0,0}& \widetilde{\vc{\eta}}_{+}
\cr
\bbD\setminus\{0\} & \widetilde{\vc{h}}_{+} & \widetilde{\vc{H}}_{+}
}.
\]
We then fix $z=1$ in (\ref{joint-transform-0}) and
(\ref{joint-transform-1}) and take the inverse of them with respect to
$\xi$. Thus
\begin{eqnarray}
\PP(\Delta B_0 \le x)
= \beta_{0}(x) + \vc{\beta}_{+} \ast \sum_{n=0}^{\infty}\vc{H}_+^{\ast n} \ast \vc{h}_{+}(x),
\label{dist-Delta B_0}
\\
\PP(\Delta B_1 \le x)
= H_{0,0}(x) + \vc{\eta}_{+} \ast \sum_{n=0}^{\infty}\vc{H}_+^{\ast n} \ast \vc{h}_{+}(x),
\label{dist-Delta B_1}
\end{eqnarray}
where the symbol $\ast$ denotes the operator of convolution and the
superscript $\ast n$ represents the $n$th-fold convolution (see
Appendix~\ref{subsec-convo-tail}), and where
\[
\vc{\beta}(x)
= 
\bordermatrix{
               & \{0\} &   \bbD\setminus\{0\}       
\cr
 & \beta_{0}(x) & \vc{\beta}_{+}(x)
},
\quad
\vc{H}(x)
= 
\bordermatrix{
               & \{0\} &   \bbD\setminus\{0\}       
\cr
\{0\} & H_{0,0}(x) & \vc{\eta}_{+}(x)
\cr
\bbD\setminus\{0\} & \vc{h}_{+}(x) & \vc{H}_{+}(x)
}.
\]
Applying
Lemma~\ref{lem-convo-tail} to (\ref{dist-Delta B_0}) and
(\ref{dist-Delta B_1}), we obtain
\begin{eqnarray*}
\limsup_{x\to\infty}
{\PP(\Delta B_0 > x) \over \PP(Y>x)}
&\le& \tilde{c}\widetilde{\beta}_{0}
+  \tilde{c}\widetilde{\vc{\beta}}_{+}
(\vc{I} - \vc{H}_{+}(\infty))^{-1}
\vc{h}_{+}(\infty)
\nonumber
\\
&& {} + \vc{\beta}_{+}(\infty)
(\vc{I} - \vc{H}_{+}(\infty))^{-1}
( \tilde{c}\widetilde{\vc{H}}_{+} )
(\vc{I} - \vc{H}_{+}(\infty))^{-1}
\vc{h}_{+}(\infty)
\nonumber
\\
&& {} + \vc{\beta}_{+}(\infty)
(\vc{I} - \vc{H}_{+}(\infty))^{-1}
(\tilde{c}\widetilde{\vc{h}}_{+})
\nonumber
\\
&=& \tilde{c} \left[\widetilde{\vc{\beta}} \vc{e} 
+ \vc{\beta}_{+}(\infty) (\vc{I} - \vc{H}_{+}(\infty))^{-1}
(\widetilde{\vc{H}}_{+} \vc{e} + \widetilde{\vc{h}}_{+}) \right] 
\le \tilde{c} C,
\\
\limsup_{x\to\infty}
{\PP(\Delta B_1 > x) \over \PP(Y>x)}
&\le& \tilde{c}\left[ (\widetilde{H}_{0,0} + \widetilde{\vc{\eta}}_{+} \vc{e} )
+ \vc{\eta}_{+}(\infty) (\vc{I} - \vc{H}_{+}(\infty))^{-1}
(\widetilde{\vc{H}}_{+} \vc{e} + \widetilde{\vc{h}}_{+})
\right]
\nonumber
\\
&=& \tilde{c}(1/\varpi_{0}) 
\left[
\varpi_{0}(\widetilde{H}_{0,0} + \widetilde{\vc{\eta}}_{+} \vc{e} )
+ \vc{\varpi}_{+}(\widetilde{\vc{H}}_{+} \vc{e} + \widetilde{\vc{h}}_{+})
\right]
\nonumber
\\
&=& \tilde{c}(1/\varpi_{0}) \vc{\varpi}\widetilde{\vc{H}} \vc{e} \le \tilde{c} C,
\end{eqnarray*}
where we use $(\vc{I} -
\vc{H}_{+}(\infty))^{-1}\vc{h}_{+}(\infty)=\vc{e}$ (which is due to
$\vc{h}_{+}(\infty) + \vc{H}_{+}(\infty)\vc{e} = \vc{e}$); and also
use $\vc{\varpi}_{+} :=(\varpi_i)_{i\in\bbD\setminus\{0\}} =
\varpi_{0}\vc{\eta}_{+}(\infty)(\vc{I} - \vc{H}_{+}(\infty))^{-1}$ and
$\varpi_{0} = 1/\EE[\Delta \tau_1]$ (see, e.g., \cite[Chapter~3,
  Theorems~2.1 and 3.2]{Brem99}).

\subsection{Proof of Lemma~\ref{add-lem-03}}\label{proof-add-lem-03}

Let $\psi_1(k,\xi)$ ($k=1,2,\dots$) denote
\[
\psi_1(k,\xi) 
= \EE[\dd{1}(\Delta \tau_1 = k)e^{\rmi \xi \Delta B_1}]
= {1 \over k!}\lim_{z\to0}{\partial^k \over \partial z^k}
\widehat{\psi}_1(z,\xi).
\]
It then follows from (\ref{joint-transform-1}) that for $k=1,2,\dots$,
\begin{eqnarray}
\psi_1(k,\xi) 
&=& \dd{1}(k=1) 
\widehat{H}_{0,0}(\xi)
+ \dd{1}(k\ge2)\widehat{\vc{\eta}}_{+}(\xi) \cdot
(\widehat{\vc{H}}_{+}(\xi))^{k-2}
\cdot
\widehat{\vc{h}}_{+}(\xi).
\label{add-eqn-43}
\end{eqnarray}
Taking the inverse of (\ref{add-eqn-43}) with respect to $\xi$ and
applying Lemma~\ref{lem-convo-tail} to the resulting equation, we
obtain
\begin{eqnarray}
\lefteqn{
\limsup_{x\to\infty}
{\PP(\Delta \tau_1 = k, \Delta B_1 > x) \over \PP(Y>x)}
}
\quad &&
\nonumber
\\
&\le& \tilde{c} \dd{1}(k=1) 
\widetilde{H}_{0,0}
+ \tilde{c} \dd{1}(k\ge2)
\left\{
\widetilde{\vc{\eta}}_{+} 
(\vc{H}_{+}(\infty))^{k-2}
\vc{h}_{+}(\infty)
\right\}
\nonumber
\\
&&{}
+ \tilde{c} \dd{1}(k\ge2)
\left\{
\vc{\eta}_{+}(\infty) \cdot
\sum_{\nu=0}^{k-3}
(\vc{H}_{+}(\infty))^{\nu}
\widetilde{\vc{H}}_{+}
(\vc{H}_{+}(\infty))^{k-\nu-3}
\cdot
\vc{h}_{+}(\infty)
\right\}
\nonumber
\\
&&{}
+ \tilde{c} \dd{1}(k\ge2)
\left\{
\vc{\eta}_{+}(\infty) 
(\vc{H}_{+}(\infty))^{k-2}
\widetilde{\vc{h}}_{+}
\right\}, \qquad \forall k=1,2,\dots.
\label{add-eqn-44}
\end{eqnarray}
Note here that $\beta_i(\infty) = 0$ (resp.\ $H_{i,j}(\infty) = 0$)
implies $\widetilde{\beta}_i = 0$ (resp.\ $\widetilde{H}_{i,j} = 0$)
and thus $\widetilde{\vc{\beta}} \le C\vc{\beta}(\infty) =
C\widehat{\vc{\beta}}(0)$ and $\widetilde{\vc{H}} \le
C\vc{H}(\infty)=C\widehat{\vc{H}}(0)$. Therefore from
(\ref{add-eqn-44}) and (\ref{dist-Delta-tau_1}), we have for all
$k=1,2,\dots$,
\begin{eqnarray*}
\lefteqn{
\limsup_{x\to\infty}
{
\PP(\Delta \tau_1 = k, \Delta B_1 > x) \over \PP(Y>x)}
}
~~ &&
\nonumber
\\
&\le& \tilde{c} C
\left[ \dd{1}(k=1) k\widehat{H}_{0,0}(0)
+ \dd{1}(k\ge2)k
\left\{
\widehat{\vc{\eta}}_{+}(0)
\left(\widehat{\vc{H}}_{+}(0)\right)^{k-2}
\widehat{\vc{h}}_{+}(0)
\right\}
\right]
= \tilde{c} C k \PP(\Delta \tau_1 = k),
\end{eqnarray*}
where $C$ is independent of $k$.

\subsection{Proof of Lemma~\ref{lem-MAdP-03}}\label{proof-lem-MAdP-03}

Since $\PP(\sum_{i=1}^0 \Delta B_i > 0 \mid N(t) = 0) =
\PP(\{\emptyset\})=0$, (\ref{add-eqn-29}) holds for all $t \ge 0$ if
$m=0$.

In what follows, we consider the case of $m\ge 1$.  Under
Assumption~\ref{assumpt-(S_n,J_n)}, $\Delta \tau_1 \ge 1$ and $N(t) =
N(\lfloor t \rfloor) \le \lfloor t \rfloor$ for all $t \ge
0$. Therefore we fix $t=n \in \{1,2,\dots\}$ without loss of
generality.

Note that $\{N(n) = m\}$ is equivalent to $\{ \sum_{i=1}^m\Delta
\tau_i \le n, \sum_{i=1}^{m+1}\Delta \tau_i > n \}$ and that $\Delta
\tau_{m+1}$ is independent of $\Delta \tau_i$ and $\Delta B_i$
($i=1,2,\dots,m$). We then have
\begin{eqnarray}
\PP\left( N(n) = m, \sum_{i=1}^m\Delta B_i > x\right)
&=& 
\PP\left(\sum_{i=1}^m\Delta \tau_i \le n, \sum_{i=1}^{m+1}\Delta \tau_i > n, 
\sum_{i=1}^m\Delta B_i > x\right)
\nonumber
\\
&=& \sum_{k=1}^n
\PP\left(\sum_{i=1}^m\Delta \tau_i=k, 
\sum_{i=1}^m\Delta B_i > x\right)
\nonumber
\\
&& {}\qquad  \times 
\PP(\Delta \tau_{m+1} > n-k).
\label{add-eqn-41}
\end{eqnarray}
Note also that $\Delta B_i$ is independent of $\Delta \tau_j$'s
($j\neq i$).  We thus have
\begin{eqnarray}
\lefteqn{
\PP\left(\sum_{i=1}^m\Delta \tau_i=k, 
\sum_{i=1}^m\Delta B_i > x\right)
}
\quad &&
\nonumber
\\
&=& \sum_{k_1+\cdots+k_m=k} \prod_{i=1}^m\PP(\Delta \tau_i=k_i)
 \cdot \PP\left(\left.
\sum_{i=1}^m\Delta B_i > x 
\,\right| \Delta \tau_i=k_i,~i=1,2,\dots,m\right)
\nonumber
\\
&=& \sum_{k_1+\cdots+k_m=k} \prod_{i=1}^m\PP(\Delta \tau_i=k_i)
 \cdot \PP\left(
\sum_{i=1}^m ( \Delta B_i\,|\,\{\Delta \tau_i=k_i\} ) > x \right), \quad
\label{add-eqn-47}
\end{eqnarray}
where $\Delta B_i\,|\,\{\Delta \tau_i=k_i\}$ denotes the conditional
random variable $\Delta B_i$ given $\Delta \tau_i=k_i$. Further it
follows from Lemmas~\ref{add-lem-03} and \ref{lem-convo-tail} that for
$(k_1,\dots,k_m)$ such that $\sum_{i=1}^m k_i = k$,
\begin{eqnarray*}
\limsup_{x\to\infty}
{\PP\left(
\sum_{i=1}^m ( \Delta B_i\,|\,\{\Delta \tau_i=k_i\} ) > x \right) 
\over \PP(Y>x)}
&\le& \tilde{c} C \cdot (k_1+\cdots+k_m)
= \tilde{c} C k. 
\end{eqnarray*}
Combining this with (\ref{add-eqn-47}) yields
\begin{eqnarray}
\lefteqn{
\limsup_{x\to\infty}
{
\PP\left(\sum_{i=1}^m\Delta \tau_i=k, 
\sum_{i=1}^m\Delta B_i > x\right)
\over \PP(Y>x)}
}
\quad &&
\nonumber
\\
&\le& 
\tilde{c} C k\sum_{k_1+\cdots+k_m=k} \prod_{i=1}^m\PP(\Delta \tau_i=k_i)
=
\tilde{c} C k \cdot \PP\left(\sum_{i=1}^m\Delta \tau_i=k\right).
\label{add-eqn-49}
\end{eqnarray}
From (\ref{add-eqn-41}) and (\ref{add-eqn-49}), we obtain for all
$n=0,1,\dots$ and $m=0,1,\dots,n$,
\begin{eqnarray*}
\limsup_{x\to\infty}
{
\PP\left( N(n) = m, \sum_{i=1}^m\Delta B_i > x\right)
\over \PP(Y>x)}
&\le& \tilde{c} C \sum_{k=1}^n k
\PP\left(\sum_{i=1}^m\Delta \tau_i=k\right)
\PP(\Delta \tau_{m+1} > n-k)
\nonumber
\\
&\le& \tilde{c} C n\sum_{k=1}^n 
\PP\left(\sum_{i=1}^m\Delta \tau_i=k\right)
\PP(\Delta \tau_{m+1} > n-k)
\nonumber
\\
&=& \tilde{c} C n
\PP\left(\sum_{i=1}^m\Delta \tau_i \le n, \sum_{i=1}^{m+1}\Delta \tau_i > n 
\right)
\nonumber
\\
&=& \tilde{c} C n \PP(N(n) = m).
\end{eqnarray*}

\subsection{Proof of Theorem~\ref{thm-MAdP-01}}\label{proof-thm-MAdP-01}

As shown later, the conditions of Theorem~\ref{thm-MAdP-01} imply
conditions (i), (ii) and (iii) of Theorem~\ref{thm-independ-02}. Thus
according to Remark~\ref{rem-R_2(x)}, we can follow the proof of
Theorem~\ref{thm-independ-02}, except for the estimation of $R_2(x)$
in (\ref{eqn-R_2(x)}). In addition, we can prove that $R_2(x) =
o(\PP(T > x))$ as follows.

From (\ref{eqn-R_2(x)}), we have
\begin{eqnarray}
R_2(x) 
&\le& \int_0^{(1-\varepsilon)x}\sum_{n\le \varepsilon x/(3\gamma)}
\PP(N(u) = n)
\PP\left(\left. \sum_{i=1}^n \Delta B_i > {\varepsilon x \over 3} 
\, \right| N(u) = n\right)
\rd \PP(T \le u). \qquad
\label{add-3rd-07}
\end{eqnarray}
Note here that condition (iii) of Theorem~\ref{thm-independ-02}
implies $\PP(\Delta B_1 > x) = o(\PP(T>x))$. Thus by using
Lemma~\ref{lem-MAdP-03} with $Y=T \in \calC$ and $\tilde{c}=0$, we
obtain
\[
\PP\left(\left. \sum_{i=1}^n \Delta B_i > {\varepsilon x \over 3} 
\, \right| N(u) = n\right)
= u \cdot o(\PP(T>x)).
\]
Substituting this into (\ref{add-3rd-07}) yields
\begin{eqnarray*}
R_2(x) 
&\le& \int_0^{(1-\varepsilon)x}\sum_{n\le \varepsilon x/(3\gamma)}
\PP(N(u) = n) u
\rd \PP(T \le u) \cdot o(\PP(T>x))
\nonumber
\\
&\le& \EE[T] \cdot o(\PP(T>x)),
\end{eqnarray*}
which implies that $R_2(x) = o(\PP(T > x))$ due to $\EE[T] < \infty$.

In what follows, we confirm that conditions (ii) and (iii) of
Theorem~\ref{thm-independ-02} are satisfied (condition (i) is
obvious). For simplicity, we assume $h = b = 1$, which does not lose
generality.

We first introduce a cumulative process $\{B^{\#}(t);t\ge0\}$ such
that $B^{\#}(t) = \sum_{n=0}^{\lfloor t \rfloor}|X_n|$ for $t \ge 0$.
Clearly, $\{B^{\#}(t)\}$ and $\{B(t)\}$ have the common regenerative
points $\tau_n$'s. Further $\{(B^{\#}(n),J_n);n=0,1,\dots\}$ is a
Markov additive process with initial distribution $\vc{\beta}^{\#}(x)$
and kernel $\vc{H}^{\#}(x)$ ($x \in \bbR$), where $\vc{\beta}^{\#}(x)
= \int_{|y| \le x}\rd \vc{\beta}(y)$ and $\vc{H}^{\#}(x) = \int_{|y|
  \le x}\rd \vc{H}(y)$.

Let $\Delta B^{\#}_n$ ($n=0,1,\dots$) denote
\[
\Delta B^{\#}_n
= 
\left\{
\begin{array}{ll}
B^{\#}(\tau_0), & n=0,
\\
B^{\#}(\tau_n) - B^{\#}(\tau_{n-1}), & n=1,2,\dots.
\end{array}
\right.
\]
We then have
\begin{eqnarray*}
\Delta B^{\#}_0
&\ge& \sup_{0 \le t \le \tau_0} |B(t)|
\ge \Delta B_0^{\ast} \ge \Delta B_0,
\\
\Delta B^{\#}_n
&\ge& \sup_{\tau_{n-1} \le t \le \tau_n} |B(t) - B(\tau_{n-1})| 
\ge \Delta B_n^{\ast} \ge \Delta B_n,\quad n=1,2,\dots.
\end{eqnarray*}
Thus, similarly to the proof of Proposition~\ref{prop-b=h}, we readily
obtain
\[
\EE\left[ \sup_{\tau_0 \le t \le \tau_1} |B(t) - B(\tau_0)| \right]
\le \EE[ \Delta B_1^{\#}]
= \vc{\varpi}\int_{x\in\bbR} |x| \rd \vc{H}(x)\vc{e}
\cdot \EE[\Delta\tau_1] < \infty,
\]
where the last inequality is due to
Assumption~\ref{assumpt-(S_n,J_n)}~(iii).  Recall here that
$\Delta\tau_n$ follows a phase-type distribution and thus
$\EE[(\Delta\tau_n)^2] < \infty$ ($n=0,1$). Therefore condition (ii)
of Theorem~\ref{thm-independ-02} is satisfied. Further following the
proof of Lemma~\ref{add-lem-MAdP-01} with $Y=T$ and $\tilde{c} = 0$,
we can prove that
\[
\PP( \Delta B_n^{\ast} > x)
\le \PP( \Delta B_n^{\#} > x) = o(\PP(T > x)),\qquad n=0,1,
\]
which shows that condition (iii) of Theorem~\ref{thm-independ-02} is
satisfied. As a result, the conditions of Theorem~\ref{thm-MAdP-01}
imply conditions (i), (ii) and (iii) of Theorem~\ref{thm-independ-02}.

\subsection{Proof of Theorem~\ref{thm-MAdP-02}}\label{proof-thm-MAdP-02}

Note that (\ref{MAdP-assumpt-02}) and (\ref{MAdP-assumpt-03}) yield
(\ref{MAdP-assumpt-01}) and thus the conditions of
Theorem~\ref{thm-MAdP-02} imply those of Theorem~\ref{thm-MAdP-01},
except for $\EE[T] < \infty$. Note also that $\EE[T] < \infty$ is not
covered by conditions (i), (ii) and (iii) of
Theorem~\ref{thm-independ-02}. Therefore the conditions of
Theorem~\ref{thm-MAdP-02} imply conditions (i), (ii) and (iii) of
Theorem~\ref{thm-independ-02} (see the proof of
Theorem~\ref{thm-MAdP-01} in subsection~\ref{proof-thm-MAdP-01}). As a
result, it suffices to prove $R_2(x) = o(\PP(T>x))$ (see
Remark~\ref{rem-R_2(x)}).

It follows from (\ref{add-3rd-07}) and Lemma~\ref{lem-MAdP-03} that
\begin{eqnarray}
R_2(x)
&\le& C\int_0^{(1-\varepsilon)x} u \rd \PP(T \le u)
\sum_{n\le \varepsilon x/(3\gamma)}
\PP(N(u) = n) \cdot \PP(Y > x),
\label{add-3rd-06}
\end{eqnarray}
where $\gamma$ is a positive number satisfying (\ref{add-eqn-31}). We
now fix $\gamma$ to be
\[
{1-\varepsilon \over \EE[\Delta \tau_1]} 
< {\varepsilon \over 3\gamma}
\le 
{1 \over \EE[\Delta \tau_1]}.
\]
As a result, from (\ref{add-3rd-06}), we have
\begin{eqnarray*}
R_2(x) 
&\le& C\int_0^{(1-\varepsilon)x} u \rd \PP(T \le u)
\sum_{n\le x/\EE[\Delta \tau_1]}
\PP(N(u) = n) \cdot \PP(Y > x)
\nonumber
\\
&=& C \EE[T \cdot \dd{1}(T\le x, N(T) \le x/\EE[\Delta \tau_1])]
\cdot \PP(Y > x)
= o(\PP(T > x)),
\end{eqnarray*}
where the last equality is due to (\ref{MAdP-assumpt-03}).


\section*{Acknowledgments}
The author thanks Professor Naoto Miyoshi for his helpful comments on
related work. The author also thanks anonymous referees for their
invaluable comments and suggestions on improving the presentation of
this paper.  This research was supported in part by Grant-in-Aid for
Young Scientists (B) of Japan Society for the Promotion of Science
under Grant No.~24710165.


\begin{thebibliography}{25}
%
\bibitem{Ales08}
A.~Ale\v{s}kevi\v{c}ien\.{e}, R.~Leipus, and J.~\v{S}iaulys:
Tail behavior of random sums under consistent variation with applications 
to the compound renewal risk model. {\it Extremes}, {\bf 11} (2008), 261--279.

\bibitem{Alfa95}
A.S.~Alfa and M.F.~Neuts: 
Modelling vehicular traffic using the discrete time Markovian arrival process. 
{\it Transportation Science}, {\bf 29} (1995), 109--117.

\bibitem{Asmu03}
S.~Asmussen: \emph{Applied Probability and Queues, 2nd~ed.}
(Springer, New York, 2003).

\bibitem{AsmuKlupSigm99}
S.~Asmussen, C.~Kl\"{u}ppelberg, and K.~Sigman:
Sampling at subexponential times, with queueing applications. 
{\it Stochastic Processes and their Applications}, {\bf 79} (1999), 265--286.

\bibitem{Baye96} 
N.~Bayer and O.J.~Boxma: 
Wiener-Hopf analysis of an M/G/1 queue with negative customers and of
a related class of random walks.
{\it Queueing Systems}, {\bf 23} (1996), 301--316.

\bibitem{Bing89} 
N.H.~Bingham, C.M.~Goldie, and J.L.~Teugels: 
{\it Regular Variation} (Cambridge University Press, Cambridge, UK, 1989).

\bibitem{Brem99}
P.~Br\'emaud:
{\it Markov Chains: Gibbs Fields, Monte Carlo Simulation, and Queues}
(Springer, New York, 1999). 

\bibitem{Breu05}
L.~Breuer and A.S.~Alfa:
An EM algorithm for platoon arrival processes in discrete time.
{\it Operations Research Letters}, {\bf 33} (2005), 535--543.

\bibitem{Cart00}
M.~Carter and B.~van Brunt:
{\it The Lebesgue-Stieltjes Integral} (Springer, New York, 2000).

\bibitem{Chis64} 
V.P.~Chistyakov:
A theorem on sums of independent positive random variables and 
its applications to branching random processes. 
{\it Theory of Probability and Its Applications}, {\bf 9} (1964), 640--648.

\bibitem{Clin94}
D.B.H.~Cline:
Intermediate regular and {$\Pi$} variation. 
{\it Proceedings of the London Mathematical Society}, {\bf 68} (1994), 594--616.

\bibitem{Embr97}
P.~Embrechts, C.~Kl\"{u}ppelberg, and T.~Mikosch:
{\it Modelling Extremal Events for Insurance and Finance}
(Springer, Berlin, 1997).

\bibitem{Embr84}
P.~Embrechts and E.~Omey:
A property of longtailed distributions. 
{\it Journal of Applied Probability}, {\bf 21} (1984), 80--87.

\bibitem{Fay06}
G.~Fa\"{y}, B.~Gonz\'{a}lez-Ar\'{e}valo, T.~Mikosch, and G.~ Samorodnitsky:
Modeling teletraffic arrivals by a Poisson cluster process. 
{\it Queueing Systems}, {\bf 54} (2006), 121--140.

\bibitem{Foss00}
S.~Foss and D.~Korshunov:
Sampling at a random time with a heavy-tailed distribution. 
{\it Markov Processes and Related Fields} {\bf 6} (2000), 543--568.

\bibitem{Foss11}
S.~Foss, D.~Korshunov, and S.~Zachary:
{\it An Introduction to Heavy-Tailed and Subexponential Distributions}
(Springer, New York, 2011).

\bibitem{Galm01}
S.~Galm\'es and R.~Puigjaner:
Performance evaluation based on an aggregate ATM model.
In {\it Proceedings of the 9th IEEE International Symposium on Modeling,
  Analysis and Simulation of Computer and Telecommunication Systems} 
(Cincinnati, OH, 2001), 399--406.

\bibitem{Galm03}
S.~Galm\'es and R.~Puigjaner:
An algorithm for computing the mean response time of a
single server queue with generalized on/off traffic arrivals.
{\it ACM SIGMETRICS Performance Evaluation Review}, {\bf 31} (2003), 37--46.
%
\bibitem{Galm05}
S.~Galm\'es and R.~Puigjaner:
The response time distribution of a discrete-time queue
under a generalized batch arrival process.
In {\it Proceedings of the 3rd International IFIP/ACM Latin American
  Conference on Networking} (Cali, Colombia, 2005),
31--39.
%
\bibitem{Gold98} C.M.~Goldie and C.~Kl\"{u}ppelberg: 
Subexponential distributions. 
In R.J.~Adler, R.E.~Feldman, and M.S.~Taqqu (eds.):  
{\it A Practical Guide to Heavy Tails: Statistical Techniques and
Applications} (Birkh\"{a}user, Boston, 1998), 435--459.

\bibitem{Jele99}
P.R.~Jelenkovi\'{c}:
Subexponential loss rates in a GI/GI/1 queue with applications. 
{\it Queueing Systems}, {\bf 33} (1999), 91--123.

\bibitem{Jele98} 
P.R.~Jelenkovi\'{c} and A.A.~Lazar:
Subexponential asymptotics of a Markov-modulated random walk with
queueing applications.
{\it Journal of Applied Probability}, {\bf 35} (1998), 325--347.

\bibitem{Jele03a}
P.R.~Jelenkovi\'{c} and P.~Mom\v{c}ilovi\'{c}:
Large deviation analysis of subexponential waiting times in a
processor-sharing queue. 
{\it Mathematics of Operations Research}, {\bf 28} (2003), 587--608.

\bibitem{Jele03b}
P.R.~Jelenkovi\'{c} and P.~Mom\v{c}ilovi\'{c}:
Asymptotic loss probability in a finite buffer
fluid queue with heterogeneous heavy-tailed on-off processes. 
{\it The Annals of Applied Probability}, {\bf 13} (2003), 576--603.

\bibitem{Jele04} 
P.R.~Jelenkovi\'{c}, P.~Mom\v{c}ilovi\'{c}, and B.~Zwart:
Reduced load equivalence under subexponentiality. 
{\it Queueing Systems}, {\bf 46} (2004), 97--112.

\bibitem{Klup88} 
C.~Kl\"{u}ppelberg: 
Subexponential distributions and integrated tails. 
{\it Journal of Applied Probability}, {\bf 25} (1988), 132--141.

\bibitem{Kors02}
D.A.~Korshunov:
Large-deviation probabilities for maxima of sums of
  independent random variables with negative mean and subexponential distribution. 
{\it Theory of Probability and Its Applications}, {\bf 46} (2002), 355--365.

\bibitem{Lato99} 
G.~Latouche and V.~Ramaswami:
{\it Introduction to Matrix Analytic Methods in Stochastic Modeling}
(ASA--SIAM, Philadelphia, PA, 1999).

\bibitem{Lin11}
Z.~Lin and X.~Shen:
Approximation of the tail probability of dependent random sums under
consistent variation and applications. 
{\it Methodology and Computing in Applied Probability}, {\bf 15} (2013), 165--186.

\bibitem{Luca91} 
D.M.~Lucantoni: 
New results on the single server queue with a batch Markovian arrival
process.  {\it Stochastic Models}, {\bf 7} (1991), 1--46.

\bibitem{Masu11}
H.~Masuyama:
Subexponential asymptotics of the stationary distributions of
M/G/1-type Markov chains. 
{\it European Journal of Operational Research}, 
{\bf 213} (2011), 509--516.

\bibitem{Masu13a} 
H.~Masuyama: 
A sufficient condition for subexponential asymptotics of GI/G/1-type
Markov chains and its application to BMAP/GI/1 queues. Preprint arXiv:1310.4590, 2013 (available online at http://arxiv.org/abs/1310.4590).

\bibitem{Masu13b}
H.~Masuyama:
Subexponential tail equivalence of the queue length distributions of
BMAP/GI/1 queues with and without retrials. Preprint arXiv:1310.4608, 2013 (available online at http://arxiv.org/abs/1310.4608).

\bibitem{Masu09}
H.~Masuyama, B.~Liu, and T.~Takine:
Subexponential asymptotics of the {BMAP/GI/1} queue. 
{\it Journal of the Operations Research Society of Japan}, 
{\bf 52} (2009), 377--401.

\bibitem{Miyo11}
N.~Miyoshi, M.~Ogura, and S.~Maruyama:
Long-tailed degree distribution of a random geometric graph
constructed by the Boolean model with spherical grains.
Research Report \#B--464, Department of
Mathematical and Computing Sciences, Tokyo Institute of Technology {\bf 2011}.

\bibitem{Mull02}
A.~M\"uller and D.~Stoyan: 
{\it Comparison Methods for Stochastic Models and Risks}
(John Wiley \& Sons, Chichester, 2002).

\bibitem{Naga77} 
A.V.~Nagaev:
On a property of sums of independent random variables.
{\it Theory of Probability and Its Applications}, {\bf 22}
(1977), 326--338.

\bibitem{Pitm80}
E.J.G.~Pitman:
Subexponential distribution functions.
{\it Journal of the Australian Mathematical Society}, {\bf A29} 
(1980), 337--347.

\bibitem{Robe08}
C.Y.~Robert and J.~Segers:
Tails of random sums of a heavy-tailed number of light-tailed terms. 
{\it Insurance: Mathematics and Economics}, {\bf 43} (2008), 85--92.

\bibitem{Shne04}
V.V.~Shneer:
Estimates for the distributions of the sums of subexponential random variables. 
{\it Siberian Mathematical Journal}, {\bf 45} (2004), 1143--1158.

\bibitem{Shne06}
V.V.~Shneer:
Estimates for interval probabilities of the sums of random variables
with locally subexponential distributions.  
{\it Siberian Mathematical Journal}, {\bf 47} (2006), 779--786.

\bibitem{Sigm99} 
K.~Sigman: 
Appendix: A primer on heavy-tailed distributions. 
{\it Queueing Systems}, {\bf 33} (1999), 261--275.

\bibitem{Taki04} 
T.~Takine: 
Geometric and subexponential asymptotics of Markov chains of M/G/1
type.
{\it Mathematics of Operations Research}, {\bf 29} (2004), 624--648.

\bibitem{Tang06}
Q.~Tang:
Insensitivity to negative dependence of the asymptotic behavior of
precise large deviations. 
{\it Electronic Journal of Probability}, {\bf 11} (2006), 107--120.

\bibitem{Will91}
D.~Williams:
{\it Probability with Martingales}
(Cambridge University Press, Cambridge, UK, 1991).

\bibitem{Wolf89}
R.W.~Wolff:
{\it Stochastic Modeling and the Theory of Queues}
(Prentice-Hall, Englewood Cliffs, NJ, 1989). 

\bibitem{Zwar00}
A.P.~Zwart:
A fluid queue with a finite buffer and subexponential input.
{\it Advances in Applied Probability}, {\bf 32} (2000), 221--243.
\end{thebibliography}

\end{document}